\documentclass[12pt]{article}
\usepackage{smile}
\usepackage{mathtools}
\usepackage{booktabs}
\usepackage[english]{babel} % English language/hyphenation
\usepackage[protrusion=true,expansion=true]{microtype} % Better typography
\usepackage{amsmath,amsfonts,amsthm}
\usepackage{ amssymb }
\usepackage{enumitem}
\usepackage{graphicx}
\usepackage{array}
\usepackage[table]{xcolor}
\usepackage[toc,page]{appendix}
\usepackage{booktabs} % Horizontal rules in tables
\usepackage{natbib}

\usepackage[margin=1in,hmarginratio=1:1,top=20mm,columnsep=20pt]{geometry} % Document margins
\usepackage{microtype} % Slightly tweak font spacing for aesthetics
\usepackage[font = small,labelfont=bf,textfont=it]{caption} % Custom captions under/above floats in tables or figures
\usepackage{footnote}
\usepackage{multirow}
\usepackage{algorithm}
\usepackage{algpseudocode}

\usepackage{color}

\title{Sobolev Calibration of Imperfect Computer Models}
\author{Qingwen Zhang\\
  Division of Emerging Interdisciplinary Areas,\\
    The Hong Kong University of Science and Technology\\
    Wenjia Wang \thanks{Corresponding author: wenjiawang@ust.hk}\\
    Data Science and Analytics Thrust, \\The Hong Kong University of Science and Technology (Guangzhou)\\
    Department of Mathematics, \\The Hong Kong University of Science and Technology}
\date{}

\begin{document}
\maketitle
% \vspace{-5mm}

\begin{abstract}
Calibration refers to the statistical estimation of unknown model parameters in computer experiments, 
such that computer experiments can match underlying physical systems. 
This work develops a new calibration method for imperfect computer models, Sobolev calibration,  which can rule out calibration parameters that generate overfitting calibrated functions.
We prove that the Sobolev calibration enjoys desired theoretical properties including fast convergence rate, asymptotic normality and semiparametric efficiency. We also demonstrate an interesting property that the Sobolev calibration can bridge the gap between two {influential} methods: $L_2$ calibration and Kennedy and O'Hagan's calibration. In addition to exploring the deterministic physical experiments, we theoretically justify that our method can transfer to the case when the physical process is indeed a Gaussian process, which follows the original idea of Kennedy and O'Hagan's. Numerical simulations as well as a real-world example illustrate the competitive performance of the proposed method.\end{abstract}

\section{Introduction}

Computer experiments are widely adopted by scientists and engineers as a simplified mathematical representation of complex physical processes \citep{fang2005design}. The applications of computer experiments have been popularized to various fields, such as those in hydrology \citep{white2005sensitivity}, biology \citep{cirovic2006computer}, and weather prediction \citep{lynch2008origins}. The mathematical models of computer experiments need to take \textit{calibration parameter} $\btheta$ as an input, which often represents inherent attributes of physical environment, and is difficult to measure with tools in reality \citep{ko,plumlee2017bayesian}. The act of estimating the calibration parameters in the sense that matching physical observations to the largest extent is termed as \textit{calibration}, which has been developed as a well-established area of statistics.

However, even with the best-tuned calibration parameters, the computer models may not align perfectly with physical experiments \citep{ko}. This is mainly because the mathematical forms of computer models are usually based on simplified assumptions \citep{tuo&wu2014,plumlee2017bayesian}. Following \cite{tuo&wu2014}, we name the computer model with inadequacy issue as \textit{imperfect}. 

\cite{ko} propose the first study to tackle model inadequacy issue by introducing the discrepancy function with the Gaussian process prior into a Bayesian calibration framework, known as Kennedy-O'Hagan calibration method (abbreviated as KO calibration). The landmark research inspired a lot of statisticians to explore under the {Bayesian} calibration framework  \citep{higdon2004combining,bayarri2007framework,qian2008bayesian,wang2009bayesian}. {Despite many successful applications under Bayesian framework, the identification issue remains a fundamentally weak spot, and may sabotage the empirical performance \citep{gu2018scaled}. Another branch of study is the frequentist framework of calibration, which tackles the problem by directly defining an identifiable calibration parameter in a discrepancy-minimal way. \cite{tuo&wu2,tuo&wu2014} propose $L_2$ calibration, which resorts to minimizing $L_2$ distance between physical process and computer model as the identifiable definition of calibration parameter $\theta^*_{L_2}$. They prove that the estimate $\hat{\theta}_{L_2}$ enjoys nice asymptotic properties including $\sqrt{n}$-consistency and semi-parametric efficiency under specific regularity conditions. \cite{wong2017frequentist} propose a similar OLS-type calibration method under fixed design, which also enjoys fast rates of convergence.} In addition to KO calibration and $L_2$ calibration, 
many efforts are made to improve or generalize these two calibration methods \citep{plumlee2017bayesian,gu2018scaled,plumlee2019computer,xie2020bayesian,plumlee2016calibrating,sung2020calibration,farah2014bayesian,gramacy2015calibrating,joseph2015engineering,storlie2015calibration}. 

Although KO calibration is a classical Bayesian method, its intrinsic nature can be examined by investigating its frequentist properties. \cite{tuo&wu2,tuo2020improved} show that the KO calibration has a simplified frequentist version (they call it \textit{frequentist KO calibration}), where the calibration parameter minimizes the distance between the physical and computer experiments, measured by the the norm in a Reproducing kernel Hilbert space (RKHS). The kernel is found identical to the one used in the kernel interpolator of discrepancy function. If the RKHS can be embedded into a Sobolev space with some smoothness $m$, the frequentist KO calibration actually emulates the system with a closer representation for the first and higher order of derivatives, which encourages function shape approximation. Compared to $L_2$ calibration which ignores the function shape, this is a good property because building calibrated computer experiments that can approximate the shape of physical process well is indeed a practical need. As an illustration, in experiments involving self-driving cars, to ensure alignment with the physical driving track, the simulation system should strive to approximate both velocity (first-order derivative) and acceleration (second-order derivative). In addition to calibration problem, the idea that connects parameter estimation to smoothing techniques also holds crucial significance in fields like solving differential equations, as demonstrated in \cite{ramsay2007parameter}.

\begin{figure}[!h]
\begin{center}
    \includegraphics[width=3.5in]{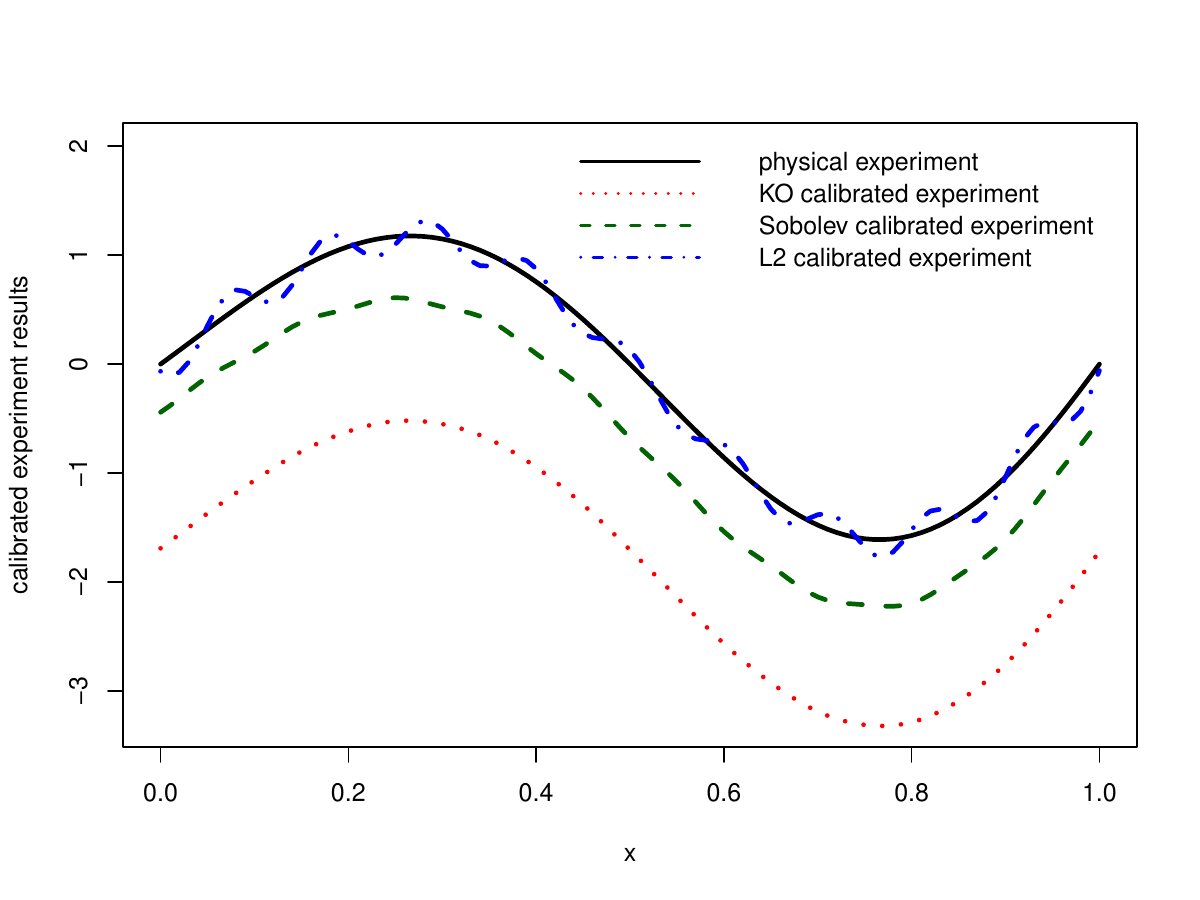}
\end{center}
    \caption{Visualization of the calibrated computer models in Example 1 of numerical experiments. The Sobolev calibration, the KO calibration, and the $L_2$ calibration are presented in dotted lines, while the physical experiment is presented in a solid line.}
    \label{intro:eg2}
\end{figure}

Despite the intriguing frequentist property of KO calibration, \cite{tuo&wu2} argues that the frequentist KO calibration can yield results that are far from minimal $L_2$ distance. This is due to the fact that, once the kernel interpolator of the discrepancy function is decided by the user (where the smoothness can be large), the chosen over-smooth RKHS forces frequentist KO method to mimic the shape even at a cost of fairly large point-wise magnitude difference, as shown in Figure \ref{intro:eg2}. In this work, we propose a novel statistical procedure of estimating calibration parameters, called Sobolev calibration, which can break the limitation of aforementioned two methods. Practitioners can select the calibration parameter by their own preference on the trade-off between magnitude approximation and shape approximation. Our proposed method can realize a more ideal calibrated experiment, as illustrated in Figure \ref{intro:eg2} with green line, which mirrors the shape of physical process with a tunable function value shift. We rigorously provide the theoretical guarantee of Sobolev calibration. Theoretical analysis shows that the Sobolev calibration not only enjoys the desirable properties as $L_2$ calibration, including fast convergence rate, asymptotic normality, and semiparametric efficiency, but also bridges the gap between $L_2$ calibration and frequentist KO calibration, in the sense that $L_2$ calibration and frequentist KO calibration are special cases of the proposed Sobolev calibration. 

In the aforementioned theoretical works of frequentist calibration methods, both the physical and computer experiments are assumed to be deterministic functions. In this work, to align with the original idea of KO calibration, we additionally study the scenario when the physical function is a random function following a Gaussian process and demonstrate the theoretical properties. Since the support of a Gaussian process is typically larger than the corresponding RKHS \citep{van2008reproducing}, the theoretical development of Gaussian process based calibration is much different with that of the deterministic function based calibration; thus we believe it has its own interest. 

The remainder of this paper is organized as follows. In Section \ref{sec:prblmst}, we formulate the Sobolev calibration method, and develop an efficient computing method for the optimization problem. In Section \ref{sec:theory_deter}, we study asymptotic behavior of Sobolev calibration under deterministic case, discuss the uncertainty quantification and build a connection between Sobolev calibration and other two important calibration methods: $L_2$ calibration and frequentist KO calibration. In Section \ref{sec:theory_gp}, we explore the Sobolev calibration under the case when the underlying physical process admits Gaussian process, and show the asymptotic results. In Section \ref{sec:eg}, numerical examples and a real-world example are used to illustrate the performance of our proposed method. Concluding remarks are left in Section \ref{sec:conclusion}. Details of the proof, additional experiment results and discussion can be found in the supplement.

\section{Sobolev Calibration}\label{sec:prblmst}

\subsection{Background and Problem Settings}\label{subsec:background}

Let $\Omega\subset\RR^d$ denote the region of control variables, which is convex and compact. Suppose a physical system $f_p(\bx)$ with $\bx\in \Omega$ is of interest. In order to study the response function $f_p(\cdot)$, physical experiments are conducted on a set of input points (or control variables) $\bX_n=\{\bx_1,...,\bx_n\}\subset\Omega$. In physical experiments, the responses are usually corrupted by noise, which may come from the natural uncertainties inherent to the complex systems, such as actuating uncertainty, controller fluctuation, and measurement error. Therefore, we assume that we observe $y_j^{(p)}$ on $\bx_j$, $j=1,...,n$, with relationship
\begin{align}\label{yp}
    y^{(p)}_j=f_p(\bx_j)+\varepsilon_j ,
\end{align}
where $\varepsilon_j$ are i.i.d. random variables with mean zero and finite variance. Model \eqref{yp} has been widely considered in calibration problems; see \cite{ko,tuo&wu2014} for example.

In many situations, actual physical experimentation can be costly and difficult, so scientists and engineers use 
computer models to study a system of interest. Let $f_s(\bx,\btheta)$ denote a computer model, where $\bx\in \Omega$ and $\btheta \in \bTheta$ with $\bTheta\subset\RR^q$. The space $\bTheta$ refers to the parameter space for the \textit{calibration parameter} $\btheta$, and we assume that $\bTheta$ is a compact region in $\RR^q$. Since computer models are constructed based on simplified assumptions, they rarely perfectly match the physical responses. One fundamental problem in computer experiments is \textit{calibration}, where the goal is to find a calibration parameter such that the computer model can approximate the physical response well. 

\cite{ko} supposes that
\begin{align}\label{eq:fpfsX}
    f_p(\bx) = f_s(\bx,\btheta_0) +\delta_0(\bx),
\end{align}
where $\btheta_0$ is the true value of the calibration parameter. 
\cite{ko} 
examines the statistical calibration problem by a Bayesian approach, where $f_s(\cdot,\cdot)$ and $\delta_0(\cdot)$ are assumed to be independent realizations of Gaussian processes. Then the calibration parameter is estimated via a Bayesian approach. However, from \eqref{eq:fpfsX} it can be seen that $\btheta_0$ is unidentifiable, because $(\btheta_0, \delta_0(\cdot))$ cannot be uniquely determined. Therefore, it is crucial to define an \textit{identifiable} calibration parameter in order to study the estimation problem. To this end, we rewrite the physical response $f_p$ as
\begin{align}\label{eq:fpfs}
    f_p(\bx) = f_s(\bx,\btheta) +\delta_{\btheta}(\bx),
\end{align}
where $\delta_{\btheta}$ is the discrepancy function. The {identifiable definition} of calibration parameter depends on the users' own criteria, and may not agree with each other.

\begin{figure}[h!]
\begin{center}
    \includegraphics[width=3.5in]{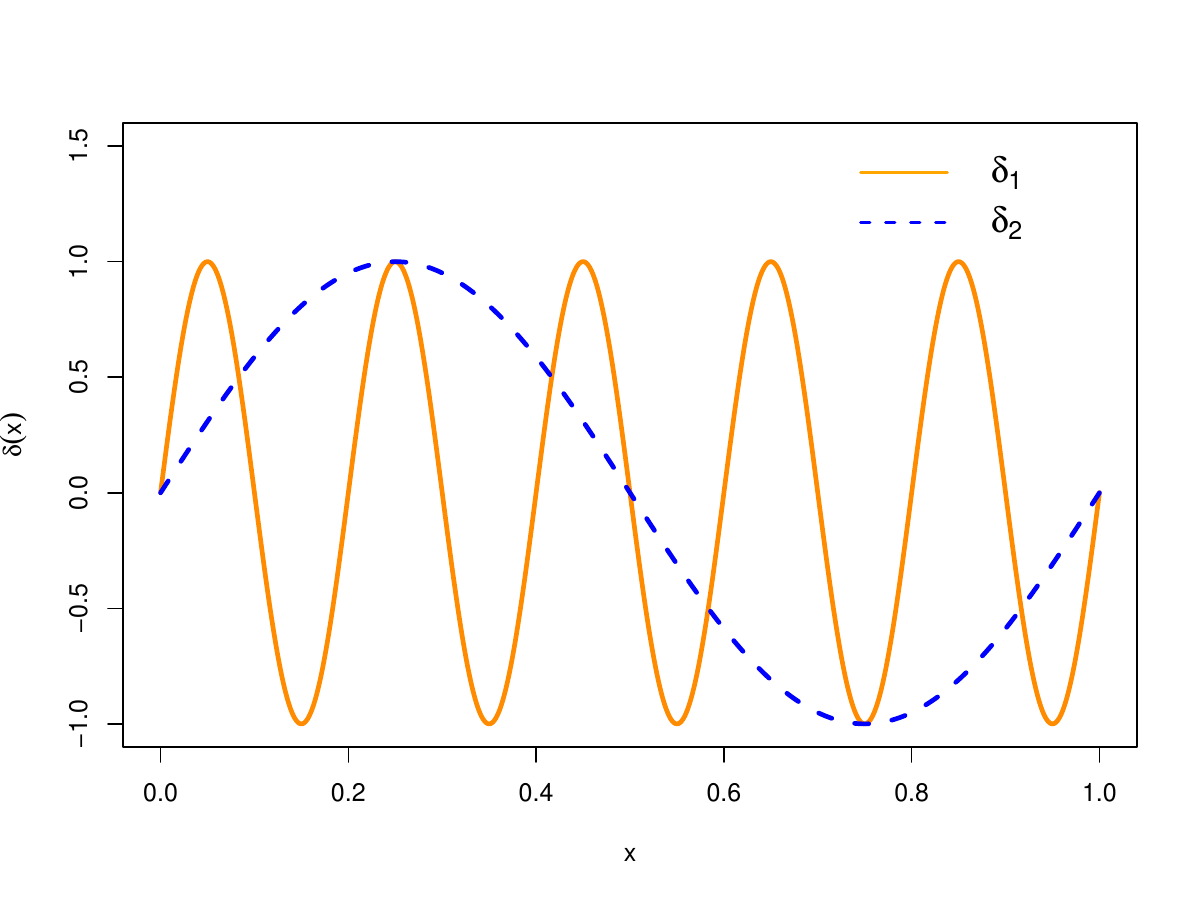}
\end{center}
    \caption{{Visualization of $\delta_1=\sin(10\pi x)$ and $\delta_2=\sin(2\pi x)$ over $\Omega=[0,1]$.}}
    \label{prop1:figure}
\end{figure}

One choice in the existing literature is the $L_2$ calibration, which minimizes the $L_2$ difference between $f_p(\cdot)$ and $f_s(\cdot,\btheta)$, as studied by \cite{tuo&wu2,tuo&wu2014,xie2020bayesian}. Consider the following example. Suppose we have two discrepancy functions $\delta_1=\sin(10\pi x)$ and $\delta_2=\sin(2\pi x)$ over $\Omega=[0,1]$, as shown in Figure \ref{prop1:figure}. Apparently, we have that $\|\delta_1\|_{L_2(\Omega)}=\|\delta_2\|_{L_2(\Omega)}=\sqrt{1/2}$, thus $L_2$ calibration provides no informative judgement between the two functions. However, we can state that for a discrepancy function, less vibration is preferred, which makes $\delta_2$ as a more desirable discrepancy function. 

{Motivated by the previous example}, we propose a novel calibration method, \textit{Sobolev calibration}, which generalizes the $L_2$ calibration and KO calibration, and can be applied to controlling both the $L_2$ norm and the vibration of the discrepancy function. Specifically, let $W^m(\Omega)$ be a Sobolev space\footnote{Although a Sobolev space is typically defined on an open set, by extension and restriction theorems \citep{devore1993besov,rychkov1999restrictions}, we can extend a Sobolev space on an open set to its closure.} with smoothness $m$. Here $m$ can be a non-integer, and the corresponding Sobolev space is called the \textit{fractional Sobolev space} \citep{adams2003sobolev}. For details, we refer to Section \ref{app:introtoSoboRKHS} in the supplement. Let $\cH_m(\Omega)$ be a function space such that $\cH_m(\Omega)=W^m(\Omega)$ with equivalent norms. We {provide the identifiable definition of }calibration parameter as 
\begin{align}\label{thetastar}
    \btheta_S^* = \argmin_{\btheta\in \bTheta}\|f_p(\cdot)-f_s(\cdot,\btheta)\|_{\cH_m(\Omega)},
\end{align}
where $\|\cdot\|_{\cH_m(\Omega)}$ is the norm of the space $\cH_m(\Omega)$. We will discuss the selection of $\cH_m(\Omega)$ in Section \ref{sec:Powerspace}. If $m=0$, then $W^m(\Omega)$ is equivalent to $L_2(\Omega)$, and \eqref{thetastar} reduces to $L_2$ calibration. In this work, we are interested in the statistical inference of $\btheta_S^*$.

\begin{remark}
{It can be shown that even for two discrepancy functions with the same $L_2$ norm, the difference of their Sobolev norms can be unlimited large, which introduces more fluctuation, as shown in the following proposition.}

\begin{proposition}\label{prop:rkhsl2infty}
For any constants $c>0$ and $m>d/2$, there exist two discrepancy functions $\delta_1,\delta_2\in W^{m}(\Omega)$ with $\frac{\|\delta_1\|_{W^m(\Omega)}}{\|\delta_2\|_{W^m(\Omega)}}>c$ while $\|\delta_1\|_{L_2(\Omega)}=\|\delta_2\|_{L_2(\Omega)}$.
\end{proposition}

This proposition implies that $L_2$ calibration can choose functions with large Sobolev norm. Since large Sobolev norm indicates substantial fluctuation which leads to distortion of original shape, this can be considered as a type of over-calibrating that Sobolev calibration can avoid. On the other hand, in Sobolev calibration, since $m$ can be flexibly chosen, it can be tuned or selected based on practicers’ needs, which avoids another risk that $L_2$ distance exceeds the satisfying level caused by overlarge $m$. In that way, our proposed method would provide additional flexibility by finer level of calibration. 
\end{remark}

\subsection{Estimation of Calibration Parameter $\btheta_S^* $}\label{sec:SobolevCali}

\subsubsection{Estimating $\btheta_S^*$ via a two-step approach}\label{subsec:SobolevCali}

With an appropriate selection of $\cH_m(\Omega)$ that is equivalent to $W^m(\Omega)$ (the selection of $\cH_m(\Omega)$ is to be introduced  in Section \ref{sec:Powerspace}), one can estimate $\btheta_S^*$ defined in \eqref{thetastar} based on the observations $(\bx_j,y_j^{(p)})$ of the physical experiments and the computer models $f_s(\cdot,\btheta)$. Following $L_2$ calibration, {we adopt a two-step approach. That is, we first estimate $f_p(\cdot)$ via some estimator, and then minimize the discrepancy between the estimated physical process and the surrogate model of computer code.} Among theoretical developments in frequentist calibration literature, it is typical to assume that $f_p(\cdot)$ is a \textit{deterministic function} \citep{tuo&wu2,tuo&wu2014,plumlee2017bayesian,sung2020calibration,tuo2019adjustments,xie2020bayesian}. For the ease of presentation, let us focus on the scenario that $f_p(\cdot)$ is a deterministic function in this section and Section \ref{sec:theory_deter}. We will consider the scenario where $f_p(\cdot)$ is a Gaussian process, which {is similar to} the original idea of KO, in Section \ref{sec:theory_gp}.

If $f_p(\cdot)$ is a deterministic function, one widely used approach to recover $f_p(\cdot)$ is \textit{kernel ridge regression} \citep{wahba1990}. Let $\Phi(\cdot,\cdot)$ be a symmetric and positive definite function and $\cN_{\Phi}(\Omega)$ be an RKHS generated by $\Phi$. For an introduction to RKHSs, we refer to Section \ref{app:introtoSoboRKHS} in the supplement. An estimator of $f_p(\cdot)$ via kernel ridge regression is given by 
\begin{align}\label{eq:krr}
\hat f_p(\cdot)\coloneqq\underset{f \in \cN_{\Phi}(\Omega)}\argmin \frac{1}{n} \sum_{j=1}^{n}\left(y_{j}^{(p)}-f\left(\bx_{j}\right)\right)^{2}+\lambda\|f\|_{\cN_{\Phi}(\Omega)}^{2},
\end{align}
where $\|\cdot\|_{\cN_{\Phi}(\Omega)}$ is the RKHS norm, and $\lambda$ is the smoothing parameter that can be chosen by certain criterion, for example, generalized cross validation \citep{wahba1990}. It follows from the representer theorem \citep{ScholkopfSmola02} that the optimal solution to \eqref{eq:krr} is 
\begin{align}\label{eq:krrest}
    \hat f_p(\cdot) = \rb(\cdot)^{\mathrm{T}}(\Rb+n\lambda \Ib_n )^{-1}\by^{(p)},
\end{align}
where $\rb(\cdot) = (\Phi(\cdot,\bx_1),\ldots,\Phi(\cdot,\bx_n))^{\mathrm{T}}$, 
$\Rb = (\Phi(\bx_j,\bx_k))_{j,k}\in\RR^{n\times n}$, $\by^{(p)}=(y_1^{(p)},...,y_n^{(p)})^{\mathrm{T}}$,
and $\Ib_n\in\RR^{n\times n}$ is an identity matrix. Since the physical experiments are limited, $n$ is not large and the closed form \eqref{eq:krrest} can be computed efficiently, which is one reason why kernel ridge regression is popular in practice.

The Sobolev calibration estimates $\btheta_S^*$ by
\begin{align}\label{eq:thetahat}
    \hat{\btheta}_S = \argmin_{\btheta\in \bTheta}\|\hat f_p(\cdot)-\hat f_s(\cdot,\btheta)\|_{\cH_m(\Omega)},
\end{align}
where $\hat{f}_s(\cdot,\btheta)$ is {a surrogate model (which is usually an approximation of $f_s(\cdot,\btheta)$)} for the computer model $f_s(\cdot, \btheta)$. Typically, computing $\hat{f}_s(\cdot,\btheta)$ is much faster than computing $f_s(\cdot, \btheta)$ (such computer experiment $f_s(\cdot, \btheta)$ is called \textit{expensive}), otherwise we can directly use $f_s(\cdot, \btheta)$ (such computer experiment $f_s(\cdot, \btheta)$ is called \textit{cheap}). We also assume that the approximation error of computer experiment is much smaller compared to that of physical experiments, which is reasonable considering the cost of conducting computer experiments is lower than the physical experiments in general. Several methods for constructing surrogate models include Gaussian process modeling \citep{santner2013design}, scattered data approximation \citep{wendland2004scattered}, and polynomial chaos approximation \citep{xiu2010numerical}. 

\subsubsection{Selection of Function Space $\cH_{m}(\Omega)$}\label{sec:Powerspace}

In order to define and estimate $\btheta_S^*$ in \eqref{thetastar}, we need to specify the function space $\cH_m(\Omega)$. We provide {three} approaches for selecting $\cH_m(\Omega)$ as follows. Again, we would like to note that since the identifiable definition of calibration parameter depends on the users' own criteria, the choice of the function space is also up to users' specific application.

\noindent\textbf{Approach 1: Setting $\cH_m(\Omega)$ as a Sobolev space.} 
If $m$ is a positive integer, one can compute the Sobolev norm by 
$$\|f\|_{\cH_m(\Omega)}^2 = \sum_{|\alpha| \leq m}\left(\int_{\Omega}\left|D^\alpha f(\bx)\right|^2 {\rm d} \bx\right),$$ where $D^\alpha f$ is the weak derivative of $f$. This approach is also considered by \cite{plumlee2017bayesian}, but asymptotic theory of calibration parameter is lacking. This approach, albeit its easy computation when $m$ is a positive integer, cannot be easily generalized to the case where $m$ is not a positive integer.

\noindent\textbf{Approach 2: Setting $\cH_m(\Omega)$ as an RKHS generated by a Mat\'ern kernel function.} Let $\Psi_m$ be the (isotropic) Mat\'ern family \citep{stein2012interpolation}. After a proper reparametrization, the Mat\'ern kernel function with $m>d/2$ is defined as 
\begin{align}\label{eq_matern}
	\Psi_m(\bx)=\frac{1}{\Gamma(m - d/2)2^{m - d/2 -1}}\|\bx\|_2^{m - d/2} B_{m - d/2}(\|\bx\|_2),
\end{align}
where $\|\cdot\|_2$ is the Euclidean distance, and $B_{m - d/2}$ is the modified Bessel function of the second kind. 
It can be shown that the Sobolev norm $W^m(\Omega)$ is equivalent to $\cN_{\Psi_m}(\Omega)$; see Corollary 10.48 of \cite{wendland2004scattered}.
Therefore, if $m >d/2$, it is natural to choose $\cH_m(\Omega)$ in \eqref{thetastar} as $\cN_{\Psi_m}(\Omega)$. 

\noindent\textbf{Approach 3: Setting $\cH_m(\Omega)$ as a power of an RKHS.} 
Approach 2 only works when the smoothness $m$ is larger than $d/2$. Therefore, we provide another approach for constructing $\cH_m(\Omega)$ such that $\cH_m(\Omega)=W^m(\Omega)$ and two norms are equivalent, which works for all $m\geq0$. This approach relies on the \textit{powers of RKHSs}. 

{Let $\cK_{m^*}(\cdot,\cdot)$ be a symmetric kernel function such that the RKHS $\cN_{\cK_{m^*}}(\Omega)$ coincides $W^{m^*}(\Omega)$ with $m^*>d/2$, and $m^*$ is an integer. Such kernel function exists. For example, we can take $\cK_{m^*}(\cdot,\cdot) = \Psi_{m^*}(\cdot-\cdot)$, where $\Psi_{m^*}(\cdot-\cdot)$ is as defined in \eqref{eq_matern}. Since $\cK_{m^*}(\cdot,\cdot)$ is positive definite, Mercer's theorem implies that it possesses an absolutely and uniformly convergent representation as
\begin{align}\label{eq:eigenck1}
    \cK_{m^*}(\bs,\bt)=& \sum_{j=1}^\infty \gamma_{\cK,j} e_j(\bs)e_j(\bt), \forall \bs,\bt\in\Omega,
\end{align}
where $\gamma_{\cK,j}$ and $e_j$ are eigenvalues and eigenfunctions of $\cK_{m^*}(\cdot,\cdot)$, respectively. The following definition can be found in \cite{steinwart2012mercer,kanagawa2018gaussian}. See Definition 4.11 of \cite{kanagawa2018gaussian} for example. } 
{\begin{definition}[Powers of RKHS]\label{def:Prkhs}
Let $0\leq\beta\leq 1$ be a constant, and assume that $\sum_{j=1}^\infty \gamma_{\cK,j}^\beta e_j(\bs)^2 <\infty$ for all $\bs\in \Omega$, where  $\gamma_{\cK,j}$ and $e_j$ are as in \eqref{eq:eigenck1}. Then the $\beta$-th power of RKHS $\cN_{\cK_{m^*}}(\Omega)$ is defined as 
\begin{align}\label{eq:prkhs}
    \cN_{\cK_{m^*}}^\beta(\Omega) = \left\{f(\cdot)=\sum_{j=1}^\infty a_j\gamma_{\cK,j}^{\beta/2}e_j(\cdot): \sum_{j=1}^\infty a_j^2 < \infty\right\},
\end{align}
with inner product defined as
\begin{align*}
    \langle f,g \rangle_{\cN_{\cK_{m^*}}^\beta(\Omega)} = \sum_{j=1}^\infty a_jb_j, \mbox{ for } f(\cdot)=\sum_{j=1}^\infty a_j\gamma_{\cK,j}^{\beta/2}e_j(\cdot), g(\cdot)=\sum_{j=1}^\infty b_j\gamma_{\cK,j}^{\beta/2}e_j(\cdot), f,g\in \cN_{\cK_{m^*}}^\beta(\Omega).
\end{align*}
The $\beta$-th power of kernel function $\cK_{m^*}(\cdot,\cdot)$ is defined by 
\begin{align*}
    \cK_{m^*}^\beta(\bs,\bt)=& \sum_{j=1}^\infty \gamma_{\cK,j}^\beta e_j(\bs)e_j(\bt), \forall \bs,\bt\in\Omega.
\end{align*}
\end{definition}}
{By selecting $\beta = m/m^*$ and $\cH_m(\Omega)$ as $\cN_{\cK_{m^*}}^\beta(\Omega)$, it can be shown that $\cN_{\cK_{m^*}}^\beta(\Omega)=W^{m}(\Omega)$ with equivalent norms. For more details, we refer to Section \ref{app:powerRKHS} in the supplement. In particular, it can be seen that when $m>d/2$, the two norms in the objective functions of Approaches 2 and 3 are equivalent. We note that in practice, the kernel functions of Approaches 2 and 3 can have some adjusted parameters, e.g., scale parameter. These parameters can be tuned by practitioners for their own use. We numerically investigate the influence of the scale parameters in Section \ref{eg:adjustment}.}

\subsubsection{Computation}\label{sec:comp}

If we choose Approach 1 to select $\cH_m(\Omega)$, and $m$ is a positive integer, the Sobolev norm can be directly calculated by $\|\hat f_p(\cdot)-\hat f_s(\cdot,\btheta)\|_{\cH_m(\Omega)}^2=\sum_{|\alpha| \leq m}\left(\int_{\Omega}\left|D^\alpha \left(\hat f_p(\bx)-\hat f_s(\bx,\btheta)\right)\right|^2 {\rm d} \bx\right)$.
If we choose Approach 2 or Approach 3 as stated in Section \ref{sec:Powerspace} to select $\cH_m(\Omega)$, there is no exact solution to compute $\|\hat f_p(\cdot)-\hat f_s(\cdot,\btheta)\|_{\cH_m(\Omega)}$. {Therefore, we approximate $\|\hat f_p(\cdot)-\hat f_s(\cdot,\btheta)\|_{\cH_m(\Omega)}$ by the following approach.} Let $\cK_{m}(\cdot,\cdot)$ be the kernel function of $\cH_m(\Omega)$. In order to approximate $\cH_m(\Omega)$, we first draw $N$ uniformly distributed points on $\Omega$, denoted by $\tilde \bX=\{\tilde\bx_1,...,\tilde\bx_N\}$. Suppose we observe $g(\tilde\bx_1,\btheta),...,g(\tilde\bx_N,\btheta)$ on $\tilde\bX$, where $g(\cdot,\btheta)\coloneqq \hat f_p(\cdot)-\hat f_s(\cdot,\btheta)$. Define
$
\mathcal{I}_{\cK_{m},\tilde \bX}g(\cdot,\btheta)=\tilde\rb(\cdot)^{\mathrm{T}} \tilde\Rb^{-1} \tilde\bg_{\btheta},
$
where $\tilde\rb(\cdot) = (\cK_{m}(\cdot,\tilde\bx_1),\ldots,\cK_{m}(\cdot,\tilde\bx_n))^{\mathrm{T}}$, 
$\tilde\Rb = (\cK_{m}(\tilde\bx_j,\tilde\bx_k))_{j,k}\in\RR^{N\times N}$, and $\tilde\bg_{\btheta}=(g(\tilde\bx_1,\btheta),...,g(\tilde\bx_N,\btheta))^{\mathrm{T}}$. 
Then we have that $\|\mathcal{I}_{\cK_{m},\tilde \bX}g(\cdot,\btheta)\|_{\cH_m(\Omega)}$ is a good approximation of $\|g(\cdot,\btheta)\|_{\cH_m(\Omega)}$  {(For theoretical justification of this statement, see Section \ref{sec:computation} in the supplement)}. Therefore, $\hat{\btheta}_S$ in \eqref{eq:thetahat} can be approximated by
$
    \tilde{\btheta}_{S} =\argmin_{\btheta\in \bTheta}\|\mathcal{I}_{\cK_{m},\tilde \bX}g(\cdot,\btheta)\|_{\cH_m(\Omega)}^2=
    \argmin_{\btheta\in \bTheta}\tilde\bg_{\btheta}^{\mathrm{T}} \tilde\Rb^{-1} \tilde\bg_{\btheta}.
$

In $L_2$ calibration, the objective function can be directly approximated via the empirical $l_2$ norm. The estimated calibration parameter in $L_2$ calibration, denoted by $\hat{\btheta}_{L_2}$, can be approximated by
$
    \tilde{\btheta}_{L_2} = \argmin_{\btheta\in \bTheta}\frac{1}{N}\sum_{j=1}^N \left(\hat f_p(\tilde\bx_j)-\hat f_s(\tilde\bx_j,\btheta)\right)^2.
$ When $\tilde\Rb=c\cdot I_{N\times N}$, we have $\tilde{\btheta}_{S}=
    \argmin_{\btheta\in \bTheta}\tilde\bg_{\btheta}^{\mathrm{T}} \tilde\bg_{\btheta}=\tilde{\btheta}_{L_2}$, which indicates that there is a natural connection between the Sobolev calibration and $L_2$ calibration. The rigorous development of the connection will be given in Section \ref{subsec:connect}. 

\begin{remark}
    One may also solve $\argmin_{\btheta\in \bTheta}\bg_{\btheta}^{\mathrm{T}} \Rb^{-1} \bg_{\btheta}$ to approximate $\hat\btheta_S$, where $\bg_{\btheta}=(g(\bx_1,\btheta),...,g(\bx_n,\btheta))^{\mathrm{T}}$ and $\Rb = (\cK_{m}(\bx_j,\bx_k))_{j,k}\in\RR^{n\times n}$, i.e., use the empirical RKHS norm based on the original data. However, it has been shown in \cite{tuo&wu2} that the empirical $L_2$ calibration is not semiparametric efficient. Therefore, we do not consider using the empirical RKHS norm to estimate $\btheta_S^*$ in the present work.
\end{remark}

\section{Theoretical Properties of the Sobolev Calibration}\label{sec:theory_deter}

In this section, we discuss the asymptotic behavior of $\hat{\btheta}_S$ in \eqref{eq:thetahat}. In the rest of this work, the following definitions are used. For two positive sequences $a_n$ and $b_n$, we write $a_n\asymp b_n$ if, for some constants $C,C'>0$, $C\leq a_n/b_n \leq C'$. We use $C,C',c_j,C_j, j\geq 0$ to denote generic positive constants, of which value can change from line to line. We assume that the observed points $\bx_1,\ldots,\bx_n$ are uniformly distributed on $\Omega$, {while we point out that our theory can be easily generalized to the case that $\bx_j$'s are independently drawn from a distribution with a density that is bounded away from zero and infinity. Such an extension may make mathematical development more involved, and may dilute our main focus in this work.} We focus on the case that $\cH_m(\Omega)$ is an RKHS generated by a symmetric kernel function $\cK_{m}(\cdot,\cdot)$, which can be ensured if we apply Approach 2 or Approach 3 as stated in Section \ref{sec:Powerspace}. Recall that $\hat f_p(\cdot)$, $\hat f_s(\cdot,\btheta)$, and $\hat \btheta_S$ can depend on $n$, where the dependency is suppressed for notational simplicity.

\subsection{Asymptotic Results of the Sobolev Calibration}\label{sec:determinaymptotic}

In this section, we show that the proposed Sobolev calibration enjoys some nice properties as $L_2$ calibration, including: 1) The estimated $\hat{\btheta}_S$ is consistent; 2) The convergence rate of $\hat{\btheta}_S$ is $\|\hat{\btheta}_S-\btheta_S^*\|_2=O_{\PP}(n^{-1/2})$; 3) $\sqrt{n}(\hat{\btheta}_S-\btheta_S^*)$ is asymptotic normal; 4) The Sobolev calibration is semiparametric efficient. For the briefness of this paper, we move all the conditions to Section \ref{app:condofdeter} of the supplement, while we point out that we consider the case where the computer model can even be \textit{rougher} than the physical model at some parameter space of $\btheta$, which has not been presented in the literature as far as we know.

We start with the consistency of the estimated calibration parameter $\hat\btheta_S$ as stated in the following proposition, whose proof is provided in Section \ref{app:pfpropconsist} of the supplement.

\begin{proposition}[Consistency of $\hat\btheta_S$]\label{prop:consist}
Suppose Conditions (C2), (C4) and (C6) are fulfilled. The estimated calibration parameter $\hat\btheta_S$ is a consistent estimator of $\btheta^*_S$, i.e., $\hat\btheta_S\rightarrow^{\PP} \btheta^*_S$.
\end{proposition}

Now we are ready to present the main theorems in this section, which state the convergence property and semiparametric efficiency of $\hat\btheta_S$. The proofs of Theorems \ref{thm:asympn} and \ref{thm:semieff} are relegated to Sections \ref{app:pfthmmain} and
\ref{app:pfthmsemieff} of the supplement, respectively.

\begin{theorem}\label{thm:asympn}
Under Conditions (C1)-(C6), we have that
\begin{align}\label{eq:thmmainconv}
\hat \btheta_S-\btheta_S^*=-2\Vb^{-1}\left(\frac{1}{n}\sum_{j=1}^{n} \varepsilon_{j}A_1(\bx_j)\right) + o_{\PP}(n^{-1/2}),
\end{align}
where $\Vb$ is as in Condition (C3), and
\begin{align*}
    A_1(\bx) = \sum_{j=1}^\infty\frac{1}{\gamma_j}\phi_j(\bx)\int_\Omega\frac{\partial f_s(\bz,{\btheta_S^*})}{\partial\btheta}\phi_j(\bz){\rm d}\bz,
\end{align*}
where $\gamma_j$ and $\phi_j$ are eigenvalues and eigenfunctions of $\cK_{m}(\bs,\bt)$, repectively.
\end{theorem}

Note that Condition (C5) implies that $A_1(\bx)<\infty$ almost everywhere. Theorem \ref{thm:asympn} directly implies the asymptotic normality of $\sqrt{n}(\hat\btheta_S-\btheta_S^*)$, provided that 
\begin{align}\label{W}
    \Wb\coloneqq\EE(A_1(\bx)A_1(\bx)^{T})
\end{align}
exists and is positive definite.  By the central limit theorem, we have
\begin{align*}
    \sqrt{n}(\hat \btheta_S-\btheta_S^*)\stackrel{d}{\rightarrow} N\left(0,4 \sigma^{2} \Vb^{-1} \Wb \Vb^{-1}\right).
\end{align*}
Theorem \ref{thm:asympn} can be regarded as an extension of Theorem 1 in \cite{tuo&wu2}, where it has been shown that $\hat \btheta_{L_2}$ obtained by the $L_2$ calibration is asymptotically normal.

\begin{theorem}\label{thm:semieff}
Under the Conditions of Theorem \ref{thm:asympn}, if $\varepsilon_i$ in \eqref{yp} follows a normal distribution for $i=1,...,n$, the Sobolev calibration $\hat\btheta_S$ is semiparametric efficient.
\end{theorem}

Calibration problem is a semiparametric model because it covers estimating the physical function given by \eqref{eq:krr} and estimating the $q$-dimensional calibration parameter given by \eqref{eq:thetahat}, where the former problem has an infinite dimensional parameter space, but the latter problem is only finite dimensional. 
For detailed discussion about semiparamteric models, we refer to \cite{bickel1993efficient}. The property of semiparametric efficiency implies that the semiparametric model has the same asymptotic variance as the statistical estimation with the same observed data in the finite parameter space. This is a desirable property because the semiparametric model shares the same estimation efficiency with the corresponding parametric model under less assumptions. 

We demonstrate that Sobolev calibration reaches the highest calibration parameter estimation efficiency under some regularity conditions, by linking the Sobolev calibration estimator to maximum likelihood estimator in the finite parameter model. Note that normality of random error is additionally required. This assumption is commonly used in modeling physical process \citep{wu2011experiments,tuo&wu2014}. 
However, if the random error is non-Guassian, the semiparametric efficiency can still be achieved by applying a similar treatment to physical process estimation as \cite{tuo&wu2014}.

\subsection{Discussion on Uncertainty Quantification}
In practice, the prediction of calibration parameters and the quantification of uncertainty are both crucial to calibration problems. Typically, the confidence intervals for $\btheta_S^*$ and point-wise confidence band for $f_s(\cdot,\btheta_S^*)$ and $f_p(\cdot)$ are of interest.
For frequentist calibration problem, which are considered as a semi-parametric problem, uncertainty is inferred from the probability distribution of the estimate $\hat{\btheta}$. Traditionally, bootstrap can be applied to estimate the distribution, as demonstrated in \cite{wong2017frequentist}. However, the authors also highlighted that the bootstrap confidence region may introduce additional bias and lead to incorrect asymptotic coverage, making it a sub-optimal approach.

Within the framework of Sobolev calibration, we can utilize the asymptotic normality of $\hat \btheta_S$ obtained in Theorem \ref{thm:asympn} to construct a confidence interval for  $\btheta_S^*$. In addition, since computer experiment is a deterministic function where the only uncertainty comes from the estimation procedure of the calibration parameter, the point-wise confidence band for $f_s(\cdot,\btheta_S^*)$ can be easily inferred at a new point $\bx$. Similarly, once ${\delta}_{\btheta}(\cdot)$ is estimated based on observations $\{(\bx_j,y_j^{(p)}-f_s(\bx_j,\btheta)\}_{j=1}^n$, the point-wise confidence band for $f_p(\cdot)$ can also be derived. Since the estimation error of ${\delta}_{\btheta}(\cdot)$ may introduce additional bias, we recommend using models like Gaussian process to estimate discrepancy function such that the true $\delta_{\btheta}(\cdot)$ can be covered with high probability. For example, if 95\% confidence intervals and bands are of interest, the details of computation steps are summarized in Algorithm \ref{alg:uq}, Section \ref{subsec_UQdetails} in the supplement, which can be easily generalized to $q$\% confidence interval. Due to semi-parametric efficiency shown in Theorem \ref{thm:semieff}, 
our proposed confidence interval is valid, which is empirically verified by numerical experiments in Section \ref{subsec_numUQ}.

\subsection{Connection with $L_2$ Calibration and KO Calibration}\label{subsec:connect}

In statistical calibration literature, two {widely used} methods are $L_2$ calibration and KO calibration, while only the former has been well studied from a theoretical perspective. In this section, we present how the Sobolev calibration serves as a bridge between these two methods. 

In $L_2$ calibration, the {identifiable} calibration parameter $\btheta_{L_2}^*$ is defined as the calibration parameter which minimizes the $L_2$ distance between the computer model $f_s(\cdot, \btheta)$ and the physical model $f_p(\cdot)$. Theorem 1 of \cite{tuo&wu2014} states that, under certain conditions, the estimate of $\btheta_{L_2}^*$ via $L_2$ calibration, denoted by $\hat\btheta_{L_2}$, satisfies
\begin{align}\label{eq:L2consist}
  \hat{\btheta}_{L_{2}}-\btheta_{L_2}^{*}=-2 \Vb_{L_2}^{-1}\left(\frac{1}{n} \sum_{j=1}^{n} \varepsilon_{j} \frac{\partial f_{s}}{\partial \btheta}\left(\bx_j, \btheta_{L_2}^{*}\right)\right)+o_{\PP}\left(n^{-1/2}\right),
\end{align}
where $$\Vb_{L_2}\coloneqq- \int_{\Omega}\frac{\partial^{2}}{\partial \btheta \partial \btheta^{\mathrm{T}}}(f_p(\bz)-f_s(\bz,\btheta_{L_2}^*))^2 {\rm d}\bz.$$ 
In fact, \eqref{eq:L2consist} can be directly obtained by taking $m=0$ in Theorem \ref{thm:asympn}, as stated in the following Corollary, whose proof is in Section \ref{app:pfofcorols} of the supplement.

\begin{corollary}[Theorem 1 of \cite{tuo&wu2014}]\label{coro:L2}
Suppose the conditions in Theorem \ref{thm:asympn} hold when taking $m=0$. Then \eqref{eq:L2consist} holds.
\end{corollary}

Now we consider KO calibration. Although KO calibration is originally established from a Bayesian perspective, \cite{tuo&wu2} {studied it from a frequentist perspective under some assumptions including the physical observations are noiseless, and called it frequentist KO calibration.}
By skipping the prior and using maximum likelihood estimation to estimate the calibration parameter, the frequentist KO calibration is given by
\begin{align}\label{eq:kodefest}
    \hat\btheta_{KO}=\argmin_{\btheta \in \bTheta}\|\hat f_p(\cdot)-\hat f_s(\cdot, \btheta)\|_{\cN_{\Phi}(\Omega)},
\end{align}
with the {identifiable calibration parameter} defined as 
\begin{align}\label{eq:kodef}
    \btheta_{KO}^*=\argmin_{\btheta \in \bTheta}\|f_p(\cdot)-f_s(\cdot, \btheta)\|_{\cN_{\Phi}(\Omega)}.
\end{align}
\cite{tuo&wu2} studied the asymptotic behavior of $\hat\btheta_{KO}$ and showed that it converges to $\btheta_{KO}^*$ in probability. \cite{tuo2020improved} further extend the frequentist KO calibration to a more general case where the physical observations can be disturbed by noise. Although the frequentist KO calibration in \cite{tuo2020improved} is slightly different with \eqref{eq:kodefest}, it still converges to \eqref{eq:kodef}; See Corollary 3.6 of \cite{tuo2020improved}. 

However, it has not been shown that whether KO calibration has similar properties as $L_2$ calibration, for example, the convergence rate $O_{\PP}(n^{-1/2})$ or the asymptotic normality. To the best of our knowledge, the only theoretical results are in \cite{tuo&wu2} and \cite{tuo2020improved}, under the assumption that $\cN_{\Phi}(\Omega)$ is equivalent to some Sobolev space $W^{m_1}(\Omega)$: \cite{tuo&wu2} showed the consistency of $\hat\btheta_{KO}$ under noiseless case, while \cite{tuo2020improved} provided a convergence rate which can be much slower than $O_{\PP}(n^{-1/2})$.

As a corollary of Theorem \ref{thm:asympn}, if $\cN_{\Phi}(\Omega)$ is equivalent to some Sobolev space $W^{m_1}(\Omega)$, the KO calibration parameter enjoys similar desired statistical properties of the $L_2$ calibration parameter. Corollary \ref{coro:ko} is a direct result of Theorem \ref{thm:asympn}, thus the proof is omitted.
In the following context of this paper, we refer to \eqref{eq:kodefest} as another version of the frequentist KO calibration, or KO calibration for brevity, since it also converges to {$\btheta^*_{KO}$} (but with a faster convergence rate).

\begin{corollary}\label{coro:ko}
Suppose the conditions in Theorem \ref{thm:asympn} hold when taking $m=m_1$. Then we have
\begin{align*}
    \hat{\btheta}_{K O}-\btheta_{K O}^*=-2\Vb_{KO}^{-1}\left(\frac{1}{n}\sum_{j=1}^{n} \varepsilon_{j}A_{KO}(\bx_j)\right) + o_{\PP}(n^{-1/2}),
\end{align*}
where
\begin{align*}
    A_{KO}(\bx) = \sum_{j=1}^\infty\frac{1}{\gamma_{\Phi,j}}\varphi_j(\bx)\int_\Omega\frac{\partial f_s(\bz,{\btheta_{KO}^*})}{\partial\btheta}\varphi_j(\bz){\rm d}\bz,
\end{align*}
and 
\begin{align*}
    \Vb_{KO}=-\sum_{j=1}^\infty \frac{1}{{\gamma_{\Phi,j}}}\frac{\partial^{2}}{\partial \btheta \partial \btheta^{\mathrm{T}}}\left(\int_{\Omega} (f_p(\bz)-f_s(\bz,\btheta_{KO}^*))\varphi_j(\bz) {\rm d} \bz\right)^2,
\end{align*}
where $\gamma_{\Phi,j}$ and $\varphi_{j}(\cdot)$ are eigenvalues and eigenfunctions of $\Phi(\cdot,\cdot)$, respectively.
\end{corollary}

Corollary \ref{coro:ko} builds a connection between the Sobolev calibration and KO calibration, and a new justification is provided on why KO calibration is widely used and has good performance. To the best of our knowledge, this is the \textit{first} result of this kind for KO calibration. 

From Corollaries \ref{coro:L2} and \ref{coro:ko}, it can be seen that both $L_2$ calibration and KO calibration are special cases of Sobolev calibration. Sobolev calibration allows the practitioners to select the intermediate calibration method of $L_2$ calibration and KO calibration, {since $m$ can be flexibly chosen between 0 and $m_1>d/2$, which avoids both the risks of violent fluctuation caused by over-small $m$ and value deviation caused by over-large $m$.} 
Thus, the Sobolev calibration is more flexible.
We believe that the Sobolev calibration is not only of theoretical interest, but also has potential wide applications in practice.

\section{Extension to Stochastic Physical Experiments}\label{sec:theory_gp}

In Section \ref{sec:theory_deter}, we consider the case that the physical experiment $f_p(\cdot)$ is a deterministic function. In this section, we extend the scenario to that $f_p(\cdot)$ admits a Gaussian process, in order to align with the original idea of KO calibration, where Gaussian process is used as the prior distribution of physical process. However, since the underlying true physical experiment is random, the Gaussian process scenario introduces additional challenges in the theoretical investigation of the Sobolev calibration. For example, the support of a Gaussian process is typically larger than the corresponding reproducing kernel Hilbert space \citep{van2008reproducing}. To the best of our knowledge, we are not aware of any theoretical investigation on the original KO calibration idea, where the underlying truth is indeed a Gaussian process. Therefore, by investigating the asymptotic behavior of Sobolev calibration (thus $L_2$ calibration and KO calibration) under the Gaussian process setting, our results not only fill the theoretical gap, but also provide theoretical justification on the use of Gaussian process model in the calibration problems.

Suppose we observe $y_j^{(p)}$ on $\bx_j\in \Omega$, $j=1,...,n$, with relationship as in \eqref{yp}, and $\varepsilon_j$ are i.i.d. random variables with mean zero and finite variance $\sigma_\epsilon^2$. Under the Gaussian process settings, the underlying true physical experiment $f_p(\cdot)$ is assumed to be a realization of a Gaussian process $Z(\cdot)$. From this point of view, we will not differentiate $f_p(\cdot)$ and $Z(\cdot)$ in this work. For the ease of mathematical treatment, we assume $Z(\cdot)$ has mean zero, variance $\sigma^2$ and correlation function $\cK(\cdot, \cdot)$, and denote it by $Z(\cdot)\sim {\rm GP}(0, \sigma^{2}\cK(\cdot, \cdot))$. We further assume $\cK(\cdot, \cdot)$ is positive definite and integrable on $\mathbb{R}^{d}$, satisfying $\cK(0)=1$, and \textit{stationary}, in the sense that the correlation between $Z(\bx)$ and $Z(\bx')$ depends only on the difference $\bx-\bx'$ between the two input variables $\bx$ and $\bx'$. Thus, we can rewrite $\cK(\cdot, \cdot)$ as $\cK(\cdot - \cdot)$. 

If $\varepsilon_j$'s are Gaussian random variables, then given $\by^{(p)}=\left(y_1^{(p)}, ..., y_n^{(p)}\right)^{\mathrm{T}}$, the conditional distribution of $Z(\bx)$ on a point $\bx\in \Omega$ is a normal distribution with conditional mean and variance given by
\begin{align}
\mathbb{E}[f_p(\bx) \mid \by^{(p)}] &=\rb_1(\bx)^\mathrm{T}\left(\Rb_1+\mu \Ib_{n}\right)^{-1} \by^{(p)}, \label{eq:fpexpectation}\\
\operatorname{Var}[f_p(\bx) \mid \by^{(p)}] &=\sigma^{2}\left(\cK(\bx-\bx)-{\rb_1(\bx)^{\mathrm{T}}}\left(\Rb_1+\mu \Ib_{n}\right)^{-1} \rb_1(\bx)\right),\label{eq:fpvar}
\end{align}
where $\rb_1({\bx})=\left(\cK\left(\bx-\bx_{1}\right), \ldots, \cK\left(\bx-\bx_{n}\right)\right)^{\mathrm{T}},\Rb_1=\left(\cK\left(\bx_{j}-\bx_{k}\right)\right)_{j k}$, $\Ib_n\in\RR^{n\times n}$ is an identity matrix, and $\mu=\sigma_\epsilon^2/\sigma^2$. The conditional mean \eqref{eq:fpexpectation} is a natural predictor of $Z(\bx)$, and the conditional variance \eqref{eq:fpvar} can be used to construct confidence intervals for statistical uncertainty quantification.

It is well-known that the conditional expectation \eqref{eq:fpexpectation} is the best linear predictor for $f_p(\bx)$, in the sense that it has the minimal mean squared prediction error \citep{gramacy2020surrogates}. Therefore, we define 
\begin{align}\label{eq:gpfhatp}
\hat f_p(\bx)\coloneqq\mathbb{E}[f_p(\bx) \mid \by^{(p)}]
\end{align}
as an estimate of $f_p(\bx)$ on a point $\bx\in\Omega$.

Although the true calibration parameter in \eqref{eq:fpfsX} $\btheta_0$ is fixed, recall that because of the identifiability issue as discussed in Section \ref{subsec:background}, our goal is to estimate {$\btheta^*_{S}$} in the Sobolev calibration as in \eqref{thetastar}. However, unlike the deterministic physical experiments, where $f_p(\cdot)$ is a deterministic function and $\btheta_S^*$ is also deterministic, in this section, we consider $f_p(\cdot)$ is a random process, and thus, $\btheta_S^*$ is a random variable.
The computer model $f_s(\cdot,\btheta)$ remains unchanged as in the deterministic case, i.e., the computer model is assumed to be a deterministic function, which is for the ease of mathematical treatment. The case that the computer model is an independent Gaussian process requires additional assumptions on $f_s(\cdot,\btheta)$, and we leave it for the future study.
The Sobolev calibration under Gaussian process settings can be further defined in the same form as \eqref{eq:thetahat}, i.e.,
\begin{align}\label{eq:gpthetahat}
    \hat{\btheta}_S = \argmin_{\btheta\in \bTheta}\|\hat f_p(\cdot)-\hat f_s(\cdot,\btheta)\|_{\cH_m(\Omega)},
\end{align}
where $\hat f_p(\cdot)$ is as in \eqref{eq:gpfhatp}, and $\hat f_s(\cdot,\btheta)$ is a surrogate model of the computer model $f_s(\cdot, \btheta)$. 

Next, we show that under Gaussian process settings and mild conditions, the proposed Sobolev calibration still enjoys desired properties, including consistency and convergence rate at $O_\PP(n^{-1/2})$. For briefness, we move all the conditions to Section \ref{app:condofgp} of the supplement.

We start with the consistency of the estimated calibration parameter $\hat\btheta_S$ as stated in the following proposition, whose proof is provided in Section \ref{app:gp+pfpropconsist} of the supplement.

\begin{proposition}[Consistency of $\hat\btheta_S$]\label{prop:gp+consist}
Suppose Conditions (C2'), (C3') and (C6) are fulfilled. The estimated calibration parameter $\hat\btheta_S$ is a consistent estimator of $\btheta^*_S$, i.e., $\hat\btheta_S\rightarrow^{\PP} \btheta^*_S$.
\end{proposition}

Next, we present the main theorem in this section, which states the convergence property of $\hat\btheta_S$. The proof of Theorem \ref{thm:gp+asympn} is in Section \ref{app:gp+asympn} of the supplement.

\begin{theorem}\label{thm:gp+asympn}
Under Conditions (C1')-(C4'), and (C6), we have that
\begin{align}\label{eq:gp+thmmainconv}
\hat \btheta_S-\btheta_S^*=-2\Vb^{-1}\left(\frac{1}{n}\sum_{j=1}^{n} \varepsilon_{j}A_1(\bx_j)\right) + O_{\PP}(n^{-1/2}).
\end{align}
where $\Vb$ is as in Condition (C3), and $A_1(\bx)$ is as in Theorem \ref{thm:asympn}.
\end{theorem}
It is worth noting that in \eqref{eq:gp+thmmainconv}, we have $O_{\PP}(n^{-1/2})$ instead of $o_{\PP}(n^{-1/2})$ in \eqref{eq:thmmainconv}, which is because of the randomness of the Gaussian process. Theorem \ref{thm:gp+asympn} directly implies the convergence rate of $\hat \btheta_S$. Specifically, we have $\|\hat{\btheta}_S-\btheta_S^*\|_2=O_{\PP}(n^{-1/2})$. Therefore, even if $f_p(\cdot)$ is a random function following a Gaussian process, we can still guarantee the convergence rate of $\hat \btheta_S$ as $O_{\PP}(n^{-1/2})$.  

\section{Numerical Experiments}\label{sec:eg}

In this section, we conduct numerical experiments to evaluate the estimation performance of our proposed method under finite samples, and compare it with the frequentist KO calibration and $L_2$ calibration. 
For the conciseness of this paper, we move all the figures to Section \ref{app:figures} in the supplement.

\subsection{Simulation Studies}\label{subsec:simu}

We consider two examples in this subsection. The underlying physical model is a deterministic function in Example \ref{dtm-case2}, and is a Gaussian process in Example \ref{gp-case3}. An additional example where the underlying physical model is a deterministic function and has the same form as in Example 2 is put in Section \ref{app_exp} in the supplement.

\begin{example}\label{dtm-case2}
In this example, we show that in some cases, the $L_2$ calibration can lead to a ``wiggly'' function, by showing an unstable case where the smoothness of calibrated function is sensitive to calibration parameter.

Suppose the computer experiment is
\begin{align}
f_s(x, \theta)=f_p(x)-1/5\sqrt{\theta^{2}-\theta+1}\left(\sin\left(2\pi x\frac{10}{\exp(\theta^2+\theta)}\right)+\frac{\cos(2\pi x 10)}{2 \exp\theta^2+1}+\theta^2\right) \nonumber.
\end{align}
For the detailed experiment setup, we refer to Section \ref{app_exp} in the supplement.

The corresponding Sobolev calibration parameter $\theta^*_{S}$ is defined as in \eqref{thetastar}, which is $\theta_S^*\approx 1.3$ by numerical optimization.
From Figure \ref{fig:Eg2} (also see Figure \ref{intro:eg2}), it can be seen that $L_2$ calibrated experiment should be extraordinarily curved in order to achieve the smallest $L_2$ distance, which may not be desired in practice, and can sabotage the coverage accuracy of $f_p(\cdot)$, as shown in Section \ref{subsec_numUQ}. KO calibration, on the other hand, keeps the shape of the physical experiment well, while leads to a larger $L_2$ discrepancy. Compared with $L_2$ calibration and KO calibration, Sobolev calibration provides a preferable choice by approximating the function value well while maintaining the shape of physical process.

The corresponding summary statistics of numerical experiment results are shown in Table \ref{tab:eg2}. It can be seen that Sobolev calibration shows a good estimation and has a small statistical estimation error, which empirically verifies the theoretical property of our proposed method. 

We further specify a range of $m$ from small to large, in order to further illustrate the flexibility of Sobolev calibration and how it generalizes the other two methods. The Sobolev calibrated experiments present proximity to $L_2$ calibration and KO calibration when $m=7/8$ and $m=9/5$, respectively, as shown in Figure \ref{fig:eg2-add}. Therefore, one can choose their preferred $m$ to attain the satisfying calibrated experiment under the theoretical guarantee. Specifically, the practitioner can choose $m$ at a relatively small value to emphasize more on point-wise value approximation, or choose a relatively large $m$ to pursue better approximation in shape.

\end{example}

\begin{table}[!h]
 \caption{\label{tab:eg2}The mean and standard deviation (SD) of estimated calibration parameters for different methods in Example \ref{dtm-case2}, where $\theta^*_S$ is 1.3.}
\begin{center}
\begin{tabular}{lccccc}
\hline
& \multicolumn{2}{c}{$\sigma^2=0.05$} &  & \multicolumn{2}{c}{$\sigma^2=0.1$} \\ \cline{2-3} \cline{5-6} 
\multicolumn{1}{c}{}                            & Mean        & SD        &  & Mean       & SD       \\ \hline
Sobolev                                         & $1.3321$           & $0.0767$         &   & $1.3279$           & $0.0921$        \\
$L_2$                                             & $0.0422$           & $0.0112$          &  & $0.0484$           & $0.0545$        \\
KO                                              & $2.1354$           & 0.0191         &  & 2.1354           & 0.0268            \\ \hline
\end{tabular}
\end{center}
\end{table}

\begin{example}\label{gp-case3}
In this example, we choose the underlying physical experiment $f_p(\cdot)$ to be a Gaussian process. In Section \ref{sec:theory_gp}, we assume that the mean function of $f_p(\cdot)$ is zero to simplify the theoretical development, while in this numerical experiment, we set the mean function dependent on the input $x$ to endow it with practical meaning. This act would not hurt the theoretical merit of Sobolev calibration.

Suppose the underlying physical process is a Gaussian process given by\\
$
    f_p(x) \sim {\rm GP}\left(\exp (\pi x / 5) \sin 2 \pi x, 0.1 \Psi_{2}\right)
$
for $x\in \Omega=[0,1]$.
The physical observations are given by 
$
    y^{(p)}_j=f_p(x_j)+\varepsilon_j,
$
with $x_j\sim$ Uniform$(0,1)$, $\varepsilon_j\sim N(0,\sigma_\epsilon^2)$ for $j=1,...,200$. Suppose the computer experiment is
\begin{align*}
f_s(x, \theta)=\exp (\pi x / 5) \sin 2 \pi x-\sqrt{\theta^{2}-\theta+1}\left(\sin 2 \pi \theta x+\exp\left(\pi\theta x/2\right)\right).
\end{align*}
$\hat f_p(x)$ in \eqref{eq:gpfhatp} is applied to estimate $f_p(x)$, and $\mu$ is chosen through a validation set, where the details are put in Section \ref{exp-gpdetails} in the supplement.

Table \ref{tab:eg3} summarizes the numerical results of Example \ref{gp-case3}. Numerical summaries present that when $f_p(\cdot)$ follows a Guassian process, the mean and standard deviation of squared estimation error are both small. From Figure \ref{fig:Eg3}, we can see that Sobolev calibration can still maintain good performance even though the underlying physical process is a random Gaussian process. 

\begin{table}[!h]
\caption{\label{tab:eg3}The mean squared errors (MSE) and standard deviation (SD) of estimated calibration parameters with respect to $\theta^*_S$ for different methods in Example \ref{gp-case3}.}
\begin{center}
\begin{tabular}{cccccc}
\hline
\multicolumn{1}{c}{\multirow{3}{*}{}} & \multicolumn{2}{c}{$\sigma^2=0.05$} &  & \multicolumn{2}{c}{$\sigma^2=0.1$} \\  \cline{2-3} \cline{5-6} 
\multicolumn{1}{c}{}                            & MSE        & SD        &  & MSE       & SD       \\ \hline
Sobolev                                         & $0.0003$           & $0.0008$         &   & $0.0010$           & $0.0058$        \\
$L_2$                                             & $  0.0857$           & $0.0921$          &  & $0.0861$           & $0.0925$        \\
KO                                              & $0.0013$           & $0.0018$         &  & $0.0039$           & 0.0248           \\ \hline
\end{tabular}
\end{center}
\end{table}

\end{example}

\subsubsection{Uncertainty Quantification}\label{subsec_numUQ}

In this subsection, we examine the empirical performance of calibration estimates by quantifying the uncertainty associated with calibration parameters $\hat{\btheta}_{L_2
}$, $\hat{\btheta}_S$, and $\hat{\btheta}_{KO}$ that experimented in the two examples of simulation studies. We compute the estimated calibration parameters along with 95\% confidence intervals, which account for uncertainty in estimation. Furthermore, we compute the 95\% confidence bands of computer experiment together with model discrepancy as the predictor of physical process, which account for uncertainty in prediction. Details are provided in Section \ref{subsec_UQdetails} in the supplement. Simulation results also include the corresponding interval lengths together with interval scores. 

Table \ref{tab:CI} summarizes the coverage results of the calibration parameters over 500 trials with $\sigma^2=0.1$. It shows that the overall coverage rates approach the nominal coverage 95\% with relatively narrow intervals, leading to small interval scores. This empirically supports that our construction of confidence interval via asymptotic properties is valid. Table \ref{tab:CI_function} presents the average point-wise coverage results of the computer and physical models, 
%point-wise average coverage rate of $f_s(\cdot,\hat{\theta})$ with respect to $f_s(\cdot,\theta^*)$ and $f_s(\cdot,\hat{\theta})+\delta(\cdot,\hat{\theta})$ with respect to $f_p(\cdot)$
computed at 100 equally spaced points on $\Omega$. Notably, all results are satisfactory, with the exception of $L_2$ calibration coverage rate for the physical process prediction in Example 1. This is due to the fact that $L_2$ calibration produces a curvy $f_s(\cdot,\hat{\theta}_{L_2})$, resulting in a similarly wiggly $\delta(\cdot,\hat{\theta}_{L_2})$. Consequently, learning a discrepancy function like that based on noisy data becomes challenging, since it is difficult to differentiate between frequent fluctuations and noise, which results in inaccurate predictions. This empirical phenomenon provides additional evidence for the importance of shape resemblance between the physical function and computer experiment, particularly in the context of physical process prediction.

\begin{table}[!h]
\caption{\label{tab:CI}Average coverage rates, interval lengths and scores of the calibration parameters $\theta^*_{L_2},\theta^*_{S}$, and $\theta^*_{KO}$, respectively.}
\centering
\begin{tabular}{llcccclccc}
\hline
                         & \multicolumn{4}{c}{{Example 1}}     &  & \multicolumn{4}{c}{Example 2}                                                          \\ \cline{2-5} \cline{7-10} 
                         &       & {coverage} & {length} & {score} &  & \multicolumn{1}{l}{}    & {coverage}           & {length}             & {score}              \\
{$L_2$}                       & {$m=0$}   & {0.9740}        & {0.0494}      & {0.3120}     &  & \multicolumn{1}{l}{{$m=0$}} & {0.9580}                  & {0.0294}                  & {0.2589}                  \\ \hline
\multirow{4}{*}{{Sobolev}} & {$m=7/8$} & {0.9760}        & {0.6960}      & {0.7435}     &  & \multirow{4}{*}{{$m=1$}}    & \multirow{4}{*}{{0.9760}} & \multirow{4}{*}{{0.1212}} & \multirow{4}{*}{{0.1652}} \\
                         & {$m=1$}   & {0.9320}        & {0.5896}      & {0.7295}     &  &                         &                    &                    &                    \\
                         & {$m=9/8$} & {0.9980}        & {0.8403}      & {0.8405}     &  &                         &                    &                    &                    \\
                         & {$m=9/5$} & {1}        & {0.3322}     & {0.3322}     &  &                         &                    &                    &                    \\ \hline
{KO}                       & {$m=2$}   & {1}        & {0.4776}      & {0.4776}     &  & \multicolumn{1}{l}{{$m=2$}} & {0.9900}                  & {0.4448}                  & {0.5735}                  \\ \hline
\end{tabular}
\end{table}

\begin{table}[!h]
 \caption{\label{tab:CI_function} Point-wise average coverage rates, interval lengths and scores of $f_s(\cdot,\theta^*)$ and $f_p(\cdot)$.}
\centering
\begin{tabular}{llccccccc}
\hline
                                                                      &          & \multicolumn{3}{c}{{Example 1}} & \multicolumn{1}{l}{} & \multicolumn{3}{c}{{Example 2}} \\ \cline{3-5} \cline{7-9} 
                                                                      &          & {$L_2$}  & {Sobolev}  & {KO} &                      & {$L_2$}  & {Sobolev}  & {KO} \\ 
\multicolumn{1}{c}{}                                                                                             &         & {$m=0$}        & {$m=1$}          & {$m=2$}     &                      & {$m=0$}        & {$m=1$}          & {$m=2$}     \\\hline
\multirow{3}{*}{{$f_s(\cdot,\theta^*)$}}                                & {coverage} & {0.9677}      &  {0.9174}               & {1}  &                      & {0.9451}      & { 0.9731}               & {0.9885}  \\
                                                                      & {length}   & {0.1755}      & {0.6996}               & {1.0687}  &                      & {0.0687}      & {0.3816}               & {1.2384}  \\
                                                                      & {score}    & {0.2497}      & {0.9079}               & {1.0687}  &                      & {0.4042}      & {0.4918}               & {1.4927}  \\ \hline
\multirow{3}{*}{\begin{tabular}[c]{@{}l@{}}{$f_p(\cdot)$} 
\end{tabular}} & {coverage} & {0.8255}      &  {0.9821}               & {0.9735}  &                      & {0.9506}      & {0.9581}               & {0.9657}  \\
                                                                      & {length}   & {0.3535}      & {0.2267}               & {0.1818}  &                      & {0.1802}      & {0.1821}               & {0.1906}  \\
                                                                      & {score}    & {0.7231}      & {0.2363}               & {0.1959}  &                      &  {0.2317}      & {0.2082}               & {0.2108}  \\ \hline
\end{tabular}
\end{table}

{\subsubsection{Adjustment to Scale Parameters}\label{eg:adjustment}}

Throughout the theoretical development of Sobolev calibration, we only require that the function space $\mathcal{H}_m(\Omega)$ is chosen as equivalent to the Sobolev space $W^m(\Omega)$, without restrictions on the specification of kernel function parameters. 
In this subsection, we numerically investigate how adjusting the kernel function's parameters affects the calibration parameter estimation in this subsection. Using the setup from Example \ref{dtm-case2} with $m=1$ and $\sigma^2=0.1$, we consider four different length-scale parameter values $\gamma=\{0.01,0.1,1,10\}$, for the Mat\'ern kernel function $\Psi(x,x^\prime)=\exp \left(-\gamma\left|x-x^\prime\right|\right)$.
The simulation results are summarized in Table \ref{tab:scale-para}, where the mean and standard deviation of $\hat{\theta}_S$ are computed over 500 trials. Numerical results suggest that under a wide range of length-scale parameters, Sobolev calibration performs reasonably well with the Sobolev norms of discrepancy $\|\delta(\cdot,{\hat{\theta}_S})\|_{W^m(\Omega)}$ lying between 0.9 to 2.5, which are much smaller than the Sobolev norm of $L_2$ calibration (6.7834). Figure \ref{fig:scale} describes the calibrated computer experiments under various length-scale parameters. In this setting, since the aim is to avoid over-calibrating the smooth underlying physical function, we would recommend choosing $\gamma=0.1$, with the smallest Sobolev norm that serves as the measurement objective for length-scale parameters.

\begin{table}[H]
 \caption{\label{tab:scale-para}The mean and standard deviation (shown in parentheses) of $\hat{\theta}_S$, RKHS norm and Sobolev norm of discrepancy function for different length-scale parameters in Example \ref{dtm-case2}.}
\centering
\begin{tabular}{lccccc}
\hline
{$\gamma$} & {$\theta^*_S$} &{$\hat{\theta}_S$}    & {$\|\delta(\cdot,{\hat{\theta}_S})\|_{\cH_m(\Omega)}$} & {$\|\delta(\cdot,{\hat{\theta}_S})\|_{W^m(\Omega)}$} \\ \hline
{0.01}     & {2.0400}     & {1.0728(0.0005)} & {11.3731}    & {1.6456}       \\
{0.1}      & {1.7800}     & {1.5172(0.2575)} &{1.5871}    & {\textbf{0.9943}}       \\
{1}        & {1.3000}     & {1.3279(0.0921)} & {0.9077}    & {1.1015}       \\
{10}       & {0.8800}     & {0.8109(0.0914)} & {0.7393}    & {2.4583}       \\ \hline
\end{tabular}
\end{table}

\subsection{Ion Channel Example}\label{subsec:real}

In this subsection, we introduce a real-world example, which is widely used for calibration experiment; see \cite{plumlee2016calibrating,plumlee2017bayesian,xie2020bayesian} for example. The dataset is collected by \cite{ednie2011modulation} for whole cell voltage clamp experiments on the sodium ion channels of cardiac cell membranes. For more details of the experiment, we refer to \cite{plumlee2016calibrating}.

In this example, we consider log scale of time as input $x$, and obtain 19 normalized current records as output. The computer model is the classical Markov model
$f_{s}(x, \boldsymbol{\theta})=\mathbf{e}_{1}^{\mathrm{T}} \exp [\exp (x) A(\boldsymbol{\theta})] \mathbf{e}_{4}$,
where the first $\exp$ is matrix exponential function, $\mathbf{e}_i$ is the column vector with one at the $i$th element and zero otherwise, and $A(\boldsymbol{\theta})$ with $\boldsymbol{\theta}=\left[\theta_1,\theta_2,\theta_3\right]^{\mathrm{T}}$ is defined as a 4 $\times$ 4 matrix, where $A(\boldsymbol{\theta})_{[1,1]}=-\theta_{2}-\theta_{3}$, $A(\boldsymbol{\theta})_{[1,2]}=\theta_{1}$, $A(\boldsymbol{\theta})_{[2,1]}=\theta_{2}$, $A(\boldsymbol{\theta})_{[2,2]}=-\theta_{1}-\theta_{2}$, $A(\boldsymbol{\theta})_{[2,3]}=\theta_{1}$, $A(\boldsymbol{\theta})_{[3,2]}=\theta_{2}$, $A(\boldsymbol{\theta})_{[3,3]}=-\theta_{1}-\theta_{2}$, $A(\boldsymbol{\theta})_{[3,4]}=\theta_{1}$, $A(\boldsymbol{\theta})_{[4,3]}=\theta_{2}$, $A(\boldsymbol{\theta})_{[4,4]}=-\theta_{1}$, and other elements are zero.

We still consider $\cH_m(\Omega)$ and $\cH_{m_1}(\Omega)$ are RKHSs generated by Mat\'ern kernel functions $\Psi_{1}(x)$ and $\Psi_{2}(x)$ respectively. Figure \ref{fig:realeg} in the supplement shows calibrated computer models by three methods and the discrepancy between observations and computer models. It can be seen that $L_2$ calibration achieves a more accurate approximation at each observation but sacrifices control of oscillations considering the entire function. Sobolev calibration obtain smoother function by specifying $m=1$ to additionally focus on the first-order derivative. KO calibration loses more approximation accuracy. Therefore, Sobolev calibration provides more flexibility for users in practice as illustrated in this real-world example. 

\section{Conclusions and Discussion}\label{sec:conclusion}

In this work, we propose a novel calibration framework, called Sobolev calibration, which allows the practitioners to adjust the identifiable calibration parameter based on practical needs. The performance of calibration can be enhanced by rebalancing the approximation in function value and function shape. We show that the estimation of the calibration parameter is consistent, asymptotically normal and semiparametric efficient by rigorous theoretical analysis. Moreover, we investigate the case where the underlying physical experiment is a realization of Gaussian process, which aligns with the original idea of KO calibration. Even if the underlying physical experiment is a random process, we prove that the estimation error still has a convergence rate $O_\PP(n^{-1/2})$. 

Our work can be extended in several ways. First, our method is based on random design, where the control variables $\bx$ are sampled from uniform distribution. Since fixed designs are popular in computer experiments \citep{wu2011experiments}, the framework under fixed design is worth investigating, which requires future work. Second, the Sobolev calibration is a frequentist approach. Note that the recent works \citep{tuo2019adjustments,xie2020bayesian} succeed in extending the $L_2$ calibration to a Bayesian framework, it is worth to explore an extension to develop a Bayesian Sobolev calibration framework in the future.

\bibliographystyle{apalike}
\bibliography{ref}

\appendix

\newpage
\setcounter{page}{1}

\numberwithin{figure}{section}
\numberwithin{table}{section}
\numberwithin{example}{section}
\section*{Supplementary Material for ``Sobolev Calibration of Imperfect Computer Models''}

In this supplement, we provide the support material to supplement the main article, including general notations used in technical proofs, introduction to reproducing kernel Hilbert spaces (RKHSs) and Sobolev spaces, powers of RKHSs, conditions listed for propositions and theorems, details of technical proofs, experimental details and figures from numerical examples. 
 
\section{Notation}
We use $\langle \cdot, \cdot \rangle_n$ to denote the empirical inner product, which is defined by $$\langle f, g \rangle_n = \frac{1}{n}\sum_{k=1}^n f(x_k)g(x_k)$$ for two functions $f$ and $g$, and let $\|g\|_n^2 = \langle g, g \rangle_n$ be the empirical norm of function $g$. In particular, let $$\langle \epsilon, f \rangle_n = \frac{1}{n}\sum_{k=1}^n\epsilon_k f(x_k)$$ for a function $f$, where $\epsilon = (\epsilon_1,\ldots,\epsilon_n)^{\mathrm{T}}$. Let $a\vee b =\max(a,b)$ for two real numbers $a,b$. We use $H(\cdot, \mathcal{F}, \|\cdot\|)$ and $H_B(\cdot,\mathcal{F},\|\cdot\|)$ to denote the entropy number and the bracket entropy number of class $\mathcal{F}$ with the (empirical) norm $\|\cdot\|$, respectively. Let $\cB_{\cF}(r) = \{f\in \cF:\|f\|_{\cF}\leq r\}$ denote the ball with radius $r$ in $\cF$, where $r>0$. We say ``$\mathcal{F}_1$ and $\mathcal{F}_2$ are equivalent'' for two Hilbert spaces $\mathcal{F}_1$ and $\mathcal{F}_2$ if $\mathcal{F}_1=\mathcal{F}_2$ with equivalent norms. Through the proof, we assume Vol$(\Omega)=1$ for the ease of notational simplicity. We use $C,C',c_j,C_j, j\geq 0$ to denote generic positive constants, of which value can change from line to line.

\section{Introduction to Reproducing Kernel Hilbert Spaces and Sobolev Spaces}\label{app:introtoSoboRKHS}

Suppose $K:\Omega \times \Omega \rightarrow \RR$ is a symmetric positive definite kernel, where $\Omega \subset\RR^d$. Define the linear space
\begin{eqnarray}\label{FPhi}
F_{K}(\Omega)=\left\{\sum_{i=1}^n\beta_i K(\cdot,\bx_i):\beta_i\in \RR,\bx_i\in \Omega,n\in\mathbb{N}\right\},
\end{eqnarray}
and equip this space with the bilinear form
\begin{eqnarray}
\left\langle\sum_{i=1}^n\beta_i K(\cdot,\bx_i),\sum_{j=1}^m\gamma_j K(\cdot, \bx'_j)\right\rangle_K:=\sum_{i=1}^n\sum_{j=1}^m\beta_i\gamma_j K(\bx_i, \bx'_j).
\label{eq1.2}
\end{eqnarray}
Then the \emph{reproducing kernel Hilbert space} (RKHS) $\mathcal{N}_{K}(\Omega)$ generated by the kernel function $K$ is defined as the closure of $F_{K}(\Omega)$ under the inner product $\langle\cdot,\cdot\rangle_{K}$, and the norm of  $\mathcal{N}_{K}(\Omega)$  is $\| f\|_{\mathcal{N}_{K}(\Omega)}=\sqrt{\langle f,f\rangle_{\mathcal{N}_{K}(\Omega)}}$, where $\langle\cdot,\cdot\rangle_{\mathcal{N}_{K}(\Omega)}$ is induced by $\langle \cdot,\cdot\rangle_{K}$. We refer the readers to \cite{wendland2004scattered_ec} for more details. In particular, we have the following theorem, which gives another characterization of the RKHS when $K$ is defined by a stationary kernel function $\Phi$, via the Fourier transform of $\Phi$. For an integrable function $f\in L_1(\mathbb{R}^d)$, its Fourier transform is defined as $$\mathcal{F}(f)(\bomega)=(2\pi)^{-d/2}\int_{\mathbb{R}^d} f(\bx) e^{-i\bx^{\mathrm{T}}\bomega}{\rm d}\bx.$$
\begin{theorem}[Theorem 10.12 of \cite{wendland2004scattered_ec}]\label{thm:NativeSpace}
Let $\Phi$ be a positive definite kernel function which is continuous and integrable in $\RR^d$. Define 
	$$\mathcal{G}:=\{f\in L_2(\RR^d)\cap C(\RR^d):\cF(f)/\sqrt{\cF(\Phi)}\in L_2(\RR^d)\},$$
	with the inner product
	$$\langle f,g\rangle_{\mathcal{N}_\Phi(\RR^d)}=(2\pi)^{-d}\int_{\RR^d}\frac{\cF(f)(\bm \omega)\overline{\cF(g)(\bomega)}}{\cF(\Phi)(\bm \omega)}d \bm \omega.$$ Then $\mathcal{G} = \mathcal{N}_\Phi(\RR^d)$, and both inner products coincide.
\end{theorem}

The Sobolev norm for functions on the whole space is
\begin{eqnarray}\label{Sobolev}
\|f\|_{W^m(\mathbb{R}^d)}=\left(\int_{\mathbb{R}^d} |\mathcal{F} (f)(\bomega)|^2 (1+\|\bomega\|_2^2)^{m} {\rm d} \bomega\right)^{1/2}.
\end{eqnarray}
Then the Sobolev space $W^m(\RR^d)$ consists of functions with finite norm defined in \eqref{Sobolev}. This function space is also called a Bessel potential space, which is a complex interpolation space of $L_2(\RR^d)$ and $W^k(\RR^d)$ with $k>m$ an integer \citep{almeida2006characterization,gurka2007bessel}. The Sobolev space on $\Omega$ can be defined via restriction. 

Since the Fourier transform of $\Psi_m$ defined in \eqref{eq_matern} satisfies \citep{tuo&wu2}
\begin{align*}
   c_1(1+\|\bomega\|_2^2)^{-m} \leq  \mathcal{F}(\Psi_m)(\bomega) \leq c_2(1+\|\bomega\|_2^2)^{-m},
\end{align*}
from Theorem \ref{thm:NativeSpace} and \eqref{Sobolev} it can be seen that the RKHS $\mathcal{N}_{\Psi_m}(\RR^d)$ generated by $\Psi_m$ is equivalent to $W^m(\RR^d)$. By some restriction technique (see \cite{adams2003sobolev_ec} and \cite{wendland2004scattered_ec}), it can be further shown that $\mathcal{N}_{\Psi_m}(\Omega) = W^m(\Omega)$ with equivalent norms.

\section{Justification of Approximation in Section \ref{sec:comp}}\label{sec:computation}
Define the integral operator $\cT:L_2(\Omega)\rightarrow L_2(\Omega)$ by
\begin{align*}
    \cT f(\bx)=\int_\Omega \cK_{m}(\bx,\bx')f(\bx'){\rm d}\bx', f\in L_2(\Omega), \bx\in \Omega.
\end{align*}
If $g\in \cT(L_2(\Omega))$, the proof of Theorem 11.23 of \cite{wendland2004scattered} implies that
\begin{align}\label{eq:compghm}
\left|\|g(\cdot,\btheta)\|_{\cH_m(\Omega)} - \|\mathcal{I}_{\cK_{m},\tilde \bX}g(\cdot,\btheta)\|_{\cH_m(\Omega)}\right|\leq 
    \|g(\cdot,\btheta) - \mathcal{I}_{\cK_{m},\tilde \bX}g(\cdot,\btheta)\|_{\cH_m(\Omega)} \leq C\|P_{\cK_{m},\tilde \bX}\|_{L_2(\Omega)},
\end{align}
where the first inequality is because of the triangle inequality. The function $P_{\cK_{m},\tilde \bX}(\bx)$ is defined and can be bounded by
\begin{align*}
    P_{\cK_{m},\tilde \bX}(\bx)^2
    = & \min_{\bu=(u_1,...,u_N)^{\mathrm{T}}\in \RR^N} \cK_{m}(\bx,\bx) - 2\sum_{j=1}^N u_j \cK_{m}(\bx,\bx_j)+ \sum_{j,k=1}^N u_ju_k\cK_{m}(\bx_j,\bx_k)\nonumber\\
    \leq & \cK_{m}(\bx,\bx) - \cK_{m}(\bx,\tilde\bx_l),
\end{align*}
where $l$ is chosen such that $\tilde\bx_l=\argmin_{\tilde\bx_j\in\tilde \bX_N}\|\bx-\tilde\bx_j\|_2$. As $N$ goes to infinity, we have that $\min_{\tilde\bx_j\in\tilde \bX_N}\|\bx-\tilde\bx_j\|_2$ converges to zero for all $\bx\in\Omega$, because $\tilde\bx_1,...,\tilde\bx_N$ are uniformly distributed, and eventually spread well in $\Omega$. Therefore, as long as $\cK_{m}(\cdot,\cdot)$ is continuous, $P_{\cK_{m},\tilde \bX}(\bx)$ converges to zero as $N$ goes to infinity. We note that under certain conditions of kernel functions (differentiability, stationarity, etc.), sharper convergence rate of  $P_{\cK_{m},\tilde \bX}(\bx)$ can be obtained; see Chapter 11 of \cite{wendland2004scattered}. 

Since $P_{\cK_{m},\tilde \bX}(\bx)$ converges to zero for all $\bx\in\Omega$, together with \eqref{eq:compghm}, it can be seen that $\|\mathcal{I}_{\cK_{m},\tilde \bX}g(\cdot,\btheta)\|_{\cH_m(\Omega)}$ is a good approximation of $\|g(\cdot,\btheta)\|_{\cH_m(\Omega)}$. One advantage of this approximation is that the RKHS norm of $\mathcal{I}_{\cK_{m},\tilde \bX}g(\cdot,\btheta)$ can be directly computed by
$
    \|\mathcal{I}_{\cK_{m},\tilde \bX}g(\cdot,\btheta)\|^2_{\cH_m(\Omega)} = \tilde\bg_{\btheta}^{\mathrm{T}} \tilde\Rb^{-1} \tilde\bg_{\btheta}.
$
Therefore, $\hat{\btheta}_S$ in \eqref{eq:thetahat} can be approximated by
$
    \tilde{\btheta}_{S} = \argmin_{\btheta\in \bTheta}\tilde\bg_{\btheta}^{\mathrm{T}} \tilde\Rb^{-1} \tilde\bg_{\btheta}.
$

\section{Powers of Reproducing kernel Hilbert spaces}\label{app:powerRKHS}

Recall that $\cK_{m^*}(\cdot,\cdot)$ is a symmetric kernel function such that the RKHS $\cN_{\cK_{m^*}}(\Omega)$ coincides $W^{m^*}(\Omega)$ with $m^*>d/2$, and $m^*$ is an integer. Also recall that by Mercer's theorem, 
it possesses an absolutely and uniformly convergent representation as
\begin{align}\label{eq:eigenck1app}
    \cK_{m^*}(\bs,\bt)=& \sum_{j=1}^\infty \gamma_{\cK,j} e_j(\bs)e_j(\bt), \forall \bs,\bt\in\Omega,
\end{align}
where $\gamma_{\cK,j}$ and $e_j$ are eigenvalues and eigenfunctions of $\cK_{m^*}(\cdot,\cdot)$, respectively. Let $\cN_{\cK_{m^*}}^\beta(\Omega)$ be the $\beta$-th power of RKHS $\cN_{\cK_{m^*}}(\Omega)$, and $\cK_{m^*}(\cdot,\cdot)$ be the $\beta$-th power of kernel function, defined as in Definition \ref{def:Prkhs}.

Proposition 4.2 of \cite{steinwart2012mercer_ec} shows that $\cN_{\cK_{m^*}}^\beta(\Omega)$ is indeed an RKHS generated by the kernel function $\cK_{m^*}^\beta(\cdot,\cdot)$. Furthermore, we can compactly embedded this space into $L_2(\Omega)$. With a slight abuse of notation, we still use $\cN_{\cK_{m^*}}^\beta(\Omega)$ to denote the embedded space, where the eigenfunctions can be different on a set with Lebesgue measure zero. As pointed by \cite{steinwart2012mercer_ec}, even if $\sum_{j=1}^\infty \gamma_{\cK,j}^\beta e_j(\bs)^2 <\infty$ for all $\bs\in \Omega$ does not hold, we can always define the space 
$\cN_{\cK_{m^*}}^\beta(\Omega)$ by the way as in \eqref{eq:prkhs}. 

An important property of $\cN_{\cK_{m^*}}^\beta(\Omega)$ is that the space $\cN_{\cK_{m^*}}^\beta(\Omega)$ is equivalent to a real interpolation space of $L_2(\Omega)$ and $W^{m^*}(\Omega)$, denoted by $\left[L_2(\Omega),W^{m^*}(\Omega)\right]_{\beta,2}$ \citep{steinwart2012mercer_ec}. The space defined by the real interpolation method is called Besov space, and denoted by $B_{2,2}^{r} = \left[L_2(\Omega),W^{m^*}(\Omega)\right]_{r/m^*,2}$ for $m^*>r>0$. In general, the real interpolation space (used for constructing Besov space) and the complex interpolation space (used for constructing Bessel potential space with norm \eqref{Sobolev}) are different, but fortunately in our case, the Sobolev space with norm \eqref{Sobolev} coincides $B_{2,2}^{m}$ with equivalent norms \citep{edmunds2008function_ec}. Putting all things together, we can see that $\cN_{\cK_{m^*}}^\beta(\Omega)=W^{m}(\Omega)$ with equivalent norms, where $\beta = m/m^*$. 

\section{Conditions in Section \ref{sec:theory_deter}}\label{app:condofdeter}

In this section, we list the conditions used in Section \ref{sec:theory_deter} with discussion. Since $\cK_{m}(\cdot,\cdot)$ is a symmetric and positive definite kernel function, Mercer's theorem implies that
\begin{align}\label{eq:eigenckm}
    \cK_{m}(\bs,\bt)=& \sum_{j=1}^\infty \gamma_{j} \phi_j(\bs)\phi_j(\bt), \forall \bs,\bt\in\Omega,
\end{align}
where $\gamma_{j}$ and $\phi_j$ are eigenvalues and eigenfunctions of $\cK_{m}(\cdot,\cdot)$, respectively, and the convergence is absolute and uniform (provided that $\cK_{m}(\bs,\bt)$ exists). 
We write $\btheta=(\theta_1,\ldots,\theta_q)^{\mathrm{T}}$.

\begin{enumerate}
	\item[(C1)] The random variables $\bx_j$ and $\varepsilon_j$ in \eqref{yp} are independent; $\bx_j$'s are i.i.d. uniformly distributed on $\Omega$; $\varepsilon_j$'s are i.i.d. random variables which are sub-Gaussian \citep{geer2000empirical_ec}, i.e., satisfying
	\begin{align*}
	C^2 (\mathbb{E} \exp\left(|\varepsilon_j|^2/C^2\right)-1)\leq C', \quad j=1,...,n.
	\end{align*}
	\item[(C2)] $\btheta^*$ is the unique solution to \eqref{thetastar}, and is an interior point of $\bTheta$.
	\item[(C3)] The matrix 
	$$\Vb\coloneqq-\sum_{j=1}^\infty \frac{1}{{\gamma_j}}\frac{\partial^{2}}{\partial \btheta \partial \btheta^{\mathrm{T}}}\left(\int_{\Omega} (f_p(\bz)-f_s(\bz,\btheta^*))\phi_j(\bz) {\rm d} \bz\right)^2$$ is invertible, where $\gamma_{j}$ and $\phi_j(\cdot)$ are as in \eqref{eq:eigenckm}.
	\item[(C4)] There exists $m_1 > d/2$ such that $m_1\geq m$ (if $m \leq d/2$, then clearly $m_1\geq m$ holds) and the RKHS $\cN_\Phi(\Omega)$ in \eqref{eq:krr} is equivalent to $W^{m_1}(\Omega)$. Furthermore, $\|f_p(\cdot)\|_{W^{m_1}(\Omega)}<\infty$, $\lambda=o_{\PP}(n^{-1/2-\epsilon_0})$ for some constant $\epsilon_0>0$, $\lambda^{-1} = O_\PP(n^{\frac{2m_1}{2m_1+d}})$, and $\left\|f_p(\cdot)-\hat f_p(\cdot)\right\|_{W^m(\Omega)}=o_{\PP}(1)$. 
	\item[(C5)] There exists  $\alpha_0 > 2m+m_1-\kappa$ with $0 \leq \kappa < \min(\frac{2m_1^2-m_1d}{2m_1+d},\frac{m_1}{1+4\epsilon_0})$ such that for all $j,k=1,...,q$ and $t=1,2,...$,
	\begin{align*}
	    & \sup _{\btheta \in \bTheta}\left\|f_s(\cdot, \btheta)\right\|_{W^{\alpha_0}(\Omega)}+\left\|\frac{\partial f_{s}(\cdot, \btheta)}{\partial \theta_{j}}\right\|_{W^{\alpha_0}(\Omega)} + \left\|\frac{\partial^{2}}{\partial \theta_j \partial \theta_k}f_s(\cdot,\btheta)\right\|_{W^{\alpha_0}(\Omega)} <\infty,\nonumber\\
	    & \frac{\partial^{2}}{\partial \theta_j \partial \theta_k}f_s(\cdot,\cdot)\in \cC^1(\Omega \times \bTheta), \phi_t(\cdot)\in W^{\max(\alpha_0,m_1)}(\Omega),
	\end{align*}
	where $m_1$ and $\epsilon_0$ are as in Condition (C4).
	\item[(C6)] $\sup _{\btheta \in \bTheta}\|f_s(\cdot,\btheta)-\hat f_s(\cdot,\btheta)\|_{W^{2m}(\Omega)}+\left\|\frac{\partial \hat f_s(\cdot,\btheta)}{\partial\theta_j}-\frac{\partial {f_s(\cdot,\btheta)}}{\partial\theta_j}\right\|_{W^{2m}(\Omega)}=o_{\PP}\left(n^{-1/2}\right)$, for all $j=1,...,q$.
\end{enumerate}

Conditions (C1)-(C3) are regularity conditions on the model. Similar conditions are also assumed in \cite{tuo&wu2_ec,xie2020bayesian_ec}. Because we use Sobolev norm instead of the $L_2$ norm used in \cite{tuo&wu2_ec}, there is an extra $\gamma_j^{-1}$ in the expression of the matrix $\Vb$ in Condition (C3). The matrix $\Vb$ reduces to the invertible matrix in Theorem 1 of \cite{tuo&wu2_ec} if $m=0$.

Condition (C4) requires that the RKHS used in kernel ridge regression \eqref{eq:krr} has a higher smoothness than $m$. By the Gagliardo–Nirenberg interpolation inequality for functions in Sobolev spaces \citep{leoni2017first_ec,BrezisMironescu19_ec}, the condition $\|f_p(\cdot)-\hat f_p(\cdot)\|_{W^m(\Omega)}=o_{\PP}(1)$ can be easily fulfilled as long as $m_1>m$ and $\|f_p(\cdot)-\hat f_p(\cdot)\|_{L_2(\Omega)}=o_{\PP}(1)$, while the later has been widely established in literature; see \cite{geer2000empirical_ec,fischer2020sobolev_ec} for example. Even if $m_1=m$, as long as $f_p(\cdot)$ has a smoothness higher than $m$, the condition $\|f_p(\cdot)-\hat f_p(\cdot)\|_{W^m(\Omega)}=o_{\PP}(1)$ can be ensured \citep{lin2017distributed_ec}. 

Condition (C5) imposes regularity conditions on the computer model, while Condition (C6) imposes conditions on the surrogate model $\hat f_s(\cdot,\cdot)$. In particular, the requirement $\phi_t(\cdot)\in W^{\max(\alpha_0,m_1)}(\Omega)$ can be fulfilled if $\cH_m(\Omega)$ is an $(m/\max(\alpha_0,m_1))$-th power of an RKHS that is equivalent to $W^{\max(\alpha_0,m_1)}(\Omega)$.
Moreover, Conditions (C5) and (C6) state that the computer model needs to be smoother than the Sobolev norm used in the Sobolev calibration, and the surrogate model needs to approximate the computer model well, in the sense that the Sobolev norm of the error is small. These conditions are reasonable considering the controllable cost of running codes. Also, these conditions reduce to the assumptions used in the $L_2$ calibration if $m=0$. Note that the assumption that the computer model lies in a function space with higher smoothness has also been used in Theorem 3.4 of \cite{tuo2020improved_ec}, which states a convergence rate of the estimated calibration parameter via KO calibration. By the Sobolev embedding theorem \citep{adams2003sobolev_ec}, Conditions (C5) and (C6) imply that for all $1\leq t<\infty$, $j,k=1,...,q$, $\phi_t(\cdot)\in W^{\alpha'}(\Omega)$ and 
	$$\sup _{\btheta \in \bTheta}\left\|f_s(\cdot, \btheta)\right\|_{W^{\alpha'}(\Omega)}+\left\|\frac{\partial f_{s}(\cdot, \btheta)}{\partial \theta_{j}}\right\|_{W^{\alpha'}(\Omega)}  + \left\|\frac{\partial^{2}}{\partial \theta_j \partial \theta_k}f_s(\cdot,\btheta)\right\|_{W^{\alpha'}(\Omega)} <\infty,$$
	for all $0\leq \alpha'\leq \alpha_0$, and 
$$\sup _{\btheta \in \bTheta}\|f_s(\cdot,\btheta)-\hat f_s(\cdot,\btheta)\|_{W^{2m'}(\Omega)}+\left\|\frac{\partial \hat f_s(\cdot,\btheta)}{\partial\theta_j}-\frac{\partial {f_s(\cdot,\btheta)}}{\partial\theta_j}\right\|_{W^{2m'}(\Omega)}=o_{\PP}\left(n^{-1/2}\right),$$ for all $j=1,...,q$, $0\leq m'\leq m$.

Condition (C5) is weaker than the conditions in \cite{tuo&wu2_ec,xie2020bayesian_ec}, even if at the case $m=0$. To see this, we set $m=0$. In \cite{tuo&wu2_ec,xie2020bayesian_ec}, it is required that $\alpha_0=m_1$, while Condition (C5) implies that $\alpha_0 > m_1-\kappa$ is sufficient. Since $\kappa$ can be larger than zero, it is possible that $\alpha_0 < m_1$. In other words, we consider the case where the computer model can even be \textit{rougher} than the physical model at some parameter space of $\btheta$, which has not been presented in the literature as far as we know.

\section{Conditions in Section \ref{sec:theory_gp}}\label{app:condofgp}

In this section, we list the conditions used in Section \ref{sec:theory_gp} with discussion. Without loss of generality, we assume that the variance of the Gaussian process satisfies $\sigma^2=1$, thus the correlation function equals the covariance function. 

\begin{enumerate}
    \item[(C1')] The random variables $\bx_j$ and $\varepsilon_j$ in \eqref{yp} are independent; $\bx_j$'s are i.i.d. uniformly distributed on $\Omega$; $\varepsilon_j$'s are i.i.d. normally distributed random variables.
	\item[(C2')] Conditions (C2)-(C3) hold with probability one for $f_p(\cdot)$.
	\item[(C3')] There exists $m_2 > \max(d, m+d/2)$ such that
	$$ c_1(1+\|\bomega\|_2^2)^{-m_2} \leq  \mathcal{F}(\cK)(\bomega)\leq c_2(1+\|\bomega\|_2^2)^{-m_2}, \forall \bomega\in\mathbb{R}^d,$$
	where $c_1$ and $c_2$ are positive constants, and $\mathcal{F}(\cK)(\bomega)$ is the Fourier transform of $\cK$. 
	\item[(C4')] Condition (C5) holds for $\alpha_0 = 2m + m_2$ and in particular we have $\phi_t(\cdot)\in W^{\alpha_0}(\Omega)$.
\end{enumerate}
Condition (C1') is a regularity condition on the model. The assumption that $\varepsilon_j$'s are normal is only for the ease of proof, and can be relaxed to the sub-Gaussian noise, with some additional mathematical treatment as in \cite{wang2021convergence_ec}. Condition (C2') is a technical assumption for the Gaussian process $f_p(\cdot)$ and the computer model $f_s(\cdot,\btheta)$. 
 Condition (C3') implies that the sample paths of $f_p(\cdot)$ lie in $W^{\tilde m}(\Omega)$ with probability one for some $\tilde m >d/2$ \citep{steinwart2019convergence_ec}. Condition (C4') is stronger than Condition (C5), which is because the Gaussian process has introduced more randomness than the deterministic function case, and it requires higher smoothness on the computer model.

\section{Technical Proofs of Propositions and Theorems}
 In this section, we provide technical proofs of propositions and theorems in the main text.

\subsection{Proof of Proposition \ref{prop:rkhsl2infty}}
Let $\cN_{\Psi_m}(\Omega)$ be the RKHS generated by $\Psi_m$ defined in \eqref{eq_matern}. Thus, as stated in Section \ref{app:introtoSoboRKHS} in the supplement, we have $\cN_{\Psi_m}(\Omega) = W^m(\Omega)$ with equivalent norms. 
By Mercer's Theorem, we have 
\begin{align*}
    \Psi_m(\bs,\bt)= & \sum_{j=1}^\infty \gamma_j \phi_j(\bs)\phi_j(\bt), \forall \bs,\bt\in\Omega,
\end{align*}
where $\phi_j$ and $\gamma_j$ are the eigenfunctions and eigenvalues of $\Psi_m(\cdot,\cdot)$, respectively, and the convergence is absolute and uniform. Furthermore, we have that $\gamma_j\rightarrow 0$ as $j\rightarrow \infty$ and $\|\phi_j\|_{L_2(\Omega)}=1$. Theorem 10.29 of \cite{wendland2004scattered_ec} implies that $\|\phi_k\|_{\cN_{\Psi_m}(\Omega)}^2=\sum_{j=1}^\infty \frac{\langle \phi_k(\cdot),\phi_j(\cdot)\rangle_{L_2(\Omega)}^2}{\gamma_j}=\frac{1}{\gamma_k}$. Also, the equivalence of 
$\|\cdot\|_{\cN_{\Psi_m}(\Omega)}$ and $\|\cdot\|_{W^{m}(\Omega)}$ implies that $c_1\|\phi_k\|_{\cN_{\Psi_m}(\Omega)}\leq \|\phi_k\|_{W^m(\Omega)} \leq c_{2}\|\phi_k\|_{\cN_{\Psi_m}(\Omega)}$. 
Thus, we can find $k$ such that $\gamma_k<\frac{c_{1}}{c_{2}}\frac{\gamma_1}{c}$. Let $\delta_2=\phi_1$ and $\delta_1=\phi_k$. It can be easily verified that $\|\delta_1\|_{L_2(\Omega)}=\|\delta_2\|_{L_2(\Omega)}=1$ and $\frac{\|\delta_1\|_{W^m(\Omega)}}{\|\delta_2\|_{W^m(\Omega)}}\geq \frac{c_{1}}{c_{2}}\frac{\|\delta_1\|_{\mathcal{N}_{\Psi_m}(\Omega)}}{\|\delta_2\|_{\mathcal{N}_{\Psi_m}(\Omega)}}=\frac{c_{1}}{c_{2}}\frac{\gamma_1}{\gamma_k}>c$. This finishes the proof.  \hfill\BlackBox
\subsection{Proof of Proposition \ref{prop:consist}}\label{app:pfpropconsist}

From the definition of $\hat\btheta_S$ and $\btheta^*_S$ in \eqref{eq:thetahat} and \eqref{thetastar}, it suffices to prove that $\|\hat f_p(\cdot)-\hat f_s(\cdot, \btheta)\|_{\cH_m(\Omega)}$ converges to $\left\|f_p(\cdot)-f_s(\cdot, \btheta)\right\|_{\cH_m(\Omega)}$ uniformly with respect to $\btheta$ in probability. This is ensured by
\begin{align*}
    &\left|\left\|\hat f_p(\cdot)-\hat f_s(\cdot, \btheta)\right\|_{\cH_m(\Omega)}-\left\|f_p(\cdot)-f_s(\cdot, \btheta)\right\|_{\cH_m(\Omega)}\right| \nonumber\\
    \le &\left\|\hat f_p(\cdot)-\hat f_s(\cdot, \btheta)-(f_p(\cdot)-f_s(\cdot, \btheta))\right\|_{\cH_m(\Omega)} \nonumber\\
    \le &\left\|\hat f_p(\cdot)-f_p(\cdot)\right\|_{\cH_m(\Omega)}+\left\|\hat f_s(\cdot, \btheta)-f_s(\cdot, \btheta)\right\|_{\cH_m(\Omega)}\nonumber\\
    \leq & C_1\left(\left\|\hat f_p(\cdot)-f_p(\cdot)\right\|_{W^m(\Omega)}+\left\|\hat f_s(\cdot, \btheta)-f_s(\cdot, \btheta)\right\|_{W^m(\Omega)}\right)\nonumber\\
    =&o_{\PP}(1),
\end{align*}
where the first and second inequalities are from the triangle inequality, the third inequality is because of the equivalence of $\|\cdot\|_{W^m(\Omega)}$ and $\|\cdot\|_{\cH_m(\Omega)}$, and the last equality is guaranteed by Conditions (C4) and (C6). This finishes the proof.

\subsection{Proof of Theorem \ref{thm:asympn}}\label{app:pfthmmain}

We first present several lemmas used in this proof. Lemma \ref{lem:innerinGeer} is Theorem 3.1 of \cite{van2014uniform_ec}. Lemma \ref{lemmaefsmall} is Lemma 8.4 of \cite{geer2000empirical_ec}. 
Lemma \ref{lem:fstar} is Lemma H.2 of \cite{wang2021convergence_ec}. Lemma \ref{lem:krr} is a direct result of Theorem 10.2 and Theorem 5.16 of \cite{geer2000empirical_ec}.
\begin{lemma}\label{lem:innerinGeer}
Let $\mathcal{F}$ and $\mathcal{G}$ be two function classes. Let $R_1\coloneqq\sup_{f\in\mathcal{F}}\|f\|_{L_2(\Omega)}, K_1\coloneqq\sup_{f\in\mathcal{F}}\|f\|_{L_\infty(\Omega)},R_2\coloneqq\sup_{g\in\mathcal{G}}\|g\|_{L_2(\Omega)}, K_2\coloneqq\sup_{g\in\mathcal{G}}\|g\|_{L_\infty(\Omega)}$. Suppose $R_1K_2\leq R_2K_1$. Consider values $t>4$ and $n$ such that
\begin{align}\label{eq:eqininnerlem}
    \left(\frac{2R_1J_\infty(K_1,\mathcal{F})+R_1K_1\sqrt{t}}{\sqrt{n}}+\frac{4J_\infty(K_1,\mathcal{F})^2+K_1^2t}{n}\right)\leq \frac{R_1^2}{C_1},\nonumber\\
    \left(\frac{2R_2J_\infty(K_2,\mathcal{G})+R_2K_2\sqrt{t}}{\sqrt{n}}+\frac{4J_\infty(K_2,\mathcal{G})^2+K_2^2t}{n}\right)\leq \frac{R_2^2}{C_1},
\end{align}
where 
\begin{align*}
{J}_{\infty}(z, \mathcal{F})\coloneqq C_{2} \inf _{\delta>0}\left[z \int_{\delta / 4}^{1} \sqrt{\mathcal{H}_{\infty}(u z / 2, \mathcal{F})} d u+\sqrt{n} \delta z\right],
\end{align*}
with $C_2$ as another constant. Denote the empirical measure of $f$ by $P_nf$ and its theoretical measure by $Pf$.
Then with probability at least $1-12\exp(-t)$,
\begin{align*}
\frac{1}{8C_1}\sup_{f\in\mathcal{F},g\in \mathcal{G}}\bigg|(P_n-P)fg\bigg|\leq \frac{R_1J_\infty(K_2,\mathcal{G})+R_2J_\infty(R_1K_2/R_2,\mathcal{F})+R_1K_2\sqrt{t}}{\sqrt{n}}+\frac{K_1K_2t}{n}.
\end{align*}
\end{lemma}

\begin{lemma}\label{lemmaefsmall}
	Let $\cK(\cdot,\cdot)$ be a kernel function such that 
	\begin{align*}
	    \mathcal{H}(\delta, \cB_{\cN_{\cK}(\Omega)}(1),\|\cdot\|_{L_\infty(\Omega)})\leq C_1\left(\frac{1}{\delta}\right)^{\frac{d}{m}}.
	\end{align*}
	Then for all $t > C$, with probability at least $1 - C_1\exp(-C_2 t^2)$,
	\begin{align*}
	\sup_{g\in \cN_{\cK}(\Omega)}\frac{|\langle \epsilon, g \rangle_n|}{\|g\|_n^{1 - \frac{d}{2m}}\|g\|_{\cN_{K}(\Omega)}^{\frac{d}{2m}}} \leq tn^{-\frac{1}{2}}.
	\end{align*}
\end{lemma}

\begin{lemma}\label{lem:fstar}
Let $m_1\geq m_0>d/2$, $\cK(\cdot,\cdot)$ be a positive definite kernel function on $\Omega$ such that $\cN_{\cK}(\Omega)$ is equivalent to $W^{m_1}(\Omega)$, and $f\in W^{m_0}(\Omega)$. Let $f^*$ be the solution to the optimization problem 
	\begin{align*}
	\min_{\tilde f \in \cN_{\cK}(\Omega)}\|f-\tilde f\|_{L_2(\Omega)}^2 + \lambda_0 \|\tilde f\|^2_{\cN_{\cK}(\Omega)}.
	\end{align*}
Then we have
\begin{align*}
	\|f-f^*\|_{L_2(\Omega)}^2 + \lambda_0 \|f^*\|^2_{\mathcal{N}_{\cK}(\Omega)} \leq C\lambda_0^{\frac{m_0}{m_1}}\|f\|_{W^{m_0}(\Omega)}^2.
	\end{align*}
\end{lemma}

\begin{lemma}\label{lem:krr}
Suppose $f\in W^{m}(\Omega)$, the reproducing kernel Hilbert space $\cN_{\Phi_1}(\Omega)$ generated by $\Phi_1(\cdot,\cdot)$ is equivalent to $W^{m}(\Omega)$, and $\lambda^{-1} = O_\PP(n^{\frac{2m}{2m+d}})$. Furthermore, assume Condition (C1) holds. Let $(\bx_j,y_j)$, $j=1,...,n$ be observations satisfying
\begin{align*}
    y_j = f(\bx_j) + \varepsilon_j.
\end{align*}
Then the estimator given by 
\begin{align*}
\hat f(\cdot)\coloneqq\underset{g \in \cN_{\Phi_1}(\Omega)}\argmin \frac{1}{n} \sum_{j=1}^{n}\left(y_{j} - g\left(\bx_{j}\right)\right)^{2}+\lambda\|g\|_{\cN_{\Phi_1}(\Omega)}^{2},
\end{align*}
satisfies
\begin{align*}
    \|\hat f - f\|_{L_2(\Omega)} = & O_{\PP}(\lambda^{\frac{1}{2}}\vee n^{-\frac{1}{2}}\lambda^{-\frac{d}{4m}}),\quad
    \|\hat f\|_{\cN_{\Phi_1}(\Omega)}= O_{\PP}(1\vee n^{-\frac{1}{2}}\lambda^{-\frac{2m + d}{4m}}).
\end{align*}
\end{lemma}

Before the proof of Theorem \ref{thm:asympn}, let us introduce some additional notation and results. Let $\alpha = \max(\alpha_0,m_1)$, where $\alpha_0$ is as in Condition (C5). Choose an RKHS that is equivalent to $W^{\alpha}(\Omega)$. Denote this RKHS as $\cN_\alpha(\Omega)$, and the corresponding kernel function as $\Psi_\alpha(\cdot,\cdot)$. Then Mercer's Theorem implies that it possesses an absolutely and uniformly convergent representations as
\begin{align*}
    \Psi_\alpha(\bs,\bt)= & \sum_{j=1}^\infty \gamma_{\Psi,j} \psi_j(\bs)\psi_j(\bt), \forall \bs,\bt\in\Omega,
\end{align*}
where $\gamma_{\Psi,j}$ and $\psi_j$ are eigenvalues and eigenfunctions of $\Psi_\alpha(\cdot,\cdot)$, respectively. Therefore, for $0\leq \beta \leq 1$, we can construct $\beta$-th powers of $\cN_\alpha(\Omega)$ as
\begin{align*}
    \cN_\alpha^\beta(\Omega) = \left\{f(\cdot)=\sum_{j=1}^\infty a_j\gamma_{\Psi,j}^{\beta/2}\psi_j(\cdot): \sum_{j=1}^\infty a_j^2 < \infty\right\},
\end{align*}
with inner product defined as
\begin{align}\label{eq:pfpbrkhsinner}
    \langle f,g \rangle_{\cN_\alpha^\beta(\Omega)} = \sum_{j=1}^\infty a_jb_j, \mbox{ for } f(\cdot)=\sum_{j=1}^\infty a_j\gamma_{\Psi,j}^{\beta/2}\psi_j(\cdot), g(\cdot)=\sum_{j=1}^\infty b_j\gamma_{\Psi,j}^{\beta/2}\psi_j(\cdot), f,g\in \cN_\alpha^\beta(\Omega).
\end{align}

In particular, \cite{steinwart2012mercer_ec} showed that $\{\gamma_{\Psi,j}^{\beta/2}\psi_j(\cdot)\}$ is an orthonormal basis of $\cN_\alpha^\beta(\Omega)$, and Proposition 4.2 of \cite{steinwart2012mercer_ec} (and the arguments thereafter) implies that $\cN_\alpha^\beta(\Omega)$ is independent of the choice of $\psi_j(\cdot)\in W^\alpha(\Omega)$, i.e., the orthonormal basis of $L_2(\Omega)$ (but needs to be in $W^\alpha(\Omega)$). By Condition (C5), $\phi_j(\cdot)\in W^\alpha(\Omega)$, thus  $\cN_\alpha^\beta(\Omega)$ can be also expressed as 
\begin{align*}
    \cN_\alpha^\beta(\Omega) = \left\{f(\cdot)=\sum_{j=1}^\infty a_j\gamma_{\Psi,j}^{\beta/2}\phi_j(\cdot): \sum_{j=1}^\infty a_j^2 < \infty\right\},
\end{align*}
with inner product defined as in \eqref{eq:pfpbrkhsinner}, by replacing $\psi_j(\cdot)$ with $\phi_j(\cdot)$. Furthermore, $\cN_\alpha^\beta(\Omega)$ is equivalent to the interpolation space of $W^{\alpha\beta}(\Omega) = [L_2(\Omega),W^{\alpha}(\Omega)]_{\beta,2}$. The space $\cN_\alpha^\beta(\Omega)$ will play an important role in the proofs of Theorems \ref{thm:asympn} and \ref{thm:gp+asympn}, as we will see later.

\textit{Proof of Theorem \ref{thm:asympn}.} By Theorem 10.29 of \cite{wendland2004scattered_ec}, for any $\btheta$, we have
\begin{align*}
    \|f_p(\cdot)-f_s(\cdot,\btheta)\|_{\cH_m(\Omega)}^2= & \sum_{j=1}^\infty \frac{\langle f_p(\cdot)-f_s(\cdot,\btheta),\phi_j(\cdot)\rangle_{L_2(\Omega)}^2}{\gamma_j},\nonumber\\
    \|\hat f_p(\cdot)-\hat f_s(\cdot,\btheta)\|_{\cH_m(\Omega)}^2=& \sum_{j=1}^\infty \frac{\langle\hat f_p(\cdot)-\hat f_s(\cdot,\btheta),\phi_j(\cdot)\rangle_{L_2(\Omega)}^2}{\gamma_j}.
\end{align*}

Since $\btheta_S^*$ minimizes \eqref{thetastar} and $\hat\btheta_S$ minimizes \eqref{eq:thetahat}, by Conditions (C2) and (C5) we have that for $t=1,...,q$,
\begin{align}\label{eq:fozerotrue} 
0 = & \frac{\partial}{\partial \theta_t}\left\|f_p(\cdot)-f_s\left(\cdot, \btheta\right)\right\|_{\cH_m(\Omega)}^{2}\bigg|_{\btheta=\btheta_S^*}\nonumber\\
= & -\sum_{j=1}^\infty \frac{2}{\gamma_j}\int_\Omega(f_p(\bz)-f_s(\bz,\btheta_S^*))\phi_j(\bz){\rm d}\bz\int_\Omega\frac{\partial f_s(\bz,\btheta_S^*)}{\partial\theta_t}\phi_j(\bz){\rm d}\bz\nonumber\\
= & -2\left\langle f_p(\cdot)-f_s(\cdot,\btheta^*_S),\frac{\partial f_s(\cdot,\btheta^*_S)}{\partial\theta_t}\right\rangle_{\cH_m(\Omega)},
\end{align}
and
\begin{align}\label{eq:fozeroest} 
0=&\frac{\partial}{\partial \theta_t}\left\|\hat f_p(\cdot)-\hat f_s\left(\cdot, \btheta\right)\right\|_{\cH_m(\Omega)}^{2}\bigg|_{\btheta=\hat\btheta_S}
\nonumber\\
=&-\sum_{j=1}^\infty \frac{2}{\gamma_j}\int_\Omega(\hat f_p(\bz)-\hat f_s(\bz,\hat\btheta_S))\phi_j(\bz){\rm d}\bz\int_\Omega\frac{\partial \hat f_s(\bz,\hat\btheta_S)}{\partial\theta_t}\phi_j(\bz){\rm d}\bz
\nonumber\\
=&-\sum_{j=1}^\infty\frac{2}{\gamma_j}\int_\Omega(\hat f_p(\bz)-{f_s}(\bz,\hat\btheta_S))\phi_j(\bz){\rm d}\bz\int_\Omega\frac{\partial {f_s}(\bz,\hat\btheta_S)}{\partial\theta_t}\phi_j(\bz){\rm d}\bz+I_{1}+I_{2},
\end{align}
where
\begin{align*}
    I_{1}=&-\sum_{j=1}^\infty\frac{2}{\gamma_j}\int_\Omega(f_s(\bz,\hat\btheta_S)-\hat f_s(\bz,\hat\btheta_S))\phi_j(\bz){\rm d}\bz\int_\Omega \frac{\partial \hat f_s(\bz,\hat\btheta_S)}{\partial\theta_t}\phi_j(\bz){\rm d}\bz,\\
    I_{2}=&-\sum_{j=1}^\infty\frac{2}{\gamma_j}\int_\Omega(\hat f_p(\bz)-f_s(\bz,\hat\btheta_S))\phi_j(\bz){\rm d}\bz\int_\Omega \left(\frac{\partial \hat f_s(\bz,\hat\btheta_S)}{\partial\theta_t}-\frac{\partial {f_s(\bz,\hat\btheta_S)}}{\partial\theta_t}\right)\phi_j(\bz){\rm d}\bz.
\end{align*}

Consider the $\beta$-th power of $\cN_\alpha(\Omega)$ with $\alpha\beta=m$. Then $\cN_\alpha^{m/\alpha}(\Omega)$ is equivalent to $W^m(\Omega)$, thus is also equivalent to $\cH_m(\Omega)$. By the equivalence of $\cN_\alpha^{m/\alpha}(\Omega)$ and $\cH_m(\Omega)$, the eigenvalues of two spaces satisfy $C_1\gamma_j \leq \gamma_{\Psi,j}^{m/\alpha}\leq C_2 \gamma_j$.

The first term $I_1$ can be bounded by
\begin{align}\label{eq:asypf1I1}
    |I_1| \leq & 2\left(\sum_{j=1}^\infty \left(\int_\Omega(f_s(z,\hat\btheta_S)-\hat f_s(z,\hat\btheta_S))\phi_j(\bz){\rm d}\bz\right)^2\right)^{1/2}\nonumber\\
    \times&\left(\sum_{j=1}^\infty  \frac{1}{\gamma_j^2}\left(\int_\Omega \frac{\partial \hat f_s(z,\hat\btheta_S)}{\partial\theta_t}\phi_j(\bz){\rm d}\bz\right)^2\right)^{1/2}\nonumber\\
    \leq & 2C_3\|f_s(\cdot,\hat\btheta_S)-\hat f_s(\cdot,\hat\btheta_S)\|_{L_2(\Omega)}\left(\sum_{j=1}^\infty  \frac{1}{\gamma_{\Psi,j}^{2m/\alpha}}\left(\int_\Omega \frac{\partial \hat f_s(z,\hat\btheta_S)}{\partial\theta_t}\phi_j(\bz){\rm d}\bz\right)^2\right)^{1/2}\nonumber\\
    = & o_{\PP}(n^{-1/2})\left\|\frac{\partial \hat f_s(\cdot,\hat\btheta_S)}{\partial\theta_t}\right\|_{\cN_\alpha^{2m/\alpha}(\Omega)},
\end{align}
where the first inequality is because of the Cauchy-Schwarz inequality, the second inequality is because $\phi_j(\cdot)$ is an orthogonal basis of $L_2(\Omega)$ and our construction of $\cN_\alpha^{2m/\alpha}(\Omega)$, and the last equality is by Condition (C6).

Because of the equivalence of $\cN_\alpha^{2m/\alpha}(\Omega)$ and $W^{2m}(\Omega)$, together with Conditions (C5) and (C6), we have
\begin{align}\label{eq:asypf1I1par1}
    \left\|\frac{\partial \hat f_s(\cdot,\hat\btheta_S)}{\partial\theta_t}\right\|_{\cN_\alpha^{2m/\alpha}(\Omega)}\leq & C_4 \left\|\frac{\partial \hat f_s(\cdot,\hat\btheta_S)}{\partial\theta_t}\right\|_{W^{2m}(\Omega)}\nonumber\\
    \leq & C_4 \left(\left\|\frac{\partial f_s(\cdot,\hat\btheta_S)}{\partial\theta_t}\right\|_{W^{2m}(\Omega)} + \left\|\frac{\partial \hat f_s(\cdot,\hat\btheta_S)}{\partial\theta_t}-\frac{\partial f_s(\cdot,\hat\btheta_S)}{\partial\theta_t}\right\|_{W^{2m}(\Omega)}\right)\nonumber\\
    \leq & C_5,
\end{align}
where the second inequality is because of the triangle inequality, and the third inequality is because of Conditions (C5) and (C6). Plugging \eqref{eq:asypf1I1par1} in \eqref{eq:asypf1I1} leads to 
\begin{align}\label{eq:asypf1I1small}
    I_1= o_{\PP}(n^{-1/2}).
\end{align}

Similarly, the second term $I_2$ can be bounded by
\begin{align}\label{eq:asypf1I2}
    |I_2| \leq & 2\left(\sum_{j=1}^\infty  \frac{1}{\gamma_j^2}\left(\int_\Omega \left(\frac{\partial \hat f_s(\bz,\hat\btheta_S)}{\partial\theta_t}-\frac{\partial {f_s(\bz,\hat\btheta_S)}}{\partial\theta_t}\right)\phi_j(\bz){\rm d}\bz\right)^2\right)^{1/2} \nonumber\\
    \times &\left(\sum_{j=1}^\infty \left(\int_\Omega(\hat f_p(\bz)-f_s(\bz,\hat\btheta_S))\phi_j(\bz){\rm d}\bz\right)^2\right)^{1/2}\nonumber\\
    \leq & 2C_6\left\|\frac{\partial \hat f_s(\cdot,\hat\btheta_S)}{\partial\theta_t}-\frac{\partial {f_s(\cdot,\hat\btheta_S)}}{\partial\theta_t}\right\|_{\cN_\alpha^{2m/\alpha}(\Omega)}\left\|\hat f_p(\cdot)-f_s(\cdot,\hat\btheta_S)\right\|_{L_2(\Omega)}\nonumber\\
    \leq & 2C_7\left\|\frac{\partial \hat f_s(\cdot,\hat\btheta_S)}{\partial\theta_t}-\frac{\partial {f_s(\cdot,\hat\btheta_S)}}{\partial\theta_t}\right\|_{W^{2m}(\Omega)}\left(\left\|\hat f_p(\cdot)\right\|_{L_2(\Omega)}+\left\|f_s(\cdot,\hat\btheta_S)\right\|_{L_2(\Omega)}\right)\nonumber\\
    = & o_{\PP}(n^{-1/2})\left(\left\|\hat f_p(\cdot)\right\|_{L_2(\Omega)}+\left\|f_s(\cdot,\hat\btheta_S)\right\|_{L_2(\Omega)}\right),
\end{align}
where the first inequality is because of the Cauchy-Schwarz inequality, the second inequality is because $\phi_j(\cdot)$ is an orthogonal basis of $L_2(\Omega)$ and because of our construction of $\cN_\alpha^{2m/\alpha}(\Omega)$, the third inequality is by the triangle inequality and the equivalence of $\cN_\alpha^{2m/\alpha}(\Omega)$ and $W^{2m}(\Omega)$, and the last equality is by Condition (C6). 

Condition (C4) implies that
\begin{align*}
    \left\|\hat f_p(\cdot)\right\|_{L_2(\Omega)} \leq \left\|f_p(\cdot)\right\|_{L_2(\Omega)} + \left\|\hat f_p(\cdot)-f_p(\cdot)\right\|_{L_2(\Omega)} \leq C_8+o_{\PP}(1),
\end{align*}
which, together with \eqref{eq:asypf1I2} and Condition (C5), implies that
\begin{align}\label{eq:asypf1I2small}
    I_2= o_{\PP}(n^{-1/2}).
\end{align}
Combining \eqref{eq:asypf1I1small} and \eqref{eq:asypf1I2small} leads to  
\begin{align}\label{eq:asypf1I1I2small}
    I_1+I_2 = o_{\PP}(n^{-1/2}).
\end{align}

Plugging \eqref{eq:asypf1I1I2small} in \eqref{eq:fozeroest} leads to
\begin{align}\label{eq:thm1pf1inners1}
    o_{\PP}(n^{-1/2}) = &  \sum_{j=1}^\infty\frac{2}{\gamma_j}\int_\Omega(\hat f_p(\bz)-{f_s}(\bz,\hat\btheta_S))\phi_j(\bz){\rm d}\bz\int_\Omega\frac{\partial {f_s}(\bz,\hat\btheta_S)}{\partial\theta_t}\phi_j(\bz){\rm d}\bz\nonumber\\
    = & \sum_{j=1}^\infty\frac{2}{\gamma_j}\int_\Omega({f_p}(\bz)-{f_s}(\bz,\hat\btheta_S))\phi_j(\bz){\rm d}\bz\int_\Omega\frac{\partial {f_s}(\bz,\hat\btheta_S)}{\partial\theta_t}\phi_j(\bz){\rm d}\bz\nonumber\\
    & + \sum_{j=1}^\infty\frac{2}{\gamma_j}\int_\Omega(\hat f_p(\bz)-{f_p}(\bz))\phi_j(\bz){\rm d}\bz\int_\Omega\frac{\partial {f_s}(\bz,\hat\btheta_S)}{\partial\theta_t}\phi_j(\bz){\rm d}\bz\nonumber\\
    = & J_1 + J_2,
\end{align}
where the change of summation order can be checked by verifying the terms in $J_1$ and $J_2$ are absolutely summable. If $J_1$ and $J_2$ are absolutely summable, then Riemann rearrangement theorem justifies the change of summation order. In fact, by Conditions (C4), (C5), and (C6), we can see that $J_1$ and $J_2$ are absolutely summable. Take $J_1$ as an example. To verify the absolute convergence, we need to bound 
\begin{align*}
    \sum_{j=1}^\infty\frac{2}{\gamma_j}\left|\int_\Omega({f_p}(\bz)-{f_s}(\bz,\hat{\btheta}_S))\phi_j(\bz){\rm d}\bz\int_\Omega\frac{\partial {f_s}(\bz,\hat{\btheta}_S)}{\partial\theta_t}\phi_j(\bz){\rm d}\bz\right|,
\end{align*}
which can be done by 
\begin{align*}
    & \sum_{j=1}^\infty\frac{2}{\gamma_j}\left|\int_\Omega({f_p}(\bz)-{f_s}(\bz,\hat\btheta_S))\phi_j(\bz){\rm d}\bz\int_\Omega\frac{\partial {f_s}(\bz,\hat\btheta_S)}{\partial\theta_t}\phi_j(\bz){\rm d}\bz\right|\nonumber\\
    \leq & 2\left(\sum_{j=1}^\infty \left(\int_\Omega({f_p}(\bz)-{f_s}(\bz,\hat\btheta_S))\phi_j(\bz){\rm d}\bz\right)^2\right)^{1/2}\left(\sum_{j=1}^\infty \frac{1}{\gamma_j^2}\left(\int_\Omega\frac{\partial {f_s}(\bz,\hat\btheta_S)}{\partial\theta_t}\phi_j(\bz){\rm d}\bz\right)^2\right)^{1/2}\nonumber\\
    \leq & 2C_9\|f_p(\cdot)-{f_s}(\cdot,\hat\btheta_S)\|_{L_2(\Omega)}\left\| \frac{\partial{f_s}(\bz,\hat\btheta_S)}{\partial\theta_t}\right\|_{W^{2m}(\Omega)}\nonumber\\
    < & \infty,
\end{align*}
where the first inequality is by the Cauchy-Schwarz inequality, the second inequality can be obtained similarly as in bounding $I_1$ in \eqref{eq:asypf1I1}, and the last inequality is by Conditions (C4) and (C5). In the rest of the Supplement, we do not provide justification of the change of summation order for the briefness of the proofs. One can check the change of summation order by a similar approach as above by verifying the absolute convergence.

Consider $J_1$. Define 
\begin{align*}
    v_t(\btheta) & = 2\left\langle f_p(\cdot)-f_s(\cdot,\btheta),\frac{\partial {f_s}(\cdot,\btheta)}{\partial\theta_t}
   \right\rangle_{\cH_m(\Omega)}.
\end{align*}
Then we have $J_1=v_t(\hat\btheta_S)$, and by \eqref{eq:fozerotrue}, $v_t(\btheta_S^*)=0$. 
By applying Taylor expansion to $v_t(\btheta)$ at point $\btheta_S^*$, 
we obtain
\begin{align*}
    J_1 & =v_t(\btheta_S^*) + \frac{\partial v_t(\tilde\btheta)}{\partial \btheta^{\mathrm{T}}}(\hat\btheta_S-\btheta_S^*) = \frac{\partial v_t(\tilde\btheta)}{\partial \btheta^{\mathrm{T}}}(\hat\btheta_S-\btheta_S^*),
\end{align*}
which, together with \eqref{eq:thm1pf1inners1}, leads to
\begin{equation}\label{eq:thm1pf1J1J2sum2}
    \frac{\partial v_t(\tilde\btheta)}{\partial \btheta^{\mathrm{T}}}(\hat\btheta_S-\btheta_S^*) + J_2=o_{\PP}(n^{-1/2}),
\end{equation}
where $\tilde\btheta$ lies between $\btheta_S^*$ and $\hat\btheta_S$. 

By the consistency of $\hat\btheta_S$ as in Proposition \ref{prop:consist}, we have $\tilde{\btheta} \stackrel{p}{\rightarrow} 
\btheta^{*}$. This implies that 
\begin{align}\label{eqv}
\frac{\partial v(\tilde\btheta)}{\partial \btheta^{\mathrm{T}}} = & -\sum_{j=1}^\infty \frac{1}{{\gamma_j}}\frac{\partial^{2}}{\partial \btheta \partial \btheta^{\mathrm{T}}}\left(\int_{\Omega} (f_p(\bz)-f_s(\bz,\tilde\btheta))\phi_j(\bz) {\rm d}\bz\right)^2 \nonumber \\
\stackrel{\PP}{\rightarrow} & -\sum_{j=1}^\infty \frac{1}{{\gamma_j}}\frac{\partial^{2}}{\partial \btheta \partial \btheta^{\mathrm{T}}}\left(\int_{\Omega} (f_p(\bz)-f_s(\bz,\btheta_S^*))\phi_j(\bz) {\rm d}\bz\right)^2\nonumber \\
= & \Vb,
\end{align}
where $v(\btheta)=2\left\langle f_p(\cdot)-f_s(\cdot,\btheta),\frac{\partial {f_s}(\cdot,\btheta)}{\partial\btheta}
   \right\rangle_{\cH_m(\Omega)}=\left(v_1(\btheta),...,v_q(\btheta)\right)^{\mathrm{T}}$ and $\Vb$ is as in Condition (C3). 

Consider $J_2$. Define 
\begin{align*}
    A_{1t}(\cdot,\btheta) = \sum_{j=1}^\infty\frac{1}{\gamma_j}\int_\Omega\frac{\partial {f_s}(\bz,{\btheta})}{\partial\theta_t}\phi_j(\bz){\rm d}\bz \cdot \phi_j(\cdot).
\end{align*}
Theorem 10.29 of \cite{wendland2004scattered_ec} implies that
\begin{align}\label{eq:pf1th1A1t1}
   \sup_{\btheta\in \bTheta} \|A_{1t}(\cdot,\btheta)\|_{\cN_\alpha^{\beta_1}(\Omega)}^2 & = \sup_{\btheta\in \bTheta}\sum_{j=1}^\infty\frac{1}{\gamma_j^2\gamma_{\Psi,j}^{\beta_1}}\left(\int_\Omega\frac{\partial {f_s}(\bz,{\btheta})}{\partial\theta_t}\phi_j(\bz){\rm d}\bz\right)^2 \leq C_{10},
\end{align}
where $\beta_1=\frac{1}{\alpha}(\frac{2m_1d}{2m_1+d} + \epsilon_1)$ with $\epsilon_1= \frac{2m_1^2-m_1d}{2m_1+d} - \kappa > 0$,
and the last inequality is by Condition (C5) and the equivalence of $\cN_\alpha^{m/\alpha}(\Omega)$ and $\cH_m(\Omega)$. Thus, $A_{1t}(\cdot,\btheta)\in W^{\alpha\beta_1}(\Omega)$ because of the equivalence of $\cN_\alpha^{\beta_1}(\Omega)$ and $W^{\alpha\beta_1}(\Omega)$. We also have
\begin{align*}
    m_1 - \kappa = \alpha\beta_1 > \frac{2m_1d}{2m_1+d} > \frac{d}{2}, \mbox{ and } 2m+\alpha\beta_1 < \alpha.
\end{align*}
Let
\begin{align*}
    \mathcal{L}_{n}(f)=\frac{1}{n} \sum_{j=1}^{n}\left(y_j-f(\bx_j)\right)^{2}+\lambda_n\|f\|_{\cN_{\Phi}(\Omega)}^{2}.
\end{align*}
Since $\hat f_p(\cdot)$ minimizes $\mathcal{L}_n$ over $\cN_{\Phi}(\Omega)$, we have
\begin{align}\label{eq:thm1pf1Ln1}
    0 = & \left.\frac{\partial}{\partial \eta} \mathcal{L}_n\left(\hat f_p(\cdot) + \eta \hat A_{1t}(\cdot,\hat\btheta_S)\right)\right|_{\eta=0} \nonumber\\
    = & \frac{2}{n} \sum_{k=1}^{n}\left(\hat f_p\left(\bx_k\right)-y_k\right) \hat A_{1t}(\bx_k,\hat\btheta_S)+2 \lambda\left\langle\hat f_p(\cdot), \hat A_{1t}(\cdot,\hat\btheta_S)\right\rangle_{\cN_{\Phi}(\Omega)}\nonumber \\
    = & \frac{2}{n} \sum_{k=1}^{n}\left(\hat f_p\left(\bx_k\right)-f_p\left(\bx_k\right)\right) A_{1t}(\bx_k,\hat\btheta_S)\nonumber \\
    + &\frac{2}{n} \sum_{k=1}^{n}\left(\hat f_p\left(\bx_k\right)-f_p\left(\bx_k\right)\right)\left(\hat A_{1t}(\bx_k,\hat\btheta_S) - A_{1t}(\bx_k,\hat\btheta_S)\right)\nonumber\\
    - & \frac{2}{n} \sum_{k=1}^{n} \varepsilon_k A_{1t}(\bx_k,\hat\btheta_S) - \frac{2}{n} \sum_{k=1}^{n} \varepsilon_k \left(\hat A_{1t}(\bx_k,\hat\btheta_S) -A_{1t}(\bx_k,\hat\btheta_S)\right)\nonumber\\
    + &2 \lambda\left\langle\hat f_p(\cdot), \hat A_{1t}(\cdot,\hat\btheta_S)\right\rangle_{\cN_{\Phi}(\Omega)}\nonumber\\
    \coloneqq & Q_{1} + Q_{2} - Q_{3} - Q_{4} + Q_{5},
\end{align}
where $\hat A_{1t}(\cdot,\hat\btheta_S)\in \cN_\Phi(\Omega)$ is a function to be chosen later.

Since $\alpha\beta_1 < m_1$, $A_{1t}(\cdot,\btheta)$ may not in $W^{m_1}(\Omega)$. Therefore, the function $\hat A_{1t}(\cdot,\btheta)$ should be a good approximation of $A_{1t}(\cdot,\btheta)$ and is in the RKHS $\cN_\Phi(\Omega)$. We select $\hat A_{1t}(\cdot,\btheta)$ as the solution to the optimization problem
\begin{align*}
	\min_{\tilde f(\cdot) \in \cN_{\Phi}(\Omega)}\|A_{1t}(\cdot,\btheta)-\tilde f(\cdot)\|_{L_2(\Omega)}^2 + \lambda_0 \|\tilde f(\cdot)\|^2_{\cN_{\Phi}(\Omega)},
\end{align*}
where $\lambda_0=C_{11}n^{-\delta}$ with $\delta =\frac{1}{2}\left(\frac{2m_1\epsilon_0}{m_1-\alpha\beta_1} + \frac{m_1}{2\alpha\beta_1}\right) > 0$. Since $m_1 > \alpha\beta_1$, we have
\begin{align*}
    & \frac{2m_1\epsilon_0}{m_1-\alpha\beta_1} \geq \frac{m_1}{2\alpha\beta_1} 
    \Leftrightarrow (1+4\epsilon_0)\alpha\beta_1\geq m_1\nonumber\\
    \Leftrightarrow & (1+4\epsilon_0)(m_1-\kappa) \geq m_1
    \Leftrightarrow \kappa \leq \frac{m_1}{1+4\epsilon_0},
\end{align*}
which holds because of Condition (C5). Thus, $\delta \in [\frac{m_1}{2\alpha\beta_1} , \frac{2m_1\epsilon_0}{m_1-\alpha\beta_1}]$.

Lemma \ref{lem:fstar} and \eqref{eq:pf1th1A1t1} imply that
\begin{align}\label{eq:A1tthetaB1}
    \sup_{\btheta\in\bTheta}\|A_{1t}(\cdot,\btheta)-\hat A_{1t}(\cdot,\btheta)\|_{L_2(\Omega)}^2 \leq C_{12}\lambda_0^{\frac{\alpha\beta_1}{m_1}}, \sup_{\btheta\in\bTheta} \|\hat A_{1t}(\cdot,\btheta)\|^2_{\cN_{\Phi}(\Omega)} \leq C_{13}\lambda_0^{\frac{\alpha\beta_1-m_1}{m_1}}.
\end{align}
In the following, we consider each term in \eqref{eq:thm1pf1Ln1}.

(i) Consider $Q_1$. Define $h_{1,\btheta}(\cdot)\coloneqq (\hat f_p(\cdot)-f_p(\cdot)) A_{1t}(\cdot,\btheta)$. Let $\cG_1=\{g(\cdot):g(\cdot)=A_{1t}(\cdot,\btheta),\btheta\in\bTheta\}$, and $\cF_1=\{f(\cdot)\in \cN_{\Phi}(\Omega): f(\cdot) = \hat f_p(\cdot)-f_p(\cdot)\}$. We will use Lemma \ref{lem:innerinGeer} to bound the difference of the empirical inner product and the $L_2$ inner product between functions  $\hat f_p(\cdot)-f_p(\cdot)$ and $A_{1t}(\cdot,\btheta)$. In order to do so, we first consider $\cG_1$. Since $A_{1t}(\cdot,\btheta)\in W^{\alpha\beta_1}(\Omega)$ and $\alpha\beta_1>d/2$, \eqref{eq:pf1th1A1t1} implies that $\cG_1\subset \cB_{W^{\alpha\beta_1}(\Omega)}(C_{10})$, where we also use the equivalence of $W^{\alpha\beta_1}(\Omega)$ and $\cN_\alpha^{\beta_1}(\Omega)$. Thus, the entropy number of $\cG_1$ can be bounded by \citep{edmunds2008function_ec}
\begin{align}\label{eq:pf1th1entroG1b1}
    H(\delta,\cG_1,\|\cdot\|_{L_\infty(\Omega)})\leq C_{14}\bigg(\frac{1}{\delta}\bigg)^{\frac{d}{\alpha\beta_1}}.
\end{align}
Furthermore, the Gagliardo--Nirenberg interpolation inequality for functions in Sobolev spaces (abbreviated as the interpolation inequality in the rest of the Supplement) \citep{BrezisMironescu19_ec} implies that
\begin{align}\label{eq:pf1th1ginG1b1}
    \|g\|_{L_\infty(\Omega)}\leq & C_{15}\|g\|_{L_2(\Omega)}^{1-\frac{d}{2\alpha\beta_1}}\|g\|_{W^{\alpha\beta_1}(\Omega)}^{\frac{d}{2\alpha\beta_1}}\leq C_{16}\|g\|_{W^{\alpha\beta_1}(\Omega)}^{\frac{d}{2\alpha\beta_1}} \leq C_{17}, \forall g\in \cG_1,
\end{align}
where the second inequality is because of the Sobolev embedding theorem, and the last inequality is because $\cG_1\subset \cB_{W^{\alpha\beta_1}(\Omega)}(C_{10})$. Therefore, \eqref{eq:pf1th1ginG1b1} implies that $R_2=\sup_{g\in \cG_1}\|g(\cdot)\|_{L_2(\Omega)}$ and $K_2=\sup_{g\in \cG_1}\|g(\cdot)\|_{L_\infty(\Omega)}$ satisfy $R_2\leq K_2\leq C_{17}$, where we recall Vol$(\Omega)=1$ through the proof. Thus, we can choose $R_2=K_2=C_{17}$ in the following calculations without changing the result. In order to apply Lemma \ref{lem:innerinGeer}, we need to compute ${J}_{\infty}(K_2, \mathcal{G}_1)$. 

Direct computation shows that
\begin{align}\label{eq:pf1th1JinftyG1}
{J}_{\infty}(K_2, \mathcal{G}_1) = & C_{18} \inf _{\delta>0}\left[K_2 \int_{\delta/4}^{1} \sqrt{\mathcal{H}_{\infty}(u K_2 / 2, \mathcal{F})} {\rm d} u+\sqrt{n} \delta K_2\right]\nonumber\\
\leq & C_{18} \inf _{\delta>0}\left[C_{19}K_2 \int_{\delta/4}^{1} (uK_2)^{-\frac{d}{2\alpha\beta_1}} {\rm d} u+\sqrt{n} \delta K_2\right]\nonumber\\
\leq & C_{20} K_2^{1-\frac{d}{2\alpha\beta_1}}
\leq C_{21},
\end{align}
where the first inequality is by \eqref{eq:pf1th1entroG1b1}, the second inequality is because $2\alpha\beta_1>d$ and $\int_0^1 u^{-\frac{d}{2\alpha\beta_1}}{\rm d}u$ is finite, and the last inequality is because we have chosen $K_2=C_{17}$.

Next, we turn to $\cF_1$. Let $R_1=\sup_{f\in \cF_1}\|f(\cdot)\|_{L_2(\Omega)}$ and $K_1=\sup_{f\in \cF_1}\|f(\cdot)\|_{L_\infty(\Omega)}$. Clearly $R_1K_2\leq K_1R_2$ because we have chosen $R_2=K_2=C_{17}$. Lemma \ref{lem:krr} implies that $\|\hat f_p(\cdot)\|_{\cN_{\Phi}(\Omega)} = O_{\PP}(1\vee n^{-1/2}\lambda^{-\frac{2m_1+d}{4m_1}})$, thus $\|\hat f_p(\cdot)\|_{W^{m_1}(\Omega)} = O_{\PP}(1\vee n^{-1/2}\lambda^{-\frac{2m_1+d}{4m_1}})$ by the equivalence of $\|\cdot\|_{\cN_{\Phi}(\Omega)}$ and $\|\cdot\|_{W^{m_1}(\Omega)}$. Furthermore, $R_1=O_{\PP}(\lambda\vee n^{-1/2}\lambda^{-\frac{d}{4m_1}})$ by Lemma \ref{lem:krr}. By a similar computation as in \eqref{eq:pf1th1JinftyG1}, we obtain
\begin{align*}
    {J}_{\infty}(R_1, \mathcal{F_1}) \leq C_{20}R_1^{1-\frac{d}{2m_1}}, {J}_{\infty}(K_1, \mathcal{F_1}) \leq C_{20}K_1^{1-\frac{d}{2m_1}},
\end{align*}
where we note that Condition (C4) implies $\cF_1\subset \cN_{\Phi}(\Omega)$. By the interpolation inequality, 
\begin{align*}
    \|\hat f_p(\cdot)-f_p(\cdot)\|_{L_\infty(\Omega)} \leq C_{21}\|\hat f_p(\cdot)-f_p(\cdot)\|_{L_2(\Omega)}^{1-\frac{d}{2m_1}}\|\hat f_p(\cdot)-f_p(\cdot)\|_{\cN_\Phi(\Omega)}^{\frac{d}{2m_1}},
\end{align*}
which implies $K_1 \leq C_{21}R_1^{1-\frac{d}{2m_1}} (\|\hat f_p(\cdot)\|_{\cN_\Phi(\Omega)}^{\frac{d}{2m_1}}+\| f_p(\cdot)\|_{\cN_\Phi(\Omega)}^{\frac{d}{2m_1}})=o_{\PP}(1)$, where we also use Condition (C4). 
Thus, it can be checked that \eqref{eq:eqininnerlem} is satisfied for sufficient large $n$. Then we can apply Lemma \ref{lem:innerinGeer} to $(\hat f_p(\cdot)-f_p(\cdot)) A_{1t}(\cdot,\btheta)$ and obtain that 
\begin{align*}
    \sup_{\btheta\in \Theta}\bigg|(P_n-P)h_{1,\btheta}(\cdot)\bigg|=o_\PP(n^{-1/2}),
\end{align*}
where $P_{n} f\coloneqq\sum_{i=1}^{n} f\left(X_{i}\right) / n$ and $P f\coloneqq\sum_{i=1}^{n} \mathbb{E} f\left(X_{i}\right) / n$.
Therefore, we have 
\begin{align*}
    \left|Q_1-2\int_\Omega \left(\hat f_p\left(\bz\right)-f_p\left(\bz\right)\right) A_{1t}(\bz,\hat\btheta_S) {\rm d}\bz \right| & \leq 2\sup_{\btheta\in \bTheta} \left|\frac{1}{n}\sum_{j=1}^n h_{1,\btheta}(\bx_j)- \int_\Omega h_{1,\btheta}(\bz){\rm d}\bz\right| \\
    &= o_{\PP}(n^{-1/2}),
\end{align*}
which implies 
\begin{align}\label{eq:pf1th1Q1}
    Q_1 = 2\int_\Omega \left(\hat f_p\left(\bz\right)-f_p\left(\bz\right)\right) A_{1t}(\bz,\hat\btheta_S) {\rm d}\bz + o_{\PP}(n^{-1/2}).
\end{align}

(ii) Consider $Q_2$. It follows from a similar approach as in (i) that
\begin{align}\label{eq:pf1th1Q21}
    Q_2 = & 2\int_\Omega \left(\hat f_p\left(\bz\right)-f_p\left(\bz\right)\right) \left(A_{1t}(\bz,\hat\btheta_S) - \hat A_{1t}(\bz,\hat\btheta_S)\right){\rm d}\bz + o_{\PP}(n^{-1/2})\nonumber\\
    \leq & 2\|\hat f_p\left(\cdot\right)-f_p\left(\cdot\right)\|_{L_2(\Omega)}\|A_{1t}(\cdot,\hat\btheta_S) - \hat A_{1t}(\cdot,\hat\btheta_S)\|_{L_2(\Omega)} + o_{\PP}(n^{-1/2})\nonumber\\
    = & o_{\PP}(n^{-1/4})O_{\PP}(\lambda_0^{\frac{\alpha\beta_1}{2m_1}})+ o_{\PP}(n^{-1/2})\nonumber\\
     = & o_{\PP}(n^{-1/4})O_{\PP}(n^{-\frac{m_1}{2\alpha\beta_1}\frac{\alpha\beta_1}{2m_1}})+ o_{\PP}(n^{-1/2})\nonumber\\
    = & o_{\PP}(n^{-1/4})O_{\PP}(n^{-1/4}) + o_{\PP}(n^{-1/2}) =o_{\PP}(n^{-1/2}),
\end{align}
where the first inequality is by the Cauchy-Schwarz inequality, and the second equality is by Lemma \ref{lem:krr} and \eqref{eq:A1tthetaB1}.

(iii) Consider $Q_4$. Define function class $\mathcal{G}_2 = \left\{g(\cdot)\bigg|g(\cdot) = A_{1t}(\cdot,\btheta) - \hat A_{1t}(\cdot,\btheta),\btheta \in \bTheta\right\}$. Since $\alpha\beta_1<m_1$, $\mathcal{G}\subset \cN_\alpha^{\beta_1}(\Omega)= W^{\alpha\beta_1}(\Omega)$. Lemma \ref{lemmaefsmall} implies
\begin{align}\label{eq:pf1th1ineQ4small}
    & \langle \bvarepsilon, A_{1t}(\cdot,\btheta) - \hat A_{1t}(\cdot,\btheta) \rangle_n\nonumber\\
    = & O_{\PP}(n^{-1/2})\|A_{1t}(\cdot,\btheta) - \hat A_{1t}(\cdot,\btheta)\|_{L_2(\Omega)}^{1-\frac{d}{2\alpha\beta_1}}\|A_{1t}(\cdot,\btheta) - \hat A_{1t}(\cdot,\btheta)\|_{W^{\alpha\beta_1}(\Omega)}^{\frac{d}{2\alpha\beta_1}}, \forall \btheta \in \bTheta.
\end{align}

By the triangle inequality,
\begin{align}\label{eq:A1tl2small}
    \sup_{\btheta \in \bTheta}\|\hat A_{1t}(\cdot,\btheta)\|_{L_2(\Omega)}\leq \sup_{\btheta \in \bTheta}\| A_{1t}(\cdot,\btheta)\|_{L_2(\Omega)} + \sup_{\btheta \in \bTheta}\|A_{1t}(\cdot,\btheta) - \hat A_{1t}(\cdot,\btheta)\|_{L_2(\Omega)} \leq C_{22},
\end{align}
where the last inequality is by \eqref{eq:pf1th1A1t1} and \eqref{eq:A1tthetaB1}. The interpolation inequality and \eqref{eq:A1tthetaB1} imply that
\begin{align}\label{eq:A1tWab1}
    \|\hat A_{1t}(\cdot,\btheta)\|_{W^{\alpha\beta_1}(\Omega)} \leq C_{23}\|\hat A_{1t}(\cdot,\btheta)\|_{L_2(\Omega)}^{1-\frac{\alpha\beta_1}{m_1}}\|\hat A_{1t}(\cdot,\btheta)\|_{W^{m_1}(\Omega)}^{\frac{\alpha\beta_1}{m_1}}\leq C_{24} \lambda_0^{\frac{(\alpha\beta_1-m_1)\alpha\beta_1}{2m_1^2}},
\end{align}
where the last inequality is by \eqref{eq:A1tthetaB1} and \eqref{eq:A1tl2small}. Plugging \eqref{eq:A1tl2small} and \eqref{eq:A1tWab1} into \eqref{eq:pf1th1ineQ4small} leads to 
\begin{align}\label{eq:pf1th1ineQ4small2}
    & \langle \bvarepsilon, A_{1t}(\cdot,\btheta) - \hat A_{1t}(\cdot,\btheta) \rangle_n\nonumber\\
    = & O_{\PP}(n^{-1/2})\lambda_0^{\left(1-\frac{d}{2\alpha\beta_1}\right)\frac{\alpha\beta_1}{2m_1}}\left(\|A_{1t}(\cdot,\btheta)\|_{W^{\alpha\beta_1}(\Omega)} + \|\hat A_{1t}(\cdot,\btheta)\|_{W^{\alpha\beta_1}(\Omega)}\right)^{\frac{d}{2\alpha\beta_1}}\nonumber\\
    = & O_{\PP}(n^{-1/2})\lambda_0^{\left(1-\frac{d}{2\alpha\beta_1}\right)\frac{\alpha\beta_1}{2m_1}}\lambda_0^{\frac{d}{2\alpha\beta_1}\left(\frac{(\alpha\beta_1-m_1)\alpha\beta_1}{2m_1^2}\right)}\nonumber\\
    = & o_{\PP}(n^{-1/2}),\quad \forall \btheta \in \bTheta,
\end{align}
where the first equality is by \eqref{eq:pf1th1ineQ4small}, the second equality is by \eqref{eq:A1tthetaB1} and interpolation inequality, and the last equality is because $\lambda_0=o(1)$ and 
\begin{align*}
    \left(1-\frac{d}{2\alpha\beta_1}\right)\frac{\alpha\beta_1}{2m_1} + \frac{d}{2\alpha\beta_1}\left(\frac{(\alpha\beta_1-m_1)\alpha\beta_1}{2m_1^2}\right) = \frac{\alpha\beta_1(2m_1+d)-2m_1d}{4m_1^2} > 0,
\end{align*}
where we use $\alpha\beta_1 > \frac{2m_1d}{2m_1+d}$.
Therefore, \eqref{eq:pf1th1ineQ4small2} implies that 
\begin{align}\label{eq:pf1th1Q4} 
    |Q_4|\leq 2\sup_{\btheta\in \bTheta}|\langle \bvarepsilon, A_{1t}(\cdot,\btheta) - \hat A_{1t}(\cdot,\btheta) \rangle_n| = o_{\PP}(n^{-1/2}).
\end{align}

(iv) Consider $Q_5$. By the Cauchy-Schwarz inequality,  
\begin{align*} 
    |Q_5|\leq & 2 \lambda\left\|\hat f_p(\cdot)\right\|_{\cN_{\Phi}(\Omega)} \left\|\hat A_{1t}(\cdot,\hat\btheta_S)\right\|_{\cN_{\Phi}(\Omega)}.
\end{align*}
By \eqref{eq:A1tthetaB1} and Condition (C4), we have
\begin{align}\label{eq:pf1th1Q5} 
    |Q_5| = & o_{\PP}(n^{-\frac{1}{2} - \epsilon_0})O_{\PP}( \lambda_0^{\frac{\alpha\beta_1-m_1}{2m_1}})\nonumber\\
    = & o_{\PP}(n^{-\frac{1}{2} - \epsilon_0})O_{\PP}( n^{\delta\frac{m_1-\alpha\beta_1}{2m_1}})\nonumber\\
    = & o_{\PP}(n^{-\frac{1}{2} - \epsilon_0})O_{\PP}( n^{\frac{2m_1\epsilon_0}{m_1-\alpha\beta_1}\frac{m_1-\alpha\beta_1}{2m_1}}) = o_{\PP}(n^{-1/2}).
\end{align}
By \eqref{eq:pf1th1Q1}, \eqref{eq:pf1th1Q21}, \eqref{eq:pf1th1Q4}, and \eqref{eq:pf1th1Q5}, we can rearrange $J_2$ as
\begin{align}\label{eq:thm1pf1J2b2}
    J_2= & 2\int_\Omega \left(\hat f_p\left(\bz\right)-f_p\left(\bz\right)\right) A_{1t}(\bz,\hat\btheta_S) {\rm d}\bz  = \frac{2}{n} \sum_{k=1}^{n} \varepsilon_k A_{1t}(\bx_k,\hat\btheta_S) + o_{\PP}(n^{-1/2})\nonumber\\
    = & \frac{2}{n} \sum_{k=1}^{n} \varepsilon_k A_{1t}(\bx_k,\btheta_S^*) + \frac{2}{n} \sum_{k=1}^{n} \varepsilon_k \left(A_{1t}(\bx_k,\hat\btheta_S) - A_{1t}(\bx_k,\btheta_S^*) \right)+ o_{\PP}(n^{-1/2}).
\end{align}

Define function class $\mathcal{G}_3 = \left\{g(\cdot)\bigg|g(\cdot) = A_{1t}(\cdot,\btheta) - A_{1t}(\cdot,\btheta_S^*),\btheta \in \bTheta\right\}$. Thus $\cG_3\subset W^{\alpha\beta_1}(\Omega)$.

Lemma \ref{lemmaefsmall} implies
\begin{align*}
    & \langle \bvarepsilon, A_{1t}(\cdot,\btheta) - A_{1t}(\cdot,\btheta_S^*) \rangle_n\nonumber\\
    = & O_{\PP}(n^{-1/2})\|A_{1t}(\cdot,\btheta) - A_{1t}(\cdot,\btheta_S^*)\|_{L_2(\Omega)}^{1-\frac{d}{2\alpha\beta_1}}\|A_{1t}(\cdot,\btheta) - A_{1t}(\cdot,\btheta_S^*)\|_{W^{\alpha\beta_1}(\Omega)}^{\frac{d}{2\alpha\beta_1}}\nonumber\\
    = & O_{\PP}(n^{-1/2})\|A_{1t}(\cdot,\btheta) - A_{1t}(\cdot,\btheta_S^*)\|_{L_2(\Omega)}^{1-\frac{d}{2\alpha\beta_1}}\left(\|A_{1t}(\cdot,\btheta)\|_{W^{\alpha\beta_1}(\Omega)} + \|A_{1t}(\cdot,\btheta_S^*)\|_{W^{\alpha\beta_1}(\Omega)}\right)^{\frac{d}{2\alpha\beta_1}}\nonumber\\
    = & O_{\PP}(n^{-1/2})\|A_{1t}(\cdot,\btheta) - A_{1t}(\cdot,\btheta_S^*)\|_{L_2(\Omega)}^{1-\frac{d}{2\alpha\beta_1}}, \forall \btheta \in \bTheta,
\end{align*}
where the second equality is because of the triangle inequality, and the last equality is by \eqref{eq:pf1th1A1t1}. Therefore,
\begin{align}\label{eq:innersmall2}
    \left|\frac{2}{n} \sum_{k=1}^{n} \varepsilon_k \left(A_{1t}(\bx_k,\hat\btheta_S) - A_{1t}(\bx_k,\btheta_S^*) \right)\right| &= O_{\PP}(n^{-1/2})\|A_{1t}(\cdot,\hat\btheta_S) - A_{1t}(\cdot,\btheta_S^*)\|_{L_2(\Omega)}^{1-\frac{d}{2\alpha\beta_1}} \nonumber\\
    &= o_{\PP}(n^{-1/2}),
\end{align}
where the last inequality is because of the consistency of $\hat\btheta_S$.

By \eqref{eq:thm1pf1J1J2sum2},\eqref{eqv},\eqref{eq:thm1pf1J2b2} and \eqref{eq:innersmall2}, we obtain that
\begin{equation*}
    (\Vb+\Eb)(\hat\btheta_S-\btheta_S^*) + \frac{2}{n} \sum_{k=1}^{n} \varepsilon_k A_{1}(\bx_k) = o_{\PP}(n^{-1/2}).
\end{equation*}
where each element in $\Eb$ is $o_{\PP}(1)$ and 
\begin{align*}
    A_{1}(\cdot)= \sum_{j=1}^\infty\frac{1}{\gamma_j}\int_\Omega\frac{\partial {f_s}(\bz,{\btheta^*_S})}{\partial\btheta}\phi_j(\bz){\rm d}\bz \cdot \phi_j(\cdot)=\left(A_{11}(\cdot,\btheta^*_S),...,A_{1q}(\cdot,\btheta^*_S)\right)^{\mathrm{T}}.
\end{align*}
Since $\Eb$ converges to zero in probability, we have
\begin{equation*}
    \hat\btheta_S-\btheta_S^* = -2\Vb^{-1} \left(\frac{1}{n} \sum_{k=1}^{n} \varepsilon_k A_{1}(\bx_k,\btheta_S^*)\right) + o_{\PP}(n^{-1/2}),
\end{equation*}
which finishes the proof. \hfill\BlackBox
\subsection{Proof of Theorem \ref{thm:semieff}}\label{app:pfthmsemieff}

The proof of Theorem \ref{thm:semieff} is along the line of the proof of Theorem 2 in \cite{tuo&wu2014_ec}, while we modify the original proof such that it fits the Sobolev calibration context. 

\textit{Proof of Theorem \ref{thm:semieff}.} For the calibration problem given by \eqref{yp} and \eqref{thetastar}, consider the $q$-dimensional parametric model indexed by $\bgamma$:
\begin{align}\label{gamma}
    f_{p\bgamma}(\cdot)=f_p(\cdot)+\bgamma^TA_1(\cdot),
\end{align}
where $A_1(\cdot)$ is defined in Theorem \ref{thm:asympn} and $\bgamma\in\RR^q$.
Combining \eqref{yp} and \eqref{gamma} yields
\begin{align}\label{lin-reg}
     y^{(p)}_j=f_p(\cdot)+\bgamma^TA_1(\cdot)+\varepsilon_j 
\end{align}
for $j=1,...,n$.
Note that \eqref{lin-reg} is a linear regression model with $\bgamma$ as parameter of interest, and the true value of $\bgamma$ is 0 regarding \eqref{yp}.

Since $\varepsilon_i$ follows a normal distribution, the maximum likelihood estimator for $\bgamma$ is equivalent to ordinary least square estimator, with regular asymptotic expression:
\begin{align}\label{gamman}
    \hat{\bgamma}_{n}=\frac{1}{n} \Wb^{-1} \sum_{i=1}^{n} \varepsilon_{i} A_1(\bx_i)+o_{p}\left(n^{-1 / 2}\right) 
\end{align}
where $\Wb$ is defined in \eqref{W}.

A natural estimator for $\btheta_S^*$ is 
\begin{align*}
    \hat\btheta_{Sn} = \argmin_{\btheta\in \bTheta}\| f_{p\hat\bgamma_n}(\cdot)-f_s(\cdot,\btheta)\|_{\cH_m(\Omega)},
\end{align*}
where we assume that $f_s(\cdot,\btheta)$ is cheap for simplicity.

Define 
\begin{align}
    \btheta(\bt) = \argmin_{\btheta\in \bTheta}\| f_{p\bt}(\cdot)-f_s(\cdot,\btheta)\|_{\cH_m(\Omega)}\label{thetat}
\end{align}
as a function of $\bt$. 
Let 
\begin{align*}
    \Upsilon(\btheta, \bt)&=-2\left\langle f_{p\bt}(\cdot)-f_s(\cdot,\btheta),\frac{\partial {f_s}(\cdot,\btheta)}{\partial\btheta}
   \right\rangle_{\cH_m(\Omega)}\nonumber\\
   &=-2\left\langle f_{p}(\cdot)-\bt^TA_1(\cdot)-f_s(\cdot,\btheta),\frac{\partial {f_s}(\cdot,\btheta)}{\partial\btheta}
   \right\rangle_{\cH_m(\Omega)}.
\end{align*}

Then \eqref{thetat} implies $\Upsilon(\btheta, \bt)=0$ for all $\bt$ near 0. From the implicit function theorem, we have
\begin{align}\label{partialt}
\frac{\partial \btheta(\bt)}{\partial \bt^{\mathrm{T}}} \bigg|_{\bt=0} &=-\left(\frac{\partial \Upsilon}{\partial \btheta^{\mathrm{T}}}\left(\btheta_S^{*}, 0\right)\right)^{-1} \frac{\partial \Upsilon}{\partial \bt^{\mathrm{T}}}\left(\btheta_S^{*}, 0\right).
\end{align}

First we calculate 
\begin{align}\label{partialtfirst}
    -\left(\frac{\partial \Upsilon}{\partial \btheta^{\mathrm{T}}}\left(\btheta_S^{*}, 0\right)\right)^{-1}&=-\left(\sum_{j=1}^\infty \frac{1}{{\gamma_j}}\frac{\partial^{2}}{\partial \btheta \partial \btheta^{\mathrm{T}}}\left(\int_{\Omega} (f_p(\bz)-f_s(\bz,\btheta^*_S))\phi_j(\bz) {\rm d} \bz\right)^2\right)^{-1} = -\Vb^{-1}.
\end{align}

Then we calculate
\begin{align}
&\frac{\partial \Upsilon}{\partial \bt^{\mathrm{T}}}\left(\btheta_S^{*}, 0\right)
\nonumber\\
=& \frac{\partial}{\partial \bt^{\mathrm{T}}}\left(-2\sum_{j=1}^\infty \frac{1}{\gamma_j}\left\langle f_p(\cdot)-\bt^\mathrm{T}A_1(\cdot)-f_s(\cdot,\btheta),\phi_j(\cdot)\right\rangle_{L_2(\Omega)}\left\langle \frac{\partial {f_s}(\cdot,\btheta)}{\partial\btheta},\phi_j(\cdot)\right\rangle_{L_2(\Omega)}\right)\bigg|_{\btheta=\btheta_S^{*}, \bt=0}\nonumber\\
=&2\sum_{j=1}^\infty \frac{1}{\gamma_j}\left\langle A_1(\cdot),\phi_j(\cdot)\right\rangle_{L_2(\Omega)}\left(\left\langle \frac{\partial {f_s}(\cdot,\btheta)}{\partial\btheta},\phi_j(\cdot)\right\rangle_{L_2(\Omega)}\right)^\mathrm{T}\bigg|_{\btheta=\btheta_S^{*}, \bt=0}\nonumber\\
=& 2\sum_{j=1}^\infty \frac{1}{\gamma_j}\nonumber\\
\times & \left\langle\sum_{i=1}^\infty\frac{1}{\gamma_i}\phi_i(\cdot)\int_\Omega\frac{\partial f_s(\bz,{\btheta})}{\partial\btheta}\phi_i(\bz){\rm d}\bz,\phi_j(\cdot)\right\rangle_{L_2(\Omega)}\left(\left\langle \frac{\partial {f_s}(\cdot,\btheta)}{\partial\btheta},\phi_j(\cdot)\right\rangle_{L_2(\Omega)}\right)^\mathrm{T}\bigg|_{\btheta=\btheta_S^{*}, \bt=0}\nonumber\\
=&2\sum_{j=1}^\infty \left(\frac{1}{\gamma_j}\int_\Omega\frac{\partial f_s(\bz,{\btheta})}{\partial\btheta}\phi_j(\bz){\rm d}\bz\right)\left(\frac{1}{\gamma_j}\int_\Omega\frac{\partial f_s(\bz,{\btheta})}{\partial\btheta}\phi_j(\bz){\rm d}\bz\right)^\mathrm{T}\bigg|_{\btheta=\btheta_S^{*}, \bt=0}\nonumber\\
=&2\EE(A_1(\bx)A_1(\bx)^\mathrm{T})\nonumber\\
=&2\Wb.\label{partialtsecond}
\end{align}

By combining \eqref{partialt}, \eqref{partialtfirst}, \eqref{partialtsecond}, we have
\begin{align}
    \frac{\partial \btheta(\bt)}{\partial \bt^{\mathrm{T}}} \bigg|_{\bt=0}=-2\Vb^{-1}\Wb.\label{partialthetat}
\end{align}

By the delta method, we obtain
\begin{align*}
    \hat{\btheta}_{Sn}-\btheta^{*}_S=\btheta\left(\hat{\gamma}_{n}\right)-\btheta(0)=\frac{\partial \btheta(\bt)}{\partial \bt^{\mathrm{T}}}\bigg|_{\bt=0} \hat{\gamma}_{n}+o_{p}\left(n^{-1 / 2}\right),
\end{align*}
which, together with \eqref{partialthetat} and \eqref{gamman}, yields
\begin{align}\label{thetahatsn}
    \hat{\btheta}_{Sn}-\btheta^{*}_S=-2\Vb^{-1}\left(\frac{1}{n}\sum_{i=1}^{n} \varepsilon_{i} A_1(\bx_i)\right)+o_{p}\left(n^{-1 / 2}\right).
\end{align}

Apparently, \eqref{thetahatsn} and \eqref{eq:thmmainconv} have the same form, which implies that the asymptotic expression of the Sobolev calibration is same as the maximum likelihood estimator for the $q$-dimensional parametric model. This finishes the proof.\hfill\BlackBox

\subsection{Proof of Corollary \ref{coro:L2}}\label{app:pfofcorols}
By taking $m=0$ in \eqref{thetastar}, the Sobolev calibration is equivalent to $L_2$ calibration. Furthermore, when $m=0$, since $\phi_j(\cdot)$'s are an orthonormal basis in  $L_2(\Omega)$, we can express $\Vb$ and $A_1(\bx)$ in Theorem \ref{thm:asympn} as
\begin{align}\label{eq:VL2}
    \Vb = & -\sum_{j=1}^\infty \frac{\partial^{2}}{\partial \btheta \partial \btheta^{\mathrm{T}}}\left(\int_{\Omega} (f_p(\bz)-f_s(\bz,\btheta^*))\phi_j(\bz) {\rm d} \bz\right)^2\nonumber\\
    = & \frac{\partial^{2}}{\partial \btheta \partial \btheta^{\mathrm{T}}}\left(-\sum_{j=1}^\infty\left(\int_{\Omega} (f_p(\bz)-f_s(\bz,\btheta^*))\phi_j(\bz) {\rm d} \bz\right)^2\right)\nonumber\\
    = & -\frac{\partial^{2}}{\partial \btheta \partial \btheta^{\mathrm{T}}}\int_{\Omega} (f_p(\bz)-f_s(\bz,\btheta^*))^2 {\rm d} \bz = \Vb_{L_2},
\end{align}
and 
\begin{align}\label{eq:A1L2}
    A_1(\bx) = \sum_{j=1}^\infty\phi_j(\bx)\int_\Omega\frac{\partial f_s(\bz,{\btheta_S^*})}{\partial\btheta}\phi_j(\bz){\rm d}\bz = \frac{\partial f_s(\bx,{\btheta_S^*})}{\partial\btheta}.
\end{align}
Thus, we can see that  Corollary \ref{coro:L2} is a direct result of Theorem \ref{thm:asympn}, \eqref{eq:VL2} and \eqref{eq:A1L2}.

\subsection{Proof of Proposition \ref{prop:gp+consist}}\label{app:gp+pfpropconsist}

Let $m_3 = (m_2-d/2+\max(m,d/2))/2$, thus we have $m_3>\max(d/2,m)$. By Condition (C3'), the sample path of $f_p(\cdot)$ lies in $W^{m_3}(\Omega)$ with probability one \citep{steinwart2019convergence_ec}. Hence,
\begin{align*}
    \|f_p(\cdot)\|_{W^{m_3}(\Omega)} <\infty 
\end{align*}
with probability one. Therefore, $\hat f_p(\cdot)$ has the same form as the kernel ridge regression estimator using an oversmoothed kernel function for $f_p(\cdot)\in W^{m_3}(\Omega)$, with $\lambda\asymp1/n$. Corollary 4.1 of \cite{fischer2020sobolev_ec} implies that $\left\|\hat f_p(\cdot)-f_p(\cdot)\right\|_{L_{2}(\Omega)}=o_\PP(1)$ and $\left\|\hat f_p(\cdot)-f_p(\cdot)\right\|_{W^{m_3}(\Omega)}=O_\PP(1)$. Thus,
\begin{align}\label{eq:GNineq}
    & \left\|\hat f_p(\cdot)-f_p(\cdot)\right\|_{\cH^{m}(\Omega)} \leq C_1\left\|\hat f_p(\cdot)-f_p(\cdot)\right\|_{W^m(\Omega)} \nonumber\\ 
    \leq & C_2 \left\|\hat f_p(\cdot)-f_p(\cdot)\right\|_{L_{2}(\Omega)}^{\frac{m_3-m}{m_3}}\left\|\hat f_p(\cdot)-f_p(\cdot)\right\|_{W^{m_3}(\Omega)}^{\frac{m}{m_3}}\nonumber\\
    =&  o_{\PP}(1).
\end{align}

From the definition of $\hat\btheta_S$ and $\btheta^*_S$ in \eqref{eq:thetahat} and \eqref{thetastar}, it suffices to prove that $\|\hat f_p(\cdot)-\hat f_s(\cdot, \btheta)\|_{\cH_m(\Omega)}$ converges to $\left\|f_p(\cdot)-f_s(\cdot, \btheta)\right\|_{\cH_m(\Omega)}$ uniformly with respect to $\btheta$ in probability. This is ensured by
\begin{align*}
    &\left|\left\|\hat f_p(\cdot)-\hat f_s(\cdot, \btheta)\right\|_{\cH_m(\Omega)}-\left\|f_p(\cdot)-f_s(\cdot, \btheta)\right\|_{\cH_m(\Omega)}\right| \nonumber\\
    \le &\left\|\hat f_p(\cdot)-\hat f_s(\cdot, \btheta)-(f_p(\cdot)-f_s(\cdot, \btheta))\right\|_{\cH_m(\Omega)} \nonumber\\
    \le &\left\|\hat f_p(\cdot)-f_p(\cdot)\right\|_{\cH_m(\Omega)}+\left\|\hat f_s(\cdot, \btheta)-f_s(\cdot, \btheta)\right\|_{\cH_m(\Omega)}\nonumber\\
    \leq & \left\|\hat f_p(\cdot)-f_p(\cdot)\right\|_{\cH_m(\Omega)}+C_3\left\|\hat f_s(\cdot, \btheta)-f_s(\cdot, \btheta)\right\|_{W^m(\Omega)}\nonumber\\
    =&o_{\PP}(1),
\end{align*}
where the first and second inequalities are from the triangle inequality, the third inequality is because of the equivalence of $\|\cdot\|_{W^m(\Omega)}$ and $\|\cdot\|_{\cH_m(\Omega)}$, and the last equality is guaranteed by \eqref{eq:GNineq} and Condition (C6).

\subsection{Proof of Theorem \ref{thm:gp+asympn}}\label{app:gp+asympn}

We first present several lemmas used in this proof. Lemma \ref{lem:gpnorm} is an intermediate step in the proof of Lemma F.8 of \cite{wang2021inference_ec}. 
Lemma \ref{Th:maximum} is a direct consequence of Theorems 1.3.3 and 2.1.1 of \cite{adler2009random_ec}.

\begin{lemma}\label{lem:gpnorm}
Suppose Conditions (C1') and (C3') hold. We have
$$
\sup_{\bx\in\Omega}\|\cK(\cdot,\bx)-r_1(\cdot)^{\mathrm{T}}\left(\Rb_1+\hat{\mu}_{n} I_{n}\right)^{-1} r_1(\bx)\|_{L_\infty(\Omega)} = O_{\PP}\left( n^{-\left(1-\frac{d}{2m_2}\right)}\right)
$$
where $r_1(\bx)=\left(\cK\left(\bx-\bx_{1}\right), \ldots, \cK\left(\bx-\bx_{n}\right)\right)^{\mathrm{T}}$, $\Rb_1=\left(\cK\left(\bx_{j}-\bx_{k}\right)\right)_{j k}$, and $\hat{\mu}_{n} \asymp C$ which is a constant.
\end{lemma}

\begin{lemma}\label{Th:maximum}
Let $Z_t$ be a centered separable Gaussian process on a $\mathfrak{d}$-compact $T$, where $\mathfrak{d}$ is the metric defined by
\begin{eqnarray*}
\mathfrak{d}(t_1,t_2)=\sqrt{\mathbb{E}(Z_{t_1}-Z_{t_2})^2}. 
\end{eqnarray*}
Then there exists a universal constant $K$ such that for all $u>0$,
\begin{align*}
\mathbb{P}\left(\sup_{t\in T} |Z_t| > K\int_0^{D/2} \sqrt{\log N(\epsilon,T,\mathfrak{d})}d\epsilon + u\right)\leq 2 e^{-u^2/2\sigma^2_T},
\end{align*}
where $\sigma^2_T = \sup_{t\in T} \mathbb{E}Z_t^2$, $N(\epsilon,T,\mathfrak{d})$ is the $\epsilon$-covering number of the metric space $(T,\mathfrak{d})$, and $D$ is the diameter of $T$.
\end{lemma}

\textit{{Proof of Theorem \ref{thm:gp+asympn}}}. 
For the briefness of the proof, we will use the same notation as in the proof of Theorem \ref{thm:asympn}, while keeping in mind that $f_p(\cdot)$ is a Gaussian process. The first half of the proof is merely repeating the proof of Theorem \ref{thm:asympn}, and the main difference comes from the second half of the proof, where we calculate $J_2$ defined in \eqref{eq:thm1pf1inners1}. For the first half of the proof, the only difference is that the term $I_2$ in \eqref{eq:fozeroest} is bounded by
\begin{align}\label{eq:gpasypf1I2}
    |I_2| = & o_{\PP}(n^{-1/2})\left(\left\|\hat f_p(\cdot)\right\|_{L_2(\Omega)}+\left\|f_s(\cdot,\hat\btheta_S)\right\|_{L_2(\Omega)}\right).
\end{align}
Since the sample path of $f_p(\cdot)$ lies in $\cH_{m_3}(\Omega)$ with probability one by Condition (C3') \citep{steinwart2019convergence_ec}, where $m_3 = (m_2-d/2+\max(m,d/2))/2$, together with Corollary 4.1 of \cite{fischer2020sobolev_ec}, we have that
\begin{align*}
    \left\|\hat f_p(\cdot)\right\|_{L_2(\Omega)} \leq \left\|\hat f_p(\cdot)-f_p(\cdot)\right\|_{L_2(\Omega)} + \left\|f_p(\cdot)\right\|_{L_2(\Omega)} =  O_{\PP}(1),
\end{align*}
which, together with \eqref{eq:gpasypf1I2} and Condition (C4'), implies that
\begin{align*}
    I_2= o_{\PP}(n^{-1/2}).
\end{align*}
Recall that 
\begin{align*}
    v_t(\btheta) & = 2\left\langle f_p(\cdot)-f_s(\cdot,\btheta),\frac{\partial {f_s}(\cdot,\btheta)}{\partial\theta_t}
   \right\rangle_{\cH_m(\Omega)},
\end{align*}
and 
\begin{align*}
J_2 = \sum_{j=1}^\infty\frac{2}{\gamma_j}\int_\Omega(\hat f_p(\bz)-{f_p}(\bz))\phi_j(\bz){\rm d}\bz\int_\Omega\frac{\partial {f_s}(\bz,\hat\btheta_S)}{\partial\theta_t}\phi_j(\bz){\rm d}\bz.
\end{align*}
Following the approach of obtaining \eqref{eq:thm1pf1J1J2sum2} as in the proof of Theorem \ref{thm:asympn}, we have
\begin{equation}\label{eq:gpthm1pf1J1J2sum2}
    \frac{\partial v_t(\tilde\btheta)}{\partial \btheta^{\mathrm{T}}}(\hat\btheta_S-\btheta_S^*) + J_2=o_{\PP}(n^{-1/2}),
\end{equation}
where $\tilde\btheta$ lies between $\btheta_S^*$ and $\hat\btheta_S$ and
\begin{align}\label{gpeqv}
\frac{\partial v(\tilde\btheta)}{\partial \btheta^{\mathrm{T}}} \stackrel{p}{\rightarrow} \Vb,
\end{align}
with $v(\btheta)=2\left\langle f_p(\cdot)-f_s(\cdot,\btheta),\frac{\partial {f_s}(\cdot,\btheta)}{\partial\btheta}\right\rangle_{\cH_m(\Omega)}=\left(v_1(\btheta),...,v_q(\btheta)\right)^{\mathrm{T}}$ and $\Vb$ as in Condition (C2'). 

Next, we consider $J_2$. Since $f_p(\cdot)$ is a Gaussian process, the calculation of $J_2$ is much different with that in the proof of Theorem \ref{thm:asympn}. We still define
\begin{align*}
    \mathcal{L}_{n}(f)=\frac{1}{n} \sum_{j=1}^{n}\left(y_j^{(p)}-f(\bx_j)\right)^{2}+\frac{\mu}{n}\|f\|_{\cN_{\cK}(\Omega)}^{2},
\end{align*}
and obtain
\begin{align}\label{eq:gpthm1pf1Ln1}
    0 = & \left.\frac{\partial}{\partial \eta} \mathcal{L}_n\left(\hat f_p(\cdot) + \eta A_{1t}(\cdot,\hat\btheta_S)\right)\right|_{\eta=0} \nonumber\\
    = & \frac{2}{n} \sum_{k=1}^{n}\left(\hat f_p\left(\bx_k\right)-f_p\left(\bx_k\right)\right) A_{1t}(\bx_k,\hat\btheta_S) -\frac{2}{n} \sum_{k=1}^{n} \varepsilon_k A_{1t}(\bx_k,\hat\btheta_S)\nonumber\\
    + &\frac{2\mu}{n}\left\langle\hat f_p(\cdot), A_{1t}(\cdot,\hat\btheta_S)\right\rangle_{\cN_{\cK}(\Omega)}\nonumber\\
    \coloneqq & Q_{1} - Q_{2} + Q_{3},
\end{align}
where we recall that $A_{1t}(\cdot,\hat\btheta_S)$ is defined by
\begin{align*}
    A_{1t}(\cdot,\btheta) = \sum_{j=1}^\infty\frac{1}{\gamma_j}\int_\Omega\frac{\partial {f_s}(\bz,{\btheta})}{\partial\theta_t}\phi_j(\bz){\rm d}\bz \cdot \phi_j(\cdot).
\end{align*}
Then $A_{1t}(\cdot,\hat\btheta_S)\in \cN_\cK(\Omega)$ is given by 
\begin{align}\label{eq:gpa1tnorm}
   \sup_{\btheta\in \bTheta} \|A_{1t}(\cdot,\btheta)\|_{\cN_\cK(\Omega)}^2 & = \sup_{\btheta\in \bTheta}\sum_{j=1}^\infty\frac{1}{\gamma_j^2\gamma_{\Psi,j}^{m_2/\alpha}}\left(\int_\Omega\frac{\partial {f_s}(\bz,{\btheta})}{\partial\theta_t}\phi_j(\bz){\rm d}\bz\right)^2 \nonumber\\
   &\leq C_1 \sup_{\btheta\in \bTheta} \left\|\frac{\partial {f_s}(\bz,{\btheta})}{\partial\theta_t}\right\|_{W^{2m+m_2}}^2 <\infty,
\end{align}
where the last inequality is guaranteed by Condition (C4').
Under the Gaussian process settings, $f_p(\cdot)$ is no longer in $\cN_{\cK}(\Omega)$. Therefore, we cannot follow the proof in Theorem \ref{thm:asympn}. In the following, we consider each term in \eqref{eq:gpthm1pf1Ln1}.

(1) Consider $Q_1$. Let
\begin{align*}
    g(\btheta) = \frac{1}{n} \sum_{k=1}^{n}\left(\hat f_p\left(\bx_k\right)-f_p\left(\bx_k\right)\right) A_{1t}(\bx_k,\btheta) - \int_\Omega \left(\hat f_p\left(\bz\right)-f_p\left(\bz\right)\right) A_{1t}(\bz,\btheta) {\rm d}\bz,
\end{align*}
which is a centered Gaussian process, since it is a linear transformation of a centered Gaussian process. Thus,
\begin{align*}
    & g(\btheta_1) - g(\btheta_2)\nonumber\\
    = & \frac{1}{n} \sum_{k=1}^{n}\left(\hat f_p\left(\bx_k\right)-f_p\left(\bx_k\right)\right) \left(A_{1t}(\bx_k,\btheta_1) - A_{1t}(\bx_k,\btheta_2)\right)\nonumber\\
    -& \int_\Omega \left(\hat f_p\left(\bz\right)-f_p\left(\bz\right)\right) \left(A_{1t}(\bz,\btheta_1) - A_{1t}(\bz,\btheta_2)\right) {\rm d}\bz.
\end{align*}
Denote
\begin{align}\label{eq:dabbre}
    &l(x,y)\nonumber\\
    =&\left(\cK(x,y) - \rb_1(x)^\mathrm{T}(\Rb_1 + \mu \Ib_n)^{-1}\rb_1(y)\right)\left(A_{1t}(x,\btheta_1) - A_{1t}(x,\btheta_2)\right)\left(A_{1t}(y,\btheta_1) - A_{1t}(y,\btheta_2)\right)
\end{align}
By \eqref{eq:dabbre}, the metric $\mathfrak{d}$ of $g(\cdot)$ can be computed by
\begin{align}\label{eq:gpmetricD}
     \mathfrak{d}(\btheta_1,\btheta_2)^2 =&  \mathbb{E}(g(\btheta_1)-g(\btheta_2))^2 \nonumber\\
    = & \frac{1}{n^2}\sum_{j=1}^n\sum_{k=1}^n l(\bx_k,\bx_j) -\frac{2}{n}\sum_{j=1}^n\int_\Omega l(\bz,\bx_j){\rm d}\bz
     + \int_\Omega\int_\Omega l(\bz_1,\bz_2) {\rm d}\bz_1{\rm d}\bz_2.
\end{align}

In order to apply Lemma \ref{Th:maximum}, we need to compute $\sigma^2_{\btheta} = \sup_{\btheta\in \bTheta} \mathbb{E}g(\btheta)^2$, the diameter $D$, and $\log N(\epsilon,T,\mathfrak{d})$.

Direct computation shows that 
\begin{align*}
    \EE g(\btheta)^2 = & \frac{1}{n^2}\sum_{j=1}^n\sum_{k=1}^n \left(\cK(\bx_k,\bx_j) - \rb_1(\bx_k)^\mathrm{T}(\Rb_1 + \mu \Ib_n)^{-1}\rb_1(\bx_j)\right)A_{1t}(\bx_k,\btheta) A_{1t}(\bx_j,\btheta)\nonumber\\
    & -\frac{2}{n}\sum_{j=1}^n\int_\Omega \left(\cK(\bz,\bx_j) - \rb_1(\bz)^\mathrm{T}(\Rb_1 + \mu \Ib_n)^{-1}\rb_1(\bx_j)\right)A_{1t}(\bz,\btheta) A_{1t}(\bx_j,\btheta){\rm d}\bz\nonumber\\
    & + \int_\Omega\int_\Omega \left(\cK(\bz_1,\bz_2) - \rb(\bz_1)^\mathrm{T}(\Rb_1 + \mu \Ib_n)^{-1}\rb_1(\bz_2)\right)A_{1t}(\bz_1,\btheta) A_{1t}(\bz_2,\btheta){\rm d}\bz_1{\rm d}\bz_2\nonumber\\
    = & \frac{1}{n}\sum_{j=1}^nA_{1t}(\bx_j,\btheta)\bigg(\frac{1}{n}\sum_{k=1}^n \left(\cK(\bx_k,\bx_j) - \rb_1(\bx_k)^\mathrm{T}(\Rb_1 + \mu \Ib_n)^{-1}\rb_1(\bx_j)\right)A_{1t}(\bx_k,\btheta) \nonumber\\
    & - \int_\Omega\left(\cK(\bz_1,\bx_j) - \rb_1(\bz_1)^\mathrm{T}(\Rb_1 + \mu \Ib_n)^{-1}\rb_1(\bx_j)\right)A_{1t}(\bz_1,\btheta_1) {\rm d}\bz_1\bigg)\nonumber\\
    & -\bigg( \int_\Omega A_{1t}(\bz_1,\btheta_1)\bigg(\frac{1}{n}\sum_{k=1}^n \left(\cK(\bz_1,\bx_k) - \rb_1(\bz_1)^\mathrm{T}(\Rb_1 + \mu \Ib_n)^{-1}\rb_1(\bx_k)\right)A_{1t}(\bx_k,\btheta_1) \nonumber\\
    & - \int_\Omega \left(\cK(\bz_1,\bz_2) - \rb_1(\bz_1)^\mathrm{T}(\Rb_1 + \mu \Ib_n)^{-1}\rb_1(\bz_2)\right) A_{1t}(\bz_2,\btheta_1){\rm d}\bz_2\bigg){\rm d}\bz_1\bigg).
\end{align*}
Let
\begin{align*}
    \cF_1 = & \bigg\{h(\cdot): h(\cdot) =  \frac{1}{n}\sum_{k=1}^n \left(\cK(\cdot,\bx_k) - \rb_1(\cdot)^\mathrm{T}(\Rb_1 + \mu \Ib_n)^{-1}\rb_1(\bx_k)\right)A_{1t}(\bx_k,\btheta) \nonumber\\
    & - \int_\Omega \left(\cK(\cdot,\bz_2) - \rb_1(\cdot)^\mathrm{T}(\Rb_1 + \mu \Ib_n)^{-1}\rb_1(\bz_2)\right)A_{1t}(\bz_2,\btheta){\rm d}\bz_2, \btheta\in\bTheta\bigg\},\nonumber\\
    \cG_1 = & \{g(\cdot):g(\cdot)=A_{1t}(\cdot,\btheta),\btheta\in\bTheta\}.
\end{align*}
Thus, $\EE g(\btheta)^2$ can be bounded by
\begin{align*}
    \EE g(\btheta)^2 \leq &  \sup_{h\in\cF_1, g\in \cG_1} \left|(P_n - P)hg\right|.
\end{align*}
Note that $\EE g(\btheta)^2$ does not involve Gaussian process, which allows us to repeat the procedure in deriving \eqref{eq:pf1th1Q1}. In order to do so, note that
\begin{align*}
    \sup_{h\in \cF_1} \|h(\cdot)\|_{L_\infty(\Omega)} \leq \sup_{h_1\in \cF_2, g\in \cG_1} \left|(P_n - P)h_1g\right|,
\end{align*}
where $\cF_2 = \{h_1(\cdot) : h_1(\cdot) = \cK(\cdot,\bx) - \rb_1(\cdot)^\mathrm{T}(\Rb_1 + \mu \Ib_n)^{-1}\rb_1(\bx), \bx\in\Omega\}$. Clearly, $\cF_2\subset \cN_{\cK}(\Omega)$, and Lemma \ref{lem:gpnorm} implies that
\begin{align*}
    \sup_{h_1\in \cF_2}\|h_1(\cdot)\|_{L_\infty} = O_{\PP}(n^{-(1-\frac{d}{2m_2})}).
\end{align*}
Thus, repeating the procedure in deriving \eqref{eq:pf1th1Q1}, we can obtain that
\begin{align}\label{eq:gpf2g1}
    \sup_{h_1\in \cF_2, g\in \cG_1} \left|(P_n - P)h_1g\right| = O_{\PP}\left(n^{-(1-\frac{d}{2m_2})-\frac{1}{2}}\right).
\end{align}
Applying the similar procedure to $\sup_{h\in\cF_1, g\in \cG_1} \left|(P_n - P)hg\right|$, we obtain that
\begin{align}\label{eq:gpf1g1}
    \sup_{h\in\cF_1, g\in \cG_1} \left|(P_n - P)hg\right| = &  O_{\PP}\left(n^{-(1-\frac{d}{2m_2})(1-\frac{d}{2m_2})-\frac{1}{2}(1-\frac{d}{2m_2})}n^{-\frac{1}{2}}\right)\nonumber\\
    = & O_{\PP}\left(n^{-1 - \frac{(4m_2 - d)(m_2-d)}{4m_2^2}}\right) = O_{\PP}\left(n^{-1 - \delta_0}\right).
\end{align}
Since $m_2 > d$, we have $\delta_0>0$. For the briefness of the proof, we omit the details of deriving \eqref{eq:gpf2g1} and \eqref{eq:gpf1g1} here.
Thus, we obtain
\begin{align*}
\sigma^2_{\btheta} = \sup_{\btheta\in \bTheta} \mathbb{E}g(\btheta)^2 = O_{\PP}\left(n^{-1 - \delta_0}\right).
\end{align*}
The diameter $D$ can be computed by
\begin{align*}
    D^2 = & \sup_{\btheta_1,\btheta_2\in\Theta}\mathfrak{d}(\btheta_1,\btheta_2)^2\nonumber\\
    \leq & 2\sup_{\btheta_1,\btheta_2\in\Theta} \EE g(\btheta_1)^2 + \EE g(\btheta_2)^2\nonumber\\
   = &O_{\PP}\left(n^{-1 - \delta_0}\right).
\end{align*}
It remains to bound $\log N(\epsilon,T,\mathfrak{d})$.By \eqref{eq:dabbre} and \eqref{eq:gpmetricD} , we have
\begin{align}\label{eq:gpmetricD2}
    \mathfrak{d}(\btheta_1,\btheta_2)^2
    \leq &\frac{2}{n^2}\sum_{j=1}^n\sum_{k=1}^n l(\bx_k,\bx_j) +2 \int_\Omega\int_\Omega l(\bz_1,\bz_2){\rm d}\bz_1{\rm d}\bz_2 \nonumber\\
    \leq &C_2\sup_{\bx\in \Omega}\left|A_{1t}(\bx,\btheta_1) - A_{1t}(\bx,\btheta_2)\right|^2,
\end{align}
where the last inequality is because $\left(\cK(\bx_k,\bx_j) - \rb_1(\bx_k)^\mathrm{T}(\Rb_1 + \mu \Ib_n)^{-1}\rb_1(\bx_j)\right) \leq 1 $ for all $\bx_k,\bx_j \in \Omega$, and 
$$\left(A_{1t}(\bx_k,\btheta_1) - A_{1t}(\bx_k,\btheta_2)\right)\left(A_{1t}(\bx_j,\btheta_1) - A_{1t}(\bx_j,\btheta_2)\right) \leq \sup_{\bx\in \Omega}\left|A_{1t}(\bx,\btheta_1) - A_{1t}(\bx,\btheta_2)\right|^2.$$
The mean value theorem implies 
\begin{align}\label{eq:gpmetricD3}
    \left|A_{1t}(\bx,\btheta_1) - A_{1t}(\bx,\btheta_2)\right| = & \left(\sum_{j=1}^\infty\frac{1}{\gamma_j}\int_\Omega\frac{\partial {f_s}(\bz,\tilde{\btheta})}{\partial\theta_t\partial\btheta}\phi_j(\bz){\rm d}\bz \cdot \phi_j(\bx)\right)^\mathrm{T}(\btheta_1-\btheta_2)\nonumber\\
    \leq & \left\|\sum_{j=1}^\infty\frac{1}{\gamma_j}\int_\Omega\frac{\partial {f_s}(\bz,\tilde{\btheta})}{\partial\theta_t\partial\btheta}\phi_j(\bz){\rm d}\bz \cdot \phi_j(\bx)\right\|_{L_2(\Omega)}\|\btheta_1-\btheta_2\|_2,
\end{align}
where $\tilde{\btheta}$ lies between $\btheta_1$ and $\btheta_2$.
Condition (C4') and Theorem 10.29 of \cite{wendland2004scattered_ec} gives us
\begin{align*}
   & \sup_{\btheta\in \bTheta} \left\|\sum_{j=1}^\infty\frac{1}{\gamma_j}\int_\Omega\frac{\partial {f_s}(\bz,\tilde{\btheta})}{\partial\theta_{t_1}\partial\theta_{t_2}}\phi_j(\bz){\rm d}\bz \cdot \phi_j(\cdot)\right\|_{\cN_\cK(\Omega)}^2\nonumber\\
    =& \sup_{\btheta\in \bTheta}\sum_{j=1}^\infty\frac{1}{\gamma_j^2\gamma_{\Psi,j}^{\beta_1}}\left(\int_\Omega\frac{\partial {f_s}(\bz,\tilde{\btheta})}{\partial\theta_{t_1}\partial\theta_{t_2}}\phi_j(\bz){\rm d}\bz\right)^2 \leq C_3 \left\|\frac{\partial {f_s}(\bz,\tilde{\btheta})}{\partial\theta_{t_1}\partial\theta_{t_2}}\right\|_{W^{2m+m_2}(\Omega)}\leq C_4,
\end{align*}
where $\beta_1=\frac{m_2}{\alpha}$,
and the last inequality is by Condition (C4') and the equivalence of $\cN_\alpha^{m/\alpha}(\Omega)$ and $\cH_m(\Omega)$. The Sobolev embedding theorem implies that 
\begin{align*}
    &\left\|\sum_{j=1}^\infty\frac{1}{\gamma_j}\int_\Omega\frac{\partial {f_s}(\bz,\tilde{\btheta})}{\partial\theta_{t_1}\partial\theta_{t_2}}\phi_j(\bz){\rm d}\bz \cdot \phi_j(\cdot)\right\|_{L_2(\Omega)}\\
    \leq & C_5 \left\|\sum_{j=1}^\infty\frac{1}{\gamma_j}\int_\Omega\frac{\partial {f_s}(\bz,\tilde{\btheta})}{\partial\theta_{t_1}\partial\theta_{t_2}}\phi_j(\bz){\rm d}\bz \cdot \phi_j(\cdot)\right\|_{W^{\alpha\beta_1}(\Omega)}\leq C_6, 
\end{align*}
which, together with \eqref{eq:gpmetricD2} and \eqref{eq:gpmetricD3}, implies 
\begin{align*}
    \mathfrak{d}(\btheta_1,\btheta_2)^2 \leq C_7 \|\btheta_1-\btheta_2\|_2^2.
\end{align*}
Therefore, the covering number is bounded above by
\begin{align}\label{e1}
\log N(\epsilon,\bTheta,\mathfrak{d}) \leq \log N\bigg(\frac{\epsilon}{C_7},\bTheta,\|\cdot\|_2\bigg).
\end{align}
The right side of \eqref{e1} is just the covering number of a Euclidean ball, which is well understood in the literature; See Lemma 4.1 of \cite{pollard1990empirical_ec}. Thus, we have
\begin{align*}
\log N(\epsilon,\bTheta,\mathfrak{d}) \leq C_8\log \bigg(1 + \frac{C_7}{\epsilon}\bigg).
\end{align*}
Therefore, the integral in Lemma \ref{Th:maximum} can be computed by
\begin{align*}
\int_0^{D/2} \sqrt{\log N(\epsilon,\bTheta,\mathfrak{d})}d\epsilon & \leq \int_0^{D/2} \sqrt{C_8\log \bigg(1 + \frac{C_7}{\epsilon}\bigg)}{\rm d}\epsilon\nonumber\\
& \leq \left(\int_0^{D/2}{\rm d}\epsilon \right)^{1/2}\left(\int_0^{D/2}C_8 \log \bigg(1 + \frac{C_7}{\epsilon}\bigg)d\epsilon\bigg)\right)^{1/2}\\
& \leq C_{9} D\sqrt{\log\left(\frac{1}{D}\right)}.
\end{align*}
Now we can apply Lemma \ref{Th:maximum} and obtain
\begin{align*}
\sup_{\btheta\in \bTheta} |g(\btheta)| = O_{\PP}(n^{\frac{-1-\delta_0}{2}}\log(n)) = o_{\PP}(n^{-1/2}).
\end{align*}
Therefore, we conclude that 
\begin{align}\label{eq:gppf1th1Q1}
    Q_1 = 2\int_\Omega \left(\hat f_p\left(\bz\right)-f_p\left(\bz\right)\right) A_{1t}(\bz,\hat\btheta_S) {\rm d}\bz + o_{\PP}(n^{-1/2}).
\end{align}

(2) Consider $Q_3$. By the Cauchy-Schwarz inequality,  
\begin{align}\label{eq:gppf1th1Q5} 
    |Q_3|\leq 2 \frac{\mu}{n}\left\|\hat f_p(\cdot)\right\|_{\cN_{\cK}(\Omega)} \left\| A_{1t}(\cdot,\hat\btheta_S)\right\|_{\cN_{\cK}(\Omega)}.
\end{align}
By \eqref{eq:gpa1tnorm}, we have that $\left\| A_{1t}(\cdot,\hat\btheta_S)\right\|_{\cN_{\cK}(\Omega)}\leq C_{10}$. Recall that 
\begin{align*}
   \hat f_p(\bx) = \rb_1(\bx)^\mathrm{T}\left(\Rb_1+\mu \Ib_{n}\right)^{-1} \by^{(p)} = \rb_1(\bx)^\mathrm{T}\left(\Rb_1+\mu \Ib_{n}\right)^{-1} \bF +\rb_1(\bx)^\mathrm{T}\left(\Rb_1+\mu \Ib_{n}\right)^{-1} \bvarepsilon,
\end{align*}
where $\bF=(f_p(\bx_1),...,f_p(\bx_n))^\mathrm{T}$ and $\bvarepsilon =(\varepsilon_1,...,\varepsilon_n)^\mathrm{T}$.
The term $\left\|\hat f_p(\cdot)\right\|_{\cN_{\cK}(\Omega)}$ can be computed by
\begin{align*}
    \left\|\hat f_p(\cdot)\right\|_{\cN_{\cK}(\Omega)}^2 = (\bF + \bvarepsilon)^\mathrm{T}\left(\Rb_1+\mu \Ib_{n}\right)^{-1}\Rb_1\left(\Rb_1+\mu \Ib_{n}\right)^{-1} (\bF + \bvarepsilon).
\end{align*}
Since both $\bF$ and $\bvarepsilon$ are normally distributed, the expectation can be computed by
\begin{align*}
    \EE\left\|\hat f_p(\cdot)\right\|_{\cN_{\cK}(\Omega)}^2 = & \EE{\rm tr}\left(\left(\Rb_1+\mu \Ib_{n}\right)^{-1}\Rb_1\left(\Rb_1+\mu \Ib_{n}\right)^{-1} (\bF + \bvarepsilon)(\bF + \bvarepsilon)^\mathrm{T}\right)\nonumber\\
    = & {\rm tr}\left(\left(\Rb_1+\mu \Ib_{n}\right)^{-1}\Rb_1\right) \leq {\rm tr}\Ib_{n} = n,
\end{align*}
and the variance can be bounded by
\begin{align*}
    {\rm Var}(\|\hat f_p(\cdot)\|_{\cN_{\cK}(\Omega)}^2) \leq & \EE\left((\bF + \bvarepsilon)^T\left(\Rb_1+\mu \Ib_{n}\right)^{-1}\Rb_1\left(\Rb_1+\mu \Ib_{n}\right)^{-1} (\bF + \bvarepsilon)\right)^2 \nonumber\\
    \leq & \EE\left((\bF + \bvarepsilon)^\mathrm{T}\left(\Rb_1+\mu \Ib_{n}\right)^{-1} (\bF + \bvarepsilon)\right)^2 \nonumber\\
    = & \EE\left(\sum_{j=1}^n Z_j^2\right)^2 \leq C_{11} n^2,
\end{align*}
where the last equality is because $\bF + \bvarepsilon$ are normally distributed with mean zero and covariance $\Rb_1+\mu \Ib_{n}$, and $Z_js$ are i.i.d standard normally distributed random variables. Thus, Chebyshev's inequality implies that 
\begin{align}\label{eq:gpth1Q5t2} 
     \left\|\hat f_p(\cdot)\right\|_{\cN_{\cK}(\Omega)}^2 = \EE\left\|\hat f_p(\cdot)\right\|_{\cN_{\cK}(\Omega)}^2 + O_{\PP}(\sqrt{{\rm Var}(\|\hat f_p(\cdot)\|_{\cN_{\cK}(\Omega)}^2)}) = O_{\PP}(n).
\end{align}
Plugging \eqref{eq:gpth1Q5t2} into \eqref{eq:gppf1th1Q5} leads to
\begin{align}\label{eq:gppf1th1Q5t} 
    |Q_3|= O_{\PP}(n^{-1}n^{1/2}) = O_{\PP}(n^{-1/2}).
\end{align}
By combining \eqref{eq:gpthm1pf1J1J2sum2},\eqref{gpeqv},\eqref{eq:gpthm1pf1Ln1}, \eqref{eq:gppf1th1Q1} and \eqref{eq:gppf1th1Q5t}, we obtain that
\begin{equation*}
    (\Vb+\Eb)(\hat\btheta_S-\btheta_S^*) + \frac{2}{n} \sum_{k=1}^{n} \varepsilon_k A_{1}(\bx_k) = O_{\PP}(n^{-1/2}).
\end{equation*}
where each element in $\Eb$ is $o_{\PP}(1)$ and 
\begin{align*}
    A_{1}(\cdot)= \sum_{j=1}^\infty\frac{1}{\gamma_j}\int_\Omega\frac{\partial {f_s}(\bz,{\btheta^*_S})}{\partial\btheta}\phi_j(\bz){\rm d}\bz \cdot \phi_j(\cdot)=\left(A_{11}(\cdot,\btheta^*_S),...,A_{1q}(\cdot,\btheta^*_S)\right)^{\mathrm{T}}.
\end{align*}

Since $\Eb$ converges to zero in probability, we have
\begin{equation*}
    \hat\btheta_S-\btheta_S^* = -2\Vb^{-1} \left(\frac{1}{n} \sum_{k=1}^{n} \varepsilon_k A_{1}(\bx_k,\btheta_S^*)\right) + O_{\PP}(n^{-1/2}),
\end{equation*}
which finishes the proof. \hfill\BlackBox

\section{Experimental Details and Additional Figures}

\subsection{Uncertainty Quantification Details}\label{subsec_UQdetails}

{In this subsection, we provide experimental details of uncertainty quantification in simulation studies. The key steps of uncertainty quantification are summarized in Algorithm \ref{alg:uq}.

\begin{algorithm}
{
\renewcommand{\algorithmicrequire}{\textbf{Input:}}
\renewcommand{\algorithmicensure}{\textbf{Output:}}
\caption{Uncertainty Quantification for Sobolev Calibration}\label{alg:uq}
\begin{algorithmic}
\Require Observations $\{(\bx_j,y_j^{(p)}\}_{j=1}^n$, a design $\tilde \bX$, estimations $\hat{\btheta}_S$ and $\hat{\sigma}^2$, approximations $\tilde{\gamma_{j}}$, $\tilde{\phi_j}$, $\hat{f}_s$, parameters $m$, $N$, $l$. 
%\State Estimate $\hat{\delta}$ based on samples $\{(\bx_j,y_j^{(p)}-\hat{f}_{s}(\bx_j,\hat{\btheta}_{S})\}_{j=1}^n$.
\State Set $\Vb \gets -\frac{\partial^{2}}{\partial \btheta \partial \btheta^{\mathrm{T}}}\|\mathcal{I}_{\cK_{m},\tilde \bX}\left(\hat{f}_p(\cdot)-\hat{f}_s(\cdot,\hat{\btheta}_S)\right)\|^2_{\cH_m(\Omega)}$; \\
\hspace{6mm}$\Wb \gets \sum_{j=1}^N\frac{1}{\tilde{\gamma_j}^2}\left(\int_\Omega\frac{\partial \hat{f}_s(\bz,{
\hat{\btheta}_S})}{\partial\btheta}\tilde{\phi_j}(\bz){\rm d}\bz\right)^2$;
\State \hspace{6mm}$\hat{\sigma}_S^2\gets 4 \hat{\sigma}^{2} \Vb^{-1} \Wb \Vb^{-1}/n$ (see Section \ref{subapp_uq_sampleV} for details).
\For{$i\leftarrow 1$ to $l$}
\begin{enumerate}
\item $\hat{f}_{s,i}(\cdot) \gets \hat{f}_{s}(\cdot,\hat{\btheta}_{S,i})$, where $\hat{\btheta}_{S,i}$ is generated from Gaussian distribution $\mathcal{N}(\hat{\btheta}_S,\hat{\sigma}_S^2)$.
\item Use Gaussian process to estimate $\delta(\cdot,\hat{\btheta}_{S,i})$ based on samples $\{(\bx_j,y_j^{(p)}-\hat{f}_{s,i}(\bx_j)\}_{j=1}^n$, and obtain the 95\% upper bound $\bar{\delta}_i(\cdot)$ and lower bound $\underline{\delta}_i(\cdot)$.
\end{enumerate}
\EndFor
\State Calculate the empirical 0.975 and 0.025 quantile $\hat{f}_{s}^{0.975}(\cdot)$ and $\hat{f}_{s}^{0.025}(\cdot)$ from $\{\hat{f}_{s,i}(\cdot)\}_{i=1}^{l}$.
\State Calculate the empirical 0.975 and 0.025 quantile $\hat{f}_{p}^{0.975}(\cdot)$ and $\hat{f}_{p}^{0.025}(\cdot)$ from $\{\hat{f}_{s,i}(\cdot)+\bar{\delta}_i(\cdot)\}_{i=1}^{l}$ and $\{\hat{f}_{s,i}(\cdot)-\underline{\delta}_i(\cdot)\}_{i=1}^{l}$ resepectively.
\Ensure \State 95\% Confidence interval of ${\btheta}_S^*$: $[\hat{\btheta}_S-1.96\times\hat{\sigma}_S, \hat{\btheta}_S+1.96\times\hat{\sigma}_S]$;
\State 95\% Point-wise confidence band of ${f}_s(\cdot,{\btheta}_S^*)$: $[\hat{f}_{s}^{0.025}(\cdot),\hat{f}_{s}^{0.975}(\cdot)]$;
\State 95\% Point-wise confidence band of $f_p(\cdot)$: $[\hat{f}_{p}^{0.025}(\cdot),\hat{f}_{p}^{0.975}(\cdot)]$.
\end{algorithmic}
}
\end{algorithm}
} 

{\subsubsection{Computation of Sample Variance of $\hat{\btheta}$ }\label{subapp_uq_sampleV}

In order to provide confidence intervals of the estimate $\hat{\btheta}$ for each trial and compute the average coverage rate, we need to numerically compute the asymptotic variance of $\hat{\btheta}$ as given in Theorem \ref{thm:asympn} and Theorem \ref{thm:gp+asympn} based on each sample set, which requires approximations to $\sigma^2$,$\Vb$ and $\Wb$. 

We here take Sobolev calibration as an example. The asymptotic variance of $\hat{\btheta}_{L_2}$ and $\hat{\btheta}_{KO}$ can be calculated in a similar manner. First, $\Vb$ can be approximated by 
$
\Vb:=-\sum_{j=1}^\infty \frac{1}{{\gamma_j}}\frac{\partial^{2}}{\partial \btheta \partial \btheta^{\mathrm{T}}}\left(\int_{\Omega} (f_p(\bz)-f_s(\bz,\btheta^*_S))\phi_j(\bz) {\rm d} \bz\right)^2
\approx -\frac{\partial^{2}}{\partial \btheta \partial \btheta^{\mathrm{T}}}\|\mathcal{I}_{\cK_{m},\tilde \bX}\left(f_p(\cdot)-f_s(\cdot,\btheta^*_S)\right)\|^2_{\cH_m(\Omega)}
\approx -\frac{\partial^{2}}{\partial \btheta \partial \btheta^{\mathrm{T}}}\|\mathcal{I}_{\cK_{m},\tilde \bX}\left(\hat{f}_p(\cdot)-\hat{f}_s(\cdot,\hat{\btheta}_S)\right)\|^2_{\cH_m(\Omega)}$, where the penultimate step is introduced in Section \ref{sec:comp}. For the computation of $\Wb$, we need to numerically calculate the eigenvalues and eigenfunctions of Mat\'ern kernel, which can be achieved by Nystr{\"o}m method \citep{williams2000using}. The details of approximation error and the number of truncations $N$ can be found in  Section \ref{app:eigenerrors}. Then we calculate 
\begin{align*}
\Wb & \coloneqq\EE(A_1(\bx)A_1(\bx)^{T})=\sum_{j=1}^\infty\frac{1}{\gamma_j^2}\left(\int_\Omega\frac{\partial f_s(\bz,{\btheta_S^*})}{\partial\btheta}\phi_j(\bz){\rm d}\bz\right)^2\nonumber\\
& \approx \sum_{j=1}^N\frac{1}{\tilde{\gamma_j}^2}\left(\int_\Omega\frac{\partial f_s(\bz,{\btheta_S^*})}{\partial\btheta}\tilde{\phi_j}(\bz){\rm d}\bz\right)^2\approx\sum_{j=1}^N\frac{1}{\tilde{\gamma_j}^2}\left(\int_\Omega\frac{\partial \hat{f}_s(\bz,{
\hat{\btheta}_S})}{\partial\btheta}\tilde{\phi_j}(\bz){\rm d}\bz\right)^2,
\end{align*} 
where $\tilde \phi_j$ and $\tilde \gamma_j$ are approximated eigenfunctions and eigenvalues, respectively. For the simplicity of computation, we assume that $\sigma^2$ and $f_s(\cdot)$ are known. If the true variance of random noise $\sigma^2$ is unknown, in practice, it can be estimated by collecting physical experiments at some fixed points repeatedly. Note that since Example \ref{gp-case3} considers $f_p(\cdot)$ as a Gaussian process, the calculation procedures of each trial are based on a random sample path of $f_p(\cdot)$. }

{\subsubsection{Computation of Point-wise
Confidence Band of ${f}_s(\cdot,{\btheta}_S^*)$ and $f_p(\cdot)$}

Recall that we have $\hat{\btheta}$ follows Gaussian distribution asymptotically. After calculating the estimated asymptotic variance, a random sample set $\{\hat{\btheta}_{1},...,\hat{\btheta}_{100}\}$ can be generated by sampling 100 times from the Gaussian distribution $\mathcal{N}(\hat{\btheta},\hat{\sigma}^2)$. Then we can plug $\{\hat{\btheta}_{1},...,\hat{\btheta}_{100}\}$ into $\hat{f}_s(\cdot,{\btheta})$ and obtain a function set $\{\hat{f}_{s,1}(\cdot),...,\hat{f}_{s,100}(\cdot)\}$, from which we can find the 0.025 and 0.975 sample quantiles and obtain the confidence band. 

The computation of point-wise confidence band of $f_p(\cdot)$ is in a similar manner, but we need to estimate the function $\delta(\cdot,\hat{\btheta})$ based on samples $\{(\bx_j,y_j^{(p)}-\hat{f}_s(\bx_j,\hat{\btheta})\}_{j=1}^n$. The estimation error can lead to the bias of confidence interval, which can have some impact especially when the original variance is small. Therefore, to guarantee the validity of confidence interval, we need to ensure the estimation of $\hat{\delta}(\cdot,\hat{\btheta})$ can cover ${\delta}(\cdot,\hat{\btheta})$. Here we adopt the Gaussian process to estimate $\delta(\cdot,\hat{\btheta})$ together with 95\% confidence region, and use upper bound and lower bound for the computation of 0.975 sample quantile and 0.025 sample quantile separately.}

{\subsubsection{Approximation Errors of Eigenvalues and Eigenfunctions}\label{app:eigenerrors}

In this section, we provide some implementation details of Nystr{\"o}m method \citep{williams2000using} and the corresponding approximation errors.

Mercer's theorem implies that
\begin{align}
    \Psi_m(\bs,\bt)=& \sum_{j=1}^\infty \gamma_{j} \phi_j(\bs)\phi_j(\bt), \forall \bs,\bt\in\Omega,\nonumber
\end{align}
where $\gamma_{j}$ and $\phi_j(\cdot)$ are the eigenvalues and eigenfunctions of $\Psi_m(\cdot,\cdot)$ respectively. We approximate the first $N$ eigenvalues and eigenfunctions, such that
\begin{align}
        \Psi_m(\bs,\bt)\approx& \sum_{j=1}^N \tilde{\gamma_{j}} \tilde{\phi_j}(\bs)\tilde{\phi_j}(\bt), \forall \bs,\bt\in\Omega.\nonumber
\end{align}
Note that empirically $N$ should decrease with the increase of smoothness $m$, because the eigendecay rate of kernel $\Psi_m(\cdot)$ would get faster.
To empirically measure the approximation error, we first choose $1,000$ equally spaced points from 0 to 1, denoted as $\{s_j,j=1,...,1,000\}$. Then we calculate $\sum_{j=1}^{1,000}e(s_j)/1,000$, where $e(s)\coloneqq|\Psi_m(s,s)-\sum_{j=1}^N \tilde{\gamma_{j}} \tilde{\phi_j}(s)\tilde{\phi_j}(s)|$. Table \ref{tab:eigenerror} summarizes truncation points $N$ we choose and approximation errors corresponding to the smoothness $m$, and we can find that the approximation error is controlled below $0.05$. 
\begin{table}[!h]
 \caption{\label{tab:eigenerror}Truncation points and approximation errors}
\centering
\begin{tabular}{l ccccc}
\hline
{$m$}        & {7/8}    & {1}    & {9/8}    & {9/5}  & {2}    \\ \hline
{$N$}        & {20}     & {12}     & {12}     &{5}      & {5}      \\
{$Error$} & {0.0466} & {0.0303} & {0.0135} & {0.0038} & {0.0020} \\ \hline
\end{tabular}
\end{table}}

\subsection{Gaussian Process Details}\label{exp-gpdetails}

To select the best $\mu$ in \eqref{eq:gpfhatp}, we randomly generate a data set with 300 samples, and divide 30\% of samples as the validation set.
During each round, a realization of $f_p(\cdot)$ is generated and physical observations are attained based on this realization. The corresponding Sobolev calibration parameter is defined as in \eqref{thetastar},
where $\cH_m(\Omega)$ is the RKHS generated by $\Psi_{1}(x)$. 

Note that $\theta^*$ is a random variable, since each realization of $f_p(\cdot)$ leads to different $\theta^*$. Therefore, we let $d_k^2=(\hat\theta_{Sk}-\theta^*_k)^2$, where the subscript $k$ defines each round of experiment, and then calculate the mean and standard deviation of $d_k^2,k=1,...,500,$ to measure the performance of different methods. 

\subsection{Additional Experiments}\label{app_exp}

\begin{example} \label{dtm-case1}
Suppose the physical experiment is
$
    f_p(x)=\exp (\pi x / 5) \sin 2 \pi x, 
$
for $x\in \Omega=[0,1]$.
The physical observations are given by 
$
    y^{(p)}_j=f_p(x_j)+\varepsilon_j,
$
with $x_j\sim$ Uniform(0,1), $\varepsilon_j\sim N(0,\sigma_\epsilon^2)$ for $j=1,...,200$. We investigate two different $\sigma_\epsilon^2$ ($\sigma_\epsilon^2=0.05$ and $\sigma_\epsilon^2=0.1$) to evaluate the robustness of different methods under different information-noise ratio.

Suppose the computer experiment is
\begin{align*}
 f_s(x, \theta)=f_p(x)-\sqrt{\theta^{2}-\theta+1}\left(\sin 2 \pi \theta x+\exp\left(\pi\theta x/2\right)\right),   
\end{align*}
which is imperfect since $\sqrt{\theta^{2}-\theta+1}$ is always positive. Therefore, the computer experiment cannot perfectly match the physical experiment. For the ease of computation, we assume the computer code is \textit{cheap}.
We follow Approach 2 to select function space $\cH_m(\Omega)$ and $\cH_{m_1}(\Omega)$ as
the RKHS generated by Mat\'ern kernel $\Psi_{1}(x)$ and $\Psi_{2}(x)$, respectively. 
Since $0$ (used in $L_2$ calibration) $<m<m_1$ (used in the frequentist KO calibration), we intuitively interpret the frequentist KO calibration and $L_2$ calibration as using relatively smoother and less smooth kernels compared with Sobolev calibration, respectively. The tuning parameter $\lambda$ in \eqref{eq:krr} of the regularized estimator $\hat f_p(x)$ is selected by the generalized cross validation \citep{wahba1990},
and $\theta^*_{S}$ in the Sobolev calibration is defined as in \eqref{thetastar}.

Figure \ref{fig:Eg1} plots the RKHS norm of the discrepancy function with respect to $\theta \in [-1,1]$. Numerical optimization shows that $\theta^*_{S}\approx -0.08$.
In Figure \ref{fig:Eg1}(b), we interpret $L_2$ calibration as an overfitting result because it only considers the best approximation in function value, and this sometimes can cause the computer model ``wiggling''; see Example \ref{dtm-case2} for more details. In contrast, by specifying an over-smoothed kernel, KO is more conservative and may result in an inaccurate calibrated experiment. Sobolev calibrated experiment serves as a flexible intermediate version of $L_2$ calibration and KO calibration, achieving a trade-off between the two methods, which further allows more potential in practical use.

Table \ref{tab:eg1} summarizes the numerical experiment result. The mean value and and the standard deviation (SD) are computed over 500 random simulations for the three methods. As shown in Table \ref{tab:eg1}, Sobolev calibration is accurate in calibration parameter estimation, which empirically verifies the theoretical results. The standard deviation of Sobolev calibration is relatively small, showing that our method is stable.

It is worth noting that the comparison of different methods is only valid on the specified settings applied to each method. Actually, since our proposed method is a generalization of other two methods, the main goal of the numerical experiment is to emphasize on the flexibility of Sobolev calibration and the merits it brings. Practitioners can select their own norm in the Sobolev calibration to reflect their preference on the function value approximation and function shape approximation.

\end{example}

\begin{table}[H]
\caption{\label{tab:eg1}The mean and standard deviation (SD) of estimated calibration parameters for different methods in Example \ref{dtm-case1}, where $\theta^*_S$ is -0.08.}
\begin{center}

\begin{tabular}{lccccc}
\hline
& \multicolumn{2}{c}{$\sigma^2=0.05$} &  & \multicolumn{2}{c}{$\sigma^2=0.1$} \\ 
 \cline{2-3} \cline{5-6} 
\multicolumn{1}{c}{}  & Mean        & SD        &  & Mean       & SD       \\ \hline
Sobolev   & $ -0.0606$           & $0.0100$         &   & $-0.0605$           & $ 0.0135$   \\
$L_2$  & $-0.4572$           & $0.0263$          &  & $-0.4491$           & $0.0493$       \\
KO  & $-0.0282$           & 0.0178         &  & -0.0280          & 0.0241 \\ \hline
\end{tabular}
\end{center}
\end{table}

\subsection{Additional Tables and Figures in Numerical Studies}\label{app:figures}
\begin{figure}[H]
\begin{center}
\subfigure[Normalized model discrepancy]{
\includegraphics[width=4in]{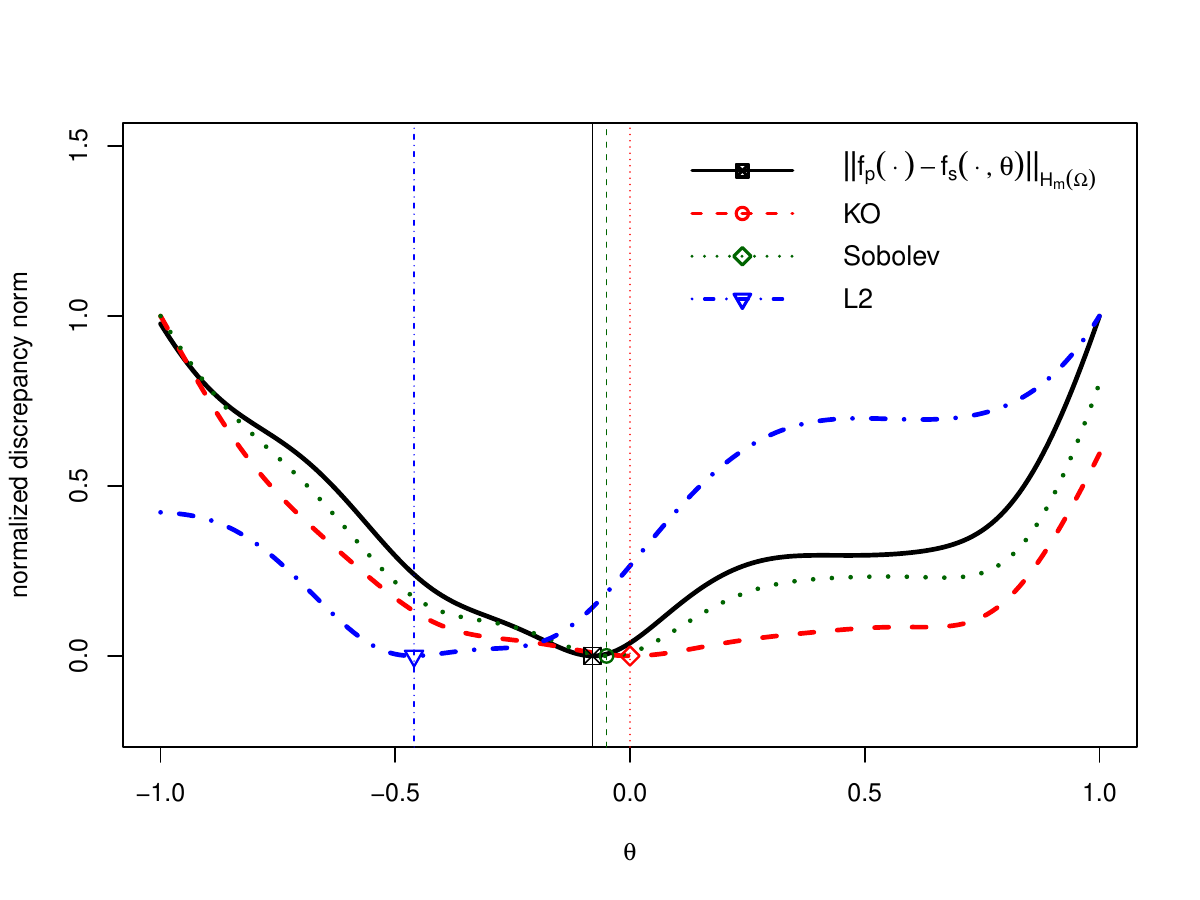}}
\end{center}
\begin{center}
\subfigure[Computer calibrated experiments]{
\includegraphics[width=4in]{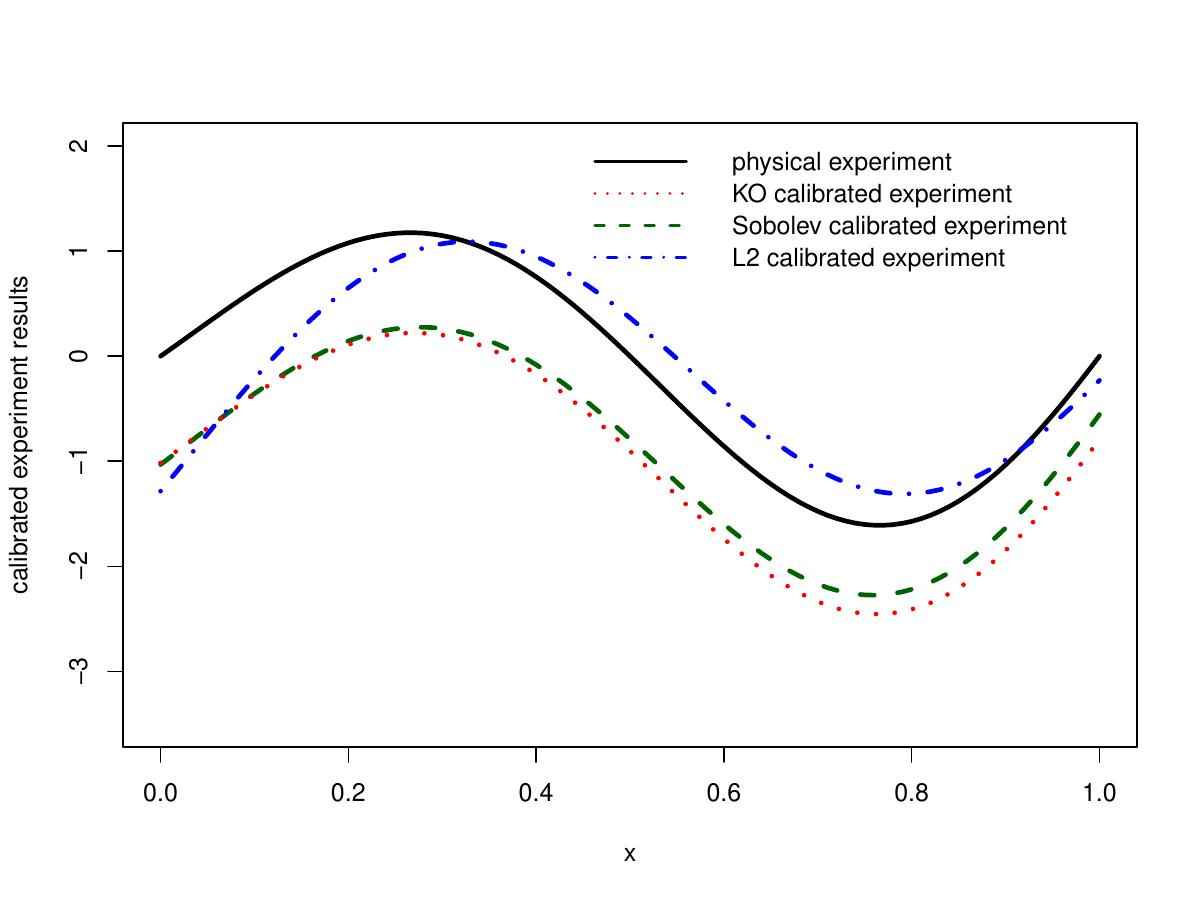}}
\end{center}
\caption{Visualization of Example \ref{dtm-case1} in the simulation studies. Panel (a) presents discrepancy measured by different normalized norms. The symbols (square, circle, triangle, etc.) together with vertical small dotted lines represent the minimizer of different norms, and the solid lines show $\|f_p(\cdot)-f_s(\cdot,\theta)\|_{\cH_m(\Omega)}$ and $\theta^*_{S}$. Panel (b) shows the calibrated computer models (in dotted lines) using the Sobolev calibration, the KO calibration, and the $L_2$ calibration, respectively, together with the physical experiment (in solid line).}
\label{fig:Eg1} 
\end{figure}

\begin{figure}[H]
\begin{center}
\subfigure[Normalized model discrepancy]{
\includegraphics[width=4in]{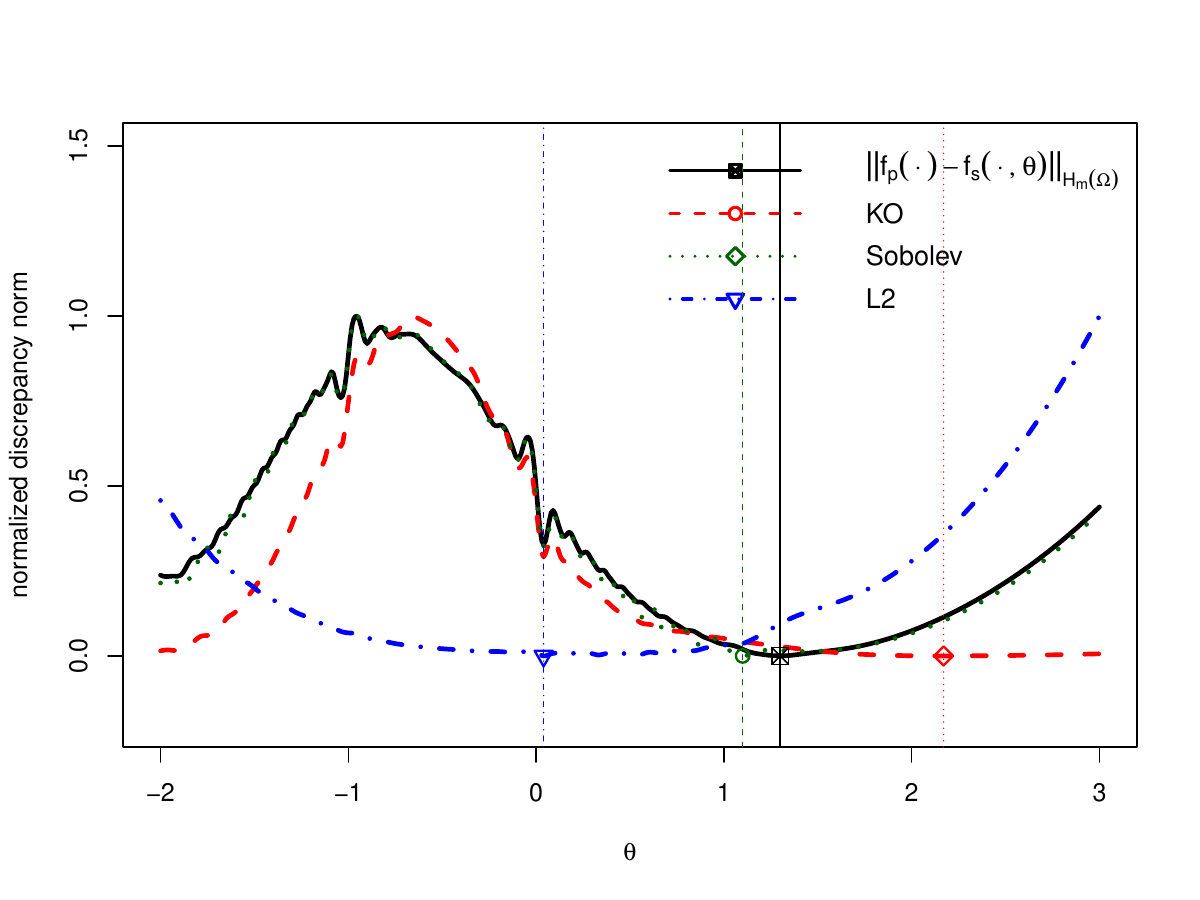}}
\end{center}
\begin{center}
\subfigure[Computer calibrated experiments]{
\includegraphics[width=4in]{pictures/fig_eg1_2.pdf}}
\end{center}
\caption{Visualization of Example \ref{dtm-case2} in the simulation studies. Panel (a) presents discrepancy measured by different normalized norms. The symbols (square, circle, triangle, etc.) together with vertical small dotted lines represent the minimizer of different norms, and the solid lines show $\|f_p(\cdot)-f_s(\cdot,\theta)\|_{\cH_m(\Omega)}$ and $\theta^*_{S}$. Panel (b) shows the calibrated computer models (in dotted lines) using the Sobolev calibration, the KO calibration, and the $L_2$ calibration, respectively, together with the physical experiment (in solid line).}
\label{fig:Eg2} 
\end{figure}

\begin{figure}[H]
\begin{center}
\includegraphics[width=4in]{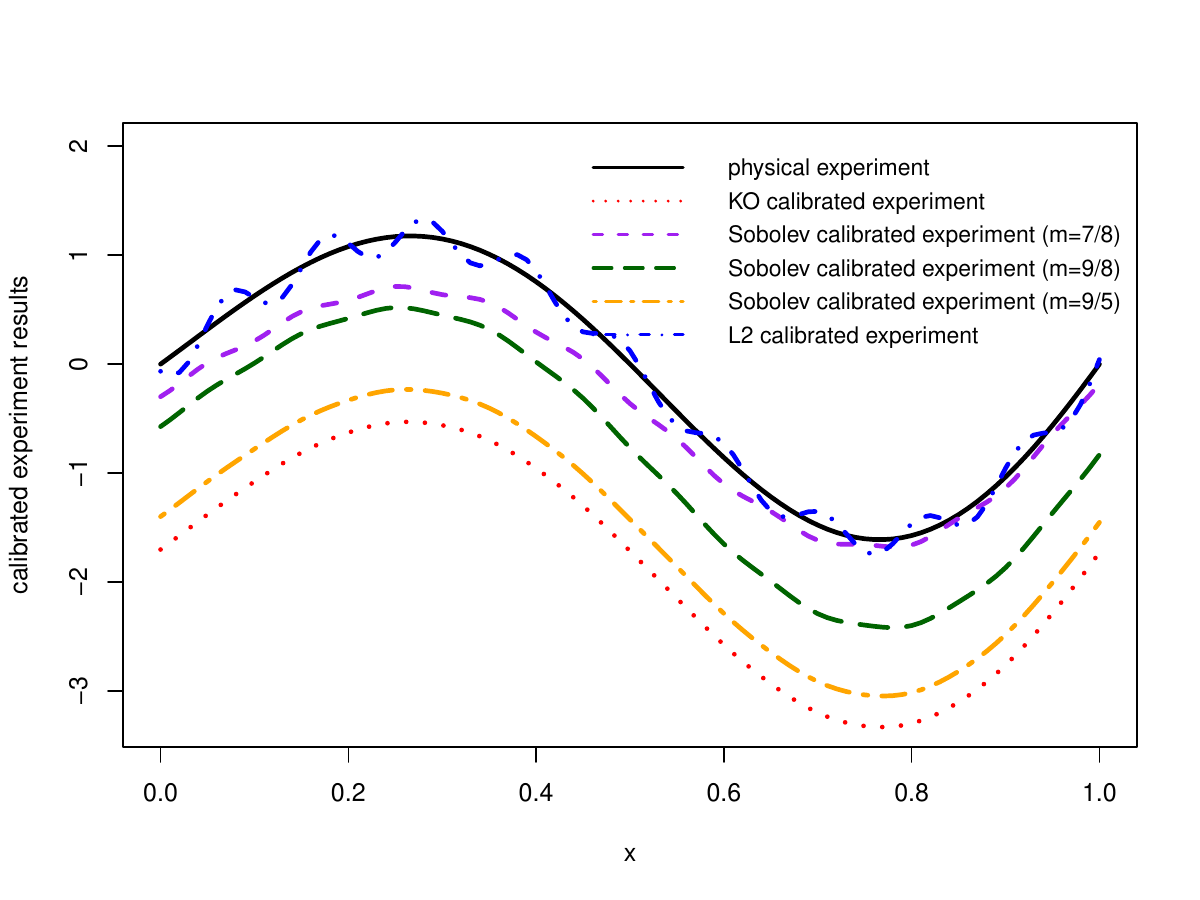}
\end{center}
\caption{Visualization of the calibrated computer models (in dotted lines) using the Sobolev calibration, the KO calibration, and the $L_2$ calibration, respectively, together with the physical experiment (in solid line).}\label{fig:eg2-add}
\end{figure}

\begin{figure}[H]
\begin{center}
\subfigure[Normalized model discrepancy]{
\includegraphics[width=4in]{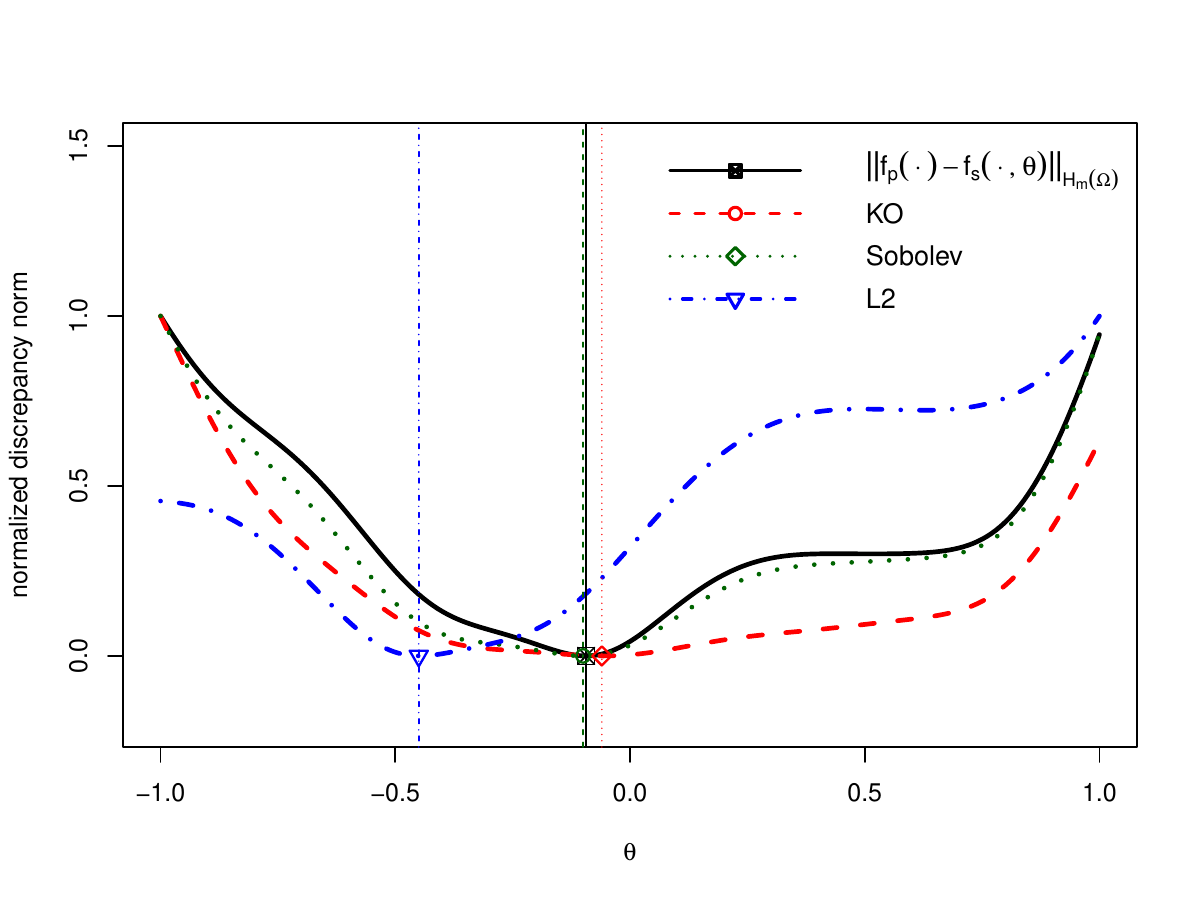}}
\end{center}
\begin{center}
\subfigure[Computer calibrated experiments]{
\includegraphics[width=4in]{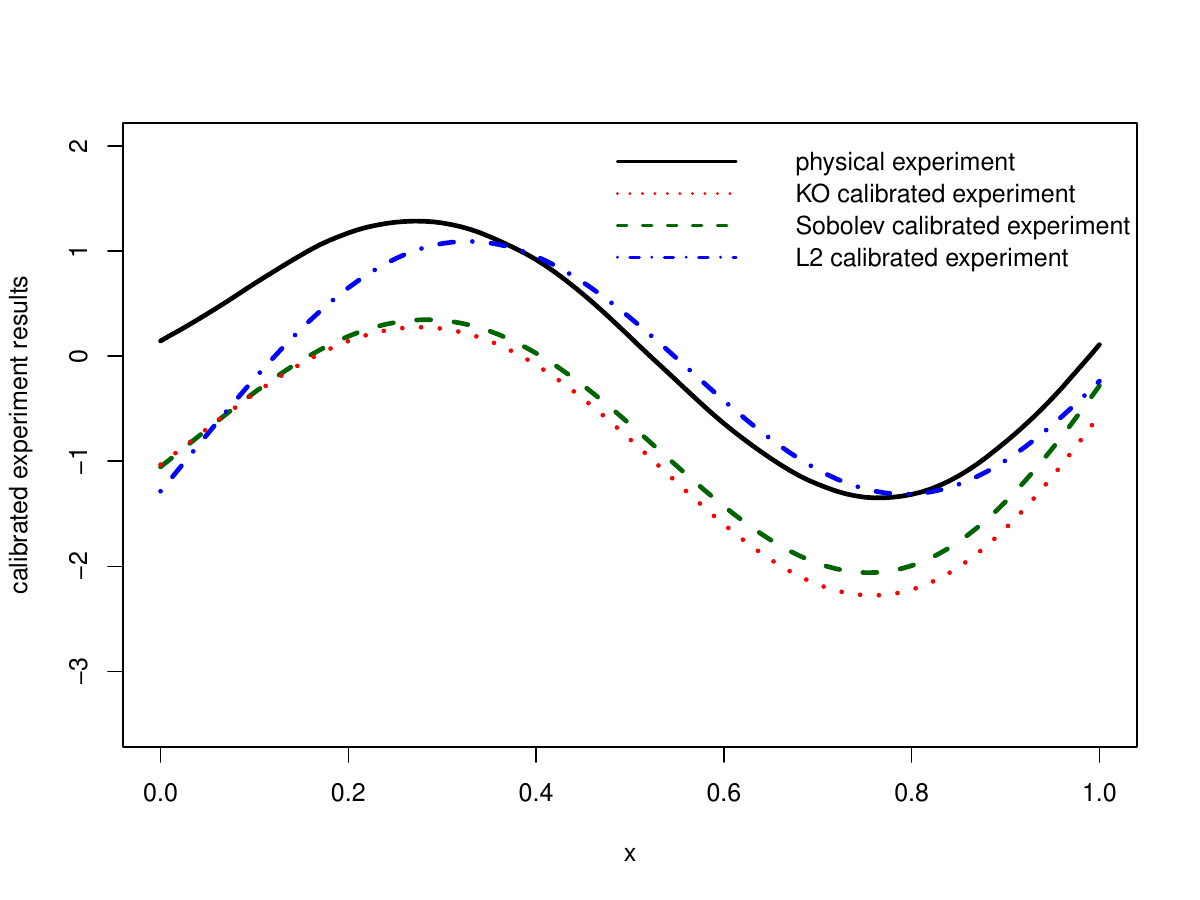}}
\end{center}
\caption{Visualization of Example \ref{gp-case3} in the simulation studies. Panel (a) presents one realization of discrepancy measured by different normalized norms (in dotted lines).
Panel (b) shows one realization of the calibrated computer models (dotted lines) using the Sobolev calibration, the KO calibration, and the $L_2$ calibration, respectively, and the physical experiment (solid line).}
\label{fig:Eg3} 
\end{figure}

\begin{figure}[H]
\begin{center}
\includegraphics[width=4in]{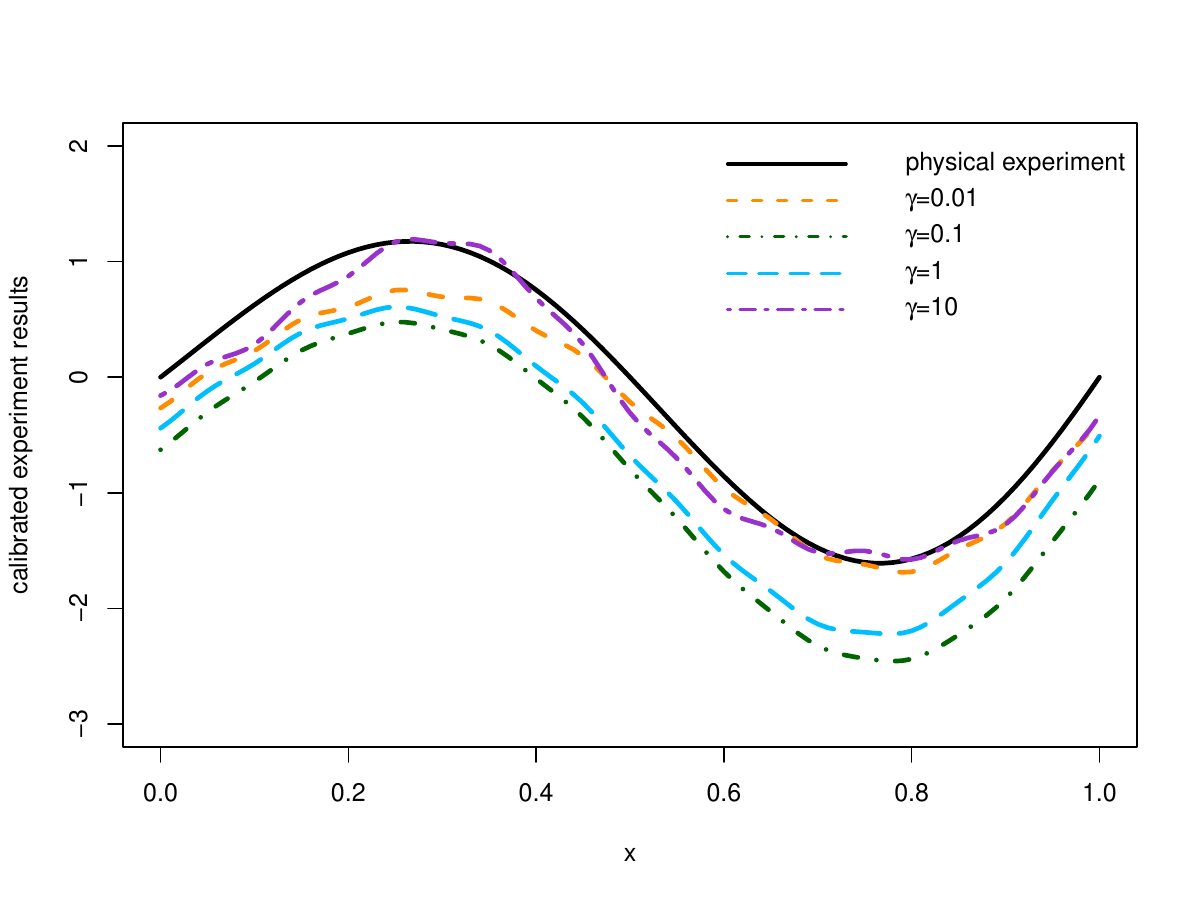}
\end{center}
\caption{{Visualization of experiments in Section \ref{eg:adjustment}. The Sobolev calibration with length-scale parameters $\gamma=0.01,0.1,1,10$ are presented in dotted lines, and the physical experiment is presented in a solid line.}}
\label{fig:scale} 
\end{figure}

\begin{figure}[H]
\begin{center}
    \includegraphics[width=5in]{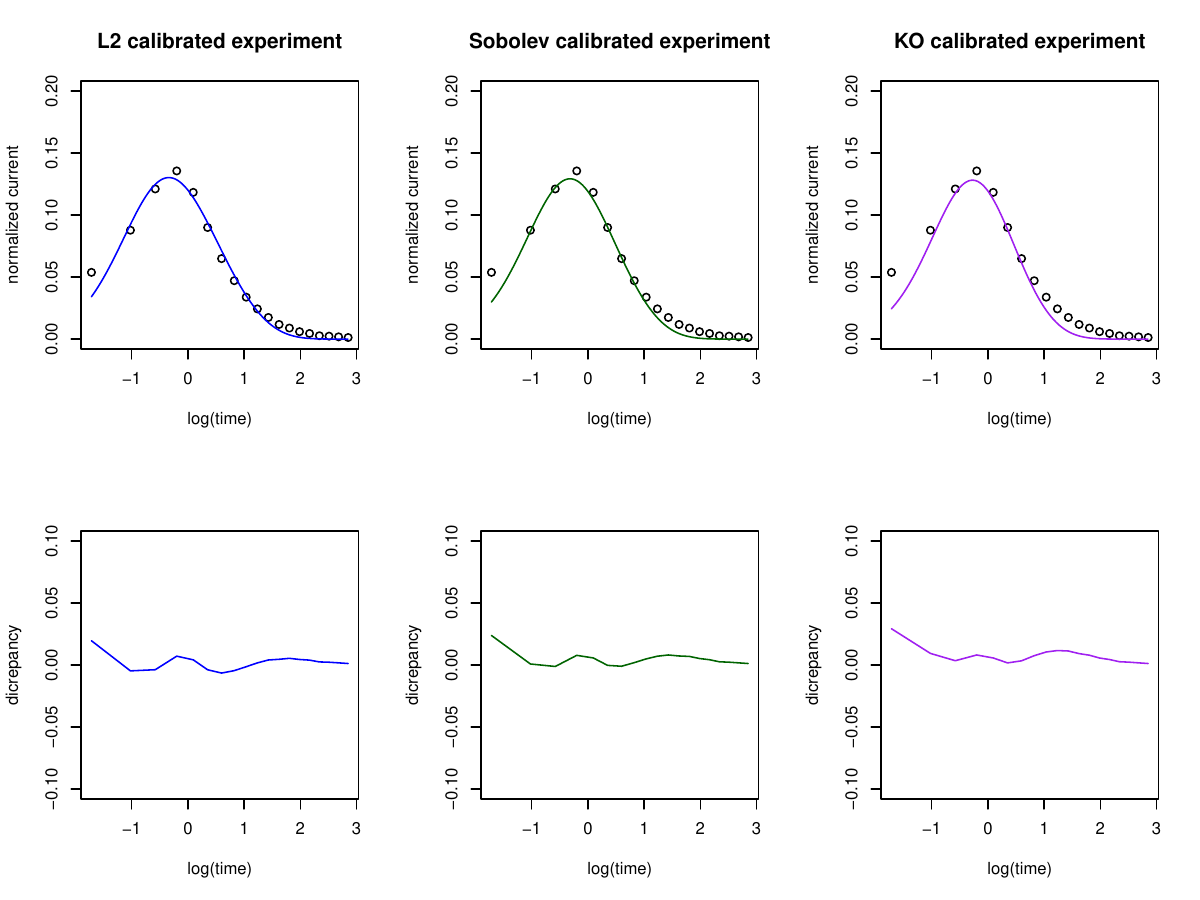}
\end{center}
    \caption{Visualization of computer model calibration for the ion channel example. The top three panels show computer experiment with calibration parameter estimated by $L_2$ calibration, Sobolev calibration, and KO calibration, respectively (in different colors), together with physcial observations (circles); The bottom three panels show discrepancy between physical experiment and computer calibrated experiments with respect to logarithm of time.}
    \label{fig:realeg}
\end{figure}

In this section, we list the figures of simulation studies owing to spatial confined. Figure \ref{fig:Eg1} is for Example \ref{dtm-case1} in Section \ref{app_exp}; Figure \ref{fig:Eg2} and \ref{fig:eg2-add} is for Example \ref{dtm-case2} in Section \ref{subsec:simu}; Figure \ref{fig:Eg3} is for Example \ref{gp-case3} in Section \ref{subsec:simu}. Generally speaking, Figure \ref{fig:Eg1}, Figure \ref{fig:Eg2} and Figure \ref{fig:Eg3} present the normalized discrepancy function and the calibrated experiment on three different settings; Figure \ref{fig:eg2-add} is used for showing the flexibility of Sobolev calibration by specifying function space with different smoothness parameters. Figure \ref{fig:scale} is for the experiments in Section \ref{eg:adjustment}, which presents the calibrated computer experiments with different length-scale parameters. Figure \ref{fig:realeg} is for the Ion Channel Example in Section \ref{subsec:real}.

\end{document}